\documentclass[leqno,12pt]{article}

\usepackage{amssymb}
\usepackage{euscript}
\usepackage{flafter}
\usepackage{pstricks}
\usepackage{graphicx}
\usepackage{caption}
\usepackage{epstopdf} %%Pour inclure des fichiers .eps%%
\usepackage{enumerate}

 \usepackage[font=scriptsize,labelfont=bf]{caption}
\usepackage{pgf,tikz}
\usepackage{mathrsfs}
\usetikzlibrary{arrows}

 \usepackage[latin1]{inputenc}

\usepackage{verbatim}
\usepackage{amsmath}
\usepackage{mathrsfs} %%% Nice caligraphic fonts

\usepackage{ulem}%barrer
\evensidemargin\oddsidemargin
\begin{document}
\baselineskip=18pt
\setcounter{page}{1}

\renewcommand{\theequation}{\thesection.\arabic{equation}}
\newtheorem{theorem}{Theorem}[section]
\newtheorem{lemma}[theorem]{Lemma}
\newtheorem{definition}[theorem]{Definition}
\newtheorem{proposition}[theorem]{Proposition}
\newtheorem{corollary}[theorem]{Corollary}
\newtheorem{remark}[theorem]{Remark}
\newtheorem{fact}[theorem]{Fact}
\newtheorem{problem}[theorem]{Problem}
\newtheorem{conjecture}[theorem]{Conjecture}
\newtheorem{claim}[theorem]{Claim}
\newtheorem{notation}[theorem]{Notation}

%%%%%%%%%%%%%%%%%%%%%%%%%%% Equation numberings
\newcommand{\eqnsection}{
\renewcommand{\theequation}{\thesection.\arabic{equation}}
    \makeatletter
    \csname  @addtoreset\endcsname{equation}{section}
    \makeatother}
\eqnsection
%%%%%%%%%%%%%%%%%%%%%%%%%%%

%%%%%%%%%%%%%% Bbb characters
%%%%%%%%%%%%%% Real numbers
\def\r{{\mathbb R}}
%%%%%%%%%%%%%% Expectation
\def\e{{\mathbb E}}
%%%%%%%%%%%%%% Probability
\def\p{{\mathbb P}}
\def\q{{\mathbb Q}}
%%%%%%%%%%%%%% Integers
\def\z{{\mathbb Z}}
%%%%%%%%%%%%%% Natural numbers
\def\N{{\mathbb N}}
%%%%%%Plane
\def\C{{\mathbb C}}
%%%%%%%%%%%%%%Filtrations
\def\F{{\mathscr F}}
%%%%%%%%%%%%%%Excursions (dÃÂÃÂ©but, fin et orientation)
\def\ggg{{\tt g}}
\def\ddd{{\tt d}}
\def\sgn{{\tt sgn}} 
\def\eee{{\mathfrak e}} %% excursion
\def\G{{\mathbb  G}} %%l'ensemble des dÃÂÃÂ©buts des excursions
%%%%%%%%%%%%%Mesure Ito
\def\nnn{{\tt n}} %%sous P%%%
%%%%%%%%%%%%% travers\'ees
\def\t{{\bf t}} %% d\'ebut
\def\s{{\bf s}} %% fin
%%%%%%%%%%%%%%%%%%%%%%
\def\D{{\mathscr D}} %%%% partition
%%%%%%%%%%%%%%%% Special symbols
%%%%%%%%%%%%%% Exponential
\def\ee{\mathrm{e}}
%%%%%%%%%%%%%% Differentiation
\def\d{\mathrm{d}}
%%%%%%%%%%%%%%
\def\law{{ \buildrel \mbox{\scriptsize \rm(law)} \over = }}
\def\wcv{{\buildrel \mbox{\scriptsize\rm (law)} \over \longrightarrow}}

%%%%Commentares
\newcommand{\note}[1]{{\textcolor{red}{[#1]}}}
%%% 
%%%%%%%%%%%%%% Beginning of the text

%%%%%%%%%%%%%% Beginning of the text

 \vglue30pt

\centerline{\large\bf   Points of infinite multiplicity of planar Brownian motion:} 
 
 \medskip
 \centerline{\large\bf   measures and local times}

\bigskip
\bigskip

 \centerline{by}

\medskip

 \centerline{Elie A\"{i}d\'ekon\let\thefootnote\relax\footnote{\scriptsize LPSM, Sorbonne Universit\'e Paris VI \& IUF, 4 place Jussieu, F-75252 Paris Cedex 05, France, {\tt elie.aidekon@upmc.fr}},   Yueyun Hu\let\thefootnote\relax\footnote{\scriptsize LAGA, Universit\'e Paris XIII, 99 avenue J-B Cl\'ement, F-93430 Villetaneuse, France, {\tt yueyun@math.univ-paris13.fr}}, 
and Zhan Shi\footnote{\scriptsize LPSM, Sorbonne Universit\'e Paris VI, 4 place Jussieu, F-75252 Paris Cedex 05, France, {\tt zhan.shi@upmc.fr}}}

\medskip

 % \centerline{\it Sorbonne Universit\'e Paris VI, Universit\'e Paris XIII, and Sorbonne Universit\'e Paris VI}

%\thankstext{t1}{project partly supported by ANR MALIN (ANR-16-CE93-0003)}
%\thankstext{t2}{partly supported by ANR GRAAL (ANR-14-CE25-0014) and  ANR Liouville (ANR-15-CE40-0013)}
%\thankstext{t3}{partly supported by ANR SWiWS (ANR-17-CE40-0032)}

\bigskip

{\leftskip=2truecm \rightskip=2truecm \baselineskip=15pt \small

 \noindent{\slshape\bfseries Summary.}    It is  well-known (see Dvoretzky,  Erd{\H o}s and  Kakutani \cite{DEK58} and   Le Gall \cite{Legall87})  that a planar Brownian motion $(B_t)_{t\ge 0}$  has points of  infinite multiplicity,  and these    points  form a dense set   on the range.  Our main result is the construction of  a family of  random measures, denoted by  $\{{\mathcal M}_{\infty}^\alpha\}_{0< \alpha<2}$, that are    supported by the set of the points of infinite multiplicity.  We prove that for any  $\alpha \in (0, 2)$,  almost surely the Hausdorff  dimension of   ${\mathcal M}_{\infty}^\alpha$  equals $2-\alpha$, and  ${\mathcal M}_{\infty}^\alpha$ is supported   by  the set of   thick  points      defined in   Bass, Burdzy and Khoshnevisan \cite{BBK}   as well as  by that defined    in    Dembo, Peres, Rosen and Zeitouni \cite{DPRZ01}. 
%     Our construction also reveals    that with probability one, ${\mathcal M}_{\infty}^\alpha(\d x)$-almost everywhere, there exists  a  continuous  nondecreasing additive functional  $({\mathfrak L}_t^x)_{t\ge 0}$, called local times at $x$,  such that  the support of $  \d {\mathfrak L}_t^x$  coincides with   the   level set $\{t: B_t=x\}$. 

  

} %%%%%% End of narrower

\bigskip
\bigskip

\section{Introduction}

\subsection{Problem and results}

It is well-known that a planar Brownian motion has points of infinite multiplicity, and these points form a dense set on the range (Dvoretzky, Erd{\H o}s and Kakutani \cite{DEK58} and   Le Gall \cite{Legall87}); see Le Gall (\cite{Legall90}, p.~204) for comments on the proof.  

In this paper, we restrict our attention to the thick points which constitute a dense  subset of points of infinite multiplicity. In the literature, there are   two  important ways  to define   thick points, one by    Bass, Burdzy and Khoshnevisan \cite{BBK}  through the number  of crossings  and another one by     Dembo, Peres, Rosen and Zeitouni \cite{DPRZ01}  through the occupation times.    

 Let     $(B_t)_{t\ge0}$  be  the standard planar Brownian motion  started at $0$ defined on a complete probability space $(\Omega, \F, (\F_t)_{t\ge0},  \p)$,  where $\F_t:= \sigma\{B_s, s\in [0, \, t]\}$,  $t\ge0$, is the natural filtration of $B$.  For     any $x \in \r^2$ and $r>0$,  denote by   ${\cal C}(x, r)$ (resp: ${\cal B}(x, r)$)    the circle  (resp:  open disc)     centered at $x$ and  with radius $r$. Let     $N_x(r)$   be  the number of crossings  from $x$ to ${\cal C}(x, r)$ by    $(B_t)$  till $T_{{\cal C}(0, 1)}$,  where $T_{{\cal C}(0, 1)}:= \inf\{t\ge0: B_t \in {\cal C}(0,1)\}$ is  the first hitting time of the unit circle ${\cal C}(0,1)$. Fix $\alpha\in (0, 2)$.  Bass, Burdzy and Khoshnevisan \cite{BBK} studied the following set:
%of the points with infinite multiplicity: 
 \begin{equation}\label{Aalpha}
 {\cal A}_\alpha:= \Big\{ x \in {\cal B}(0, 1): \lim_{r\to0^+} \frac{N_x(r)}{\log 1/r}= \alpha\Big\}. 
 \end{equation}

 For any   Borel   measure $\beta$ on $\r^2$,    let   $\mathrm{Dim}(\beta)$  be   the Hausdorff dimension of $\beta$: $$
\mathrm{Dim}(\beta):= \inf\left\{r:\mbox{$\exists$      Borel set $ A  $ such that $\beta(A^c)=0$ and $\mathrm{dim_H}(A )=r$}\right\},$$ 

\noindent  where  $\mathrm{dim_H}(A)$  denotes  the Hausdorff dimension of the set $A$.  The main result in   \cite{BBK}   can be stated   as follows:
  
  \medskip
 
 \noindent {\bf Theorem  A (Bass, Burdzy and Khoshnevisan \cite{BBK})}. {\it  Let $\alpha\in (0, \frac12)$. Almost surely  there exists a  measure $\beta_\alpha$ carried  by ${\cal A}_\alpha$. Moreover, } $$  \mathrm{Dim}(\beta_\alpha)= 2-\alpha, \qquad \mbox{a.s.} $$

The random measure $\beta_\alpha$ in Theorem A gives important information on the set of points to which planar Brownian motion makes ``many visits". For example, it plays a crucial role in Cammarota and M\"{o}rters \cite{CammarotaMortes} who gave a characterization of the gauge functions for the level sets of planar Brownian motion, confirming a conjecture of Taylor~\cite{taylor}.

Another definition of thick points was given by Dembo, Peres, Rosen and Zeitouni \cite{DPRZ01}:   A point     $x\in \r^2$  is called  an {\it $\alpha$-thick point } (called perfectly  thick point in  \cite{DPRZ01}) if 
\begin{equation}\label{perfectlythick}
\lim_{r\to 0^+} \frac1{r^2 (\log r)^2} \, \int_0^{T_{{\cal C}(0, 1)}} 1_{\{ B_s \in {\cal B}(x, r)\}} \d s = \alpha,
 \end{equation}
 \noindent where, as pointed out in \cite{DPRZ01},  $T_{{\cal C}(0, 1)}$ can be replaced by any positive and finite $(\F_t)$-stopping  time. 
% \noindent where $T_{{\cal C}(0, 1)}$ denotes the hitting time of ${\cal C}(0,1)$ and ${\cal B}(x, r)$ the disc of radius $r$ and centered at $x$.  

 \medskip
 The fractal measure of the  $\alpha$-thick points  was studied in depth in \cite{DPRZ01};  in particular the following result   holds: 

\medskip
 \noindent {\bf Theorem B (Dembo, Peres, Rosen and Zeitouni \cite{DPRZ01})} {\it For $0\le  \alpha \le 2$,}  
 $$
\mathrm{dim_H}\left\{ \alpha\mbox{-thick points}\right\}=2-\alpha, \qquad   \mbox{a.s.} 
 $$

To the best of our best knowledge, it is still an open question whether the two ways of defining thick points are equivalent or not. For instance, we have not been able to determine the Hausdorff dimension of ${\cal A}_\alpha$.

The starting point of our study  is to generalize Theorem A to all parameters $\alpha \in (0, 2)$.  Our construction of random measures, different from that in \cite{BBK} where the authors utilized  the   local times on   circles,   relies  on a   change of measures involving the  excursions around   points with infinite multiplicity.  This change of measures is given by \cite{BBK}, Theorem 5.2, where it plays a crucial role.

We consider Brownian motion inside a domain. By a domain in $ \r^2$ we mean an open,   connected  and bounded subset of $ \r^2$.      Given a domain $D$, a boundary point $z\in \partial D$ is said to be {\it nice} if there exists a one-to-one analytic function $f$ from the unit disc such that $f(0)=z$ and the image of the set of points of the unit disc with positive imaginary part is the intersection of $D$ with the image of the unit disc by $f$ .  A domain $D$ is {\it nice} if every boundary point, except perhaps a finite number of them, is {\it nice} (in Lawler \cite{Lawler}, p.48, $\partial D$ is said to be piecewise analytic). Let $\overline D:= D \cup \partial D$. We say that a point $z \in \overline D$ is nice if either $z \in D$ or $z$ is a nice boundary point.

%A domain  $D$   is said to be     {\it nice}     if  its boundary $\partial D$ is a finite union of analytic curves.  A boundary point is {\it nice}  if the boundary is  an analytic curve locally around the point.    Let $\overline D:= D \cup \partial D$ and we say that a point $z \in \overline D$ is nice if either $z \in D$ or $z$ is a nice boundary point. 

Denote by   $\p^{z, z'}_D$   the probability law  of a Brownian excursion inside $D$ from $z$ to $z'$,  whose definition is given  in Section \ref{s:pzd}.  

The main result of this paper  reads  as follows:

\begin{theorem}
\label{t:main} 
Let $D$ be a simply connected nice domain. Let $z$ and $z'$ be distinct nice points of $\overline D$. Fix $\alpha\in (0, \, 2)$. With $\p^{z, z'}_D$-probability one, there exists a random finite measure ${\mathcal M}_{\infty}^\alpha$ carried by  ${\cal A}_\alpha$ as well as by the set of $\alpha$-thick points; moreover,  
 \begin{equation}
  \label{carryingdimension}
     \mathrm{Dim} ({\mathcal M}_{\infty}^\alpha)
     =
     2-\alpha, \qquad \hbox{\rm $\p^{z, z'}_D$-a.s.} 
 \end{equation}
\end{theorem}

We mention that ${\mathcal M}_{\infty}^\alpha$ is uniquely determined by the forthcoming  formula \eqref{spine}, the latter being analogous to \cite{BBK}, Theorem 5.2. Since ${\mathcal M}_\infty^\alpha$  is supported by  the set of $\alpha$-thick points,  \eqref{carryingdimension} yields that  almost surely,  $\mathrm{dim_H}\left\{ \alpha\mbox{-thick points}\right\}\ge 2-\alpha$, giving a new proof of the lower bound in Theorem B. 

We will see that equation \eqref{spine} %Our construction of ${\mathcal M}_\infty^\alpha$ will 
also yields the existence,\footnote{This was a private question by Chris Burdzy.} for ${\mathcal M}_{\infty}^\alpha(\d x)$-almost every $x$, of a  continuous and  non-decreasing additive functional   $({\mathfrak L}_t^x)_{t\ge 0}$, called local times at the point $x$.  
%Fix two distinct nice points   $z, z'$   of $\overline D$.  Let $\p^{z, z'}_D$   be  as before the probability law  of a Brownian excursion from $z$ to $z'$. 
Define
$$
T_x:=\inf\{t>0\,:\, B_t = x \} \, ,
\qquad
x\in \r^2 \, .
$$
  
%\noindent and for any point $x$, $T_{x}=T_{\{x\}}$.  
 
 \begin{theorem}
 \label{t:main2}  
 Let $D$ be a simply connected nice domain. Let $z$ and $z'$ be distinct nice points of $\overline D$. Fix $\alpha\in (0, \, 2)$. With $\p^{z, z'}_D$-probability one, ${\mathcal M}_{\infty}^\alpha(\d x)$-almost everywhere,  there exists  a continuous    Brownian  additive functional  $({\mathfrak L}_t^x)_{t\ge 0}$ such that 
 \begin{equation}
 \label{localtime}
 {\mathfrak L}_t^x= \lim_{r\to 0^+} \frac1{r^2 (\log r)^2} \, \int_0^t 1_{\{|B_s-x |< r\}} \d s, \qquad \forall\,  t\ge 0 \, ;
 \end{equation}
   moreover %$ {\mathfrak L}_{T_{z'}}^x=\alpha$, and  
   the support of $  \d {\mathfrak L}_t^x$ is  identified as   the   level set  $\{t\ge 0: \, B_t=x\}$. 
 \end{theorem}

In Theorem \ref{t:main2}, under $\p^{z, z'}_D$, the process is killed upon first hitting $z'$, i.e.,  at time $T_{z'} :=\inf\{t>0\,:\, B_t = z' \}$. So the local time ${\mathfrak L}_t^x$ is defined for $t\in [0, \, T_{z'}]$. We are going to see that ${\mathfrak L}_{T_{z'}}^x=\alpha$ (with $\p^{z, z'}_D$-probability one, ${\mathcal M}_{\infty}^\alpha(\d x)$-almost everywhere). 

% Our construction of ${\mathcal M}_{\infty}^\alpha$ also allows us to get  % from the capacity of ${\mathcal M}_{\infty}^\alpha$ 
%the lower bound for the Hausdorff dimension of $\alpha$-thick points, see Remark \ref{lowerdprz} for more details. 

\subsection{A brief description of the construction of ${\mathcal M}_{\infty}^\alpha$}
\label{subs:Q}

Let $D$ be a nice domain, and  let  $x \in  D$ and $z, z'$  distinct nice points of $\overline D$, different from $x$. Let as before $T_x:=\inf\{t>0\,:\, B_t = x \}$. We consider  a probability measure $\q_{x,D}^{z, z', \alpha}$  similar to the measure ${\bf Q}^x_\alpha$  introduced in Bass, Burdzy and Khoshnevisan (\cite{BBK}, p.~606): Under $\q_{x,D}^{z, z', \alpha}$,  $(B_t)_{t\ge0}$  is  split into three parts:   \begin{enumerate}
\item Until time $T_{x}$, $B$ is a Brownian motion  starting from  $z$ and conditioned at hitting  $x$ before $\partial D$, whose law  $\p^{z, x}_D$   is  defined in  Notation \ref{pzd}.
\item After $T_{x}$, the trajectory is a concatenation of Brownian loops generated by a Poisson point process $({\mathfrak e}_s)_{s\ge 0}$ with intensity $1_{[0,\alpha]}\d t$ in time and $\nu_{D}(x,x)$ in space, where $\nu_D(x,x)$ denotes the law of Brownian loops in $D$  at $x$ (see \eqref{def-muDxx} for the definition).\footnote{When $\alpha=0$, the second part is reduced to  the single point $\{x\}$. \goodbreak    Let us say a few words  on the concatenation of the loops $({\mathfrak e}_s)$.  Denote by $\zeta({\mathfrak e})$ the lifetime of a loop  ${\mathfrak e}$. Remark   that  $ \sum_{s \le \alpha} \zeta({\mathfrak e}_s)< \infty,  $ a.s., see  \eqref{lawofzeta} for the law of $\zeta$.  For $T_x < t \le T_x + \sum_{s \le \alpha} \zeta({\mathfrak e}_s)$,   let $s\in  (0, \alpha]$ be  such that $t-T_x \in [ \sum_{u < s } \zeta({\mathfrak e}_u), \sum_{u \le s } \zeta({\mathfrak e}_u) )$.   We define $B_t:= {\mathfrak e}_s(t-T_x-\sum_{u < s } \zeta({\mathfrak e}_u)).$      }
\item The last part of the trajectory is a standard  Brownian motion in $D$, started from $x$ and conditioned to hit $z'$ if $z' \in D$ or to exit $D$ at $z'$ if $z' \in \partial D$. The  law $\p^{x, z'}_D$ of this  process   is defined in \eqref{htransform2}.
\end{enumerate}

Consider $D_{1}\subset D$   a nice domain   containing  $x$.  The  baseline  in our construction of ${\mathcal M}_{\infty}^\alpha$,  stated as Corollary \ref{c:radon}, is the absolute continuity of $\q_{x,D}^{z, z', \alpha}$ with respect to $\p^{z,z'}_D$ considered at  $\F^+_{D_1}$, where $\F^+_{D_1}$ denotes  the sigma-algebra generated by the excursions outside $D_1$ together with the order of their appearances (see Notation \ref{notation:tribu2}). Denote by $ M_{D_{1}}(x,\alpha)$ the Radon-Nikodym derivative (up to a renormalization factor $M_{D}(x,\alpha)$). Then $M_{D_{1}}(x,\alpha)$   satisfies a certain restriction property (see Corollary  \ref{r:martingale}). 

Now we   construct ${\mathcal M}_{\infty}^\alpha$ as follows: Let $D$ be a simply connected nice domain and  let     $\D_{n}$ be the connected components of $D$ minus a  grid of mesh size $2^{-n}$.   By using Corollary  \ref{r:martingale}, we may construct a sequence of random measures $\mathcal{M}_{\D_{n}}^\alpha$  (defined in  \eqref{measuredna}). % which is in fact  a  measure-valued martingale. {\tt (xxxx D\'efinition.)} 
The measure ${\mathcal M}_{\infty}^\alpha$ is nothing but the  (nonnegative martingale) limit of $\mathcal{M}_{\D_{n}}^\alpha$ as $n \to \infty$, and the  limit is not trivial and defines a (finite) measure thanks to the uniform integrability of $(\mathcal{M}_{\D_{n}}^\alpha)_{n\ge 1}$  (see Theorem \ref{t:conv-p}). This gives the construction of ${\mathcal M}_{\infty}^\alpha$. % in Theorem \ref{t:main}. 

\subsection{Comparison with previous works} 

In Dembo, Peres, Rosen and Zei\-touni \cite{DPRZ01}, the main interest was focused on the size of the set of the thick points; as a tool, the authors constructed a measure on the thick points as a limit along subsequences of a tight family of measures, for all $\alpha\in (0, \, 2)$. In Bass, Burdzy and Khoshnevisan \cite{BBK}, a measure (in fact, the measure $\beta_\alpha$ recalled in Theorem A) on the thick points was constructed by means of an $L^2$-limit, for all $\alpha\in (0, \, \frac12)$; a formula characterizing the measure $\beta_\alpha$ will be recalled in our paper (see Proposition \ref{p:spine}).

Our work makes a link between \cite{BBK} and \cite{DPRZ01}. Inspired by the conceptual approach of Lyons~\cite{lyons} and Lyons, Pemantle and Peres~\cite{lyons-pemantle-peres} for branching processes, we extend the measure $\beta_\alpha$ in \cite{BBK} to all $\alpha\in (0, \, 2)$, relying on the martingale structure instead of computations of moments. It is to the measure $\beta_\alpha$ in a similar spirit as Lyons~\cite{lyons} to the Biggins martingale convergence theorem for spatial branching processes originally established by Biggins~\cite{biggins-mart-cvg}. 

{\it After this work was completed, a paper of Jego \cite{jego} managed to construct the measure $\beta_\alpha$ for all $\alpha \in (0,2)$ by approximating it via the thick points, as in \cite{BBK}, relying on a truncated second moment method.}

\bigskip

We close this introduction by describing the organization of this paper: 
\begin{itemize}
\item Section \ref{s:pre}: We collect some results on the Brownian  measure  $\p^{z, z'}_D$ on the   paths inside a nice domain $D$ from $z$ to $z'$ and  on the $\sigma$-finite measure  $\nu_D(x,x)$ on the Brownian loops at $x$.
\item Section \ref{s:changeofmeasure}: We characterize the law of the macroscopic excursions at $x$ under $\q_{x,D}^z$ (Proposition \ref{p:crossing}) and get Corollary \ref{c:radon}. We also establish Proposition \ref{p:martingale2} which plays a key role in the proof of Theorem \ref{t:conv-p}. 
\item Section \ref{s:construction}:  A standard argument says that  the proof of  Theorem \ref{t:conv-p}   boils down to the study of ${\mathcal M}_{\infty}^\alpha$ under  $\q_{x,D}^z$ (Proposition \ref{p:Mfini}).   The proof of Proposition \ref{p:Mfini}   is based on a truncation argument  in the computation of moments  under  $\q_{x,D}^z$,  which     also provides the needed capacity  estimates (see \eqref{c_5(N)})  in the proof of  \eqref{carryingdimension}.
\item Section \ref{s:conformal}: The measure  ${\mathcal M}_{\infty}^\alpha$ satisfies a certain conformal invariance, up to a multiplicative factor, which implies its almost sure positivity. 
\item
Section \ref{s:thickpoints}: We  prove that  ${\mathcal M}_{\infty}^\alpha$ is  supported by     ${\cal A}_\alpha$ as well as by  the set of $\alpha$-thick points (Corollary \ref{c:nxr} and Theorem \ref{p:dprz}) and   complete   the proof of Theorem \ref{t:main}.  The proof of Theorem \ref{t:main2} is also  given in Section \ref{s:thickpoints}.
\item Section \ref{s:discussion}: We identify the measure ${\mathcal M}_{\infty}^\alpha$ in the case $\alpha=0$ and discuss its relationship with the intersection local times of independent Brownian motions. 
\end{itemize}

For notational convenience,  we write $\e[X, A ] := \e[X 1_A]$ when $X$ is a random variable and $A$ is an event (and we write $\e[X, A_1, A_2 ]$ for $ \e[X, A ] $  if $A=A_1\cap A_2$).  When $\nu$ is a positive measure, we write $\nu(X):= \int X \d \nu$ and use the similar notation $\nu(X, A)$ for $\int_A X \d \nu$.  Finally, by $f(x) \sim g(x)$ as $ x\to x_0$ we mean that $\lim_{x \to x_0} \frac{f(x)}{g(x)} =1$. 
%%We denote by $(c_i, 1\le i \le 23)$ some positive constants. 

\section{Preliminaries on Brownian excursions}\label{s:pre}
At first we recall some facts, taken from    Lawler (\cite{Lawler}, Chapter 2),  on the Green function and the Poisson kernel. 
Let $D$ be      a  nice domain and consider    $x \in D$.  The harmonic measure $\p^x(B_{T_{\partial D}}\in \bullet)$ is absolutely continuous with respect to the one-dimensional Lebesgue measure, whose density is called the Poisson kernel and is  denoted by $H_D(x, y)$, $y \in \partial D$.  For any nice point $y \in \partial D$, $H_D(\bullet, y)$ is harmonic in $D$.

Let $p_D(t, x, y)$, $x\in D$, $y\in D$ be the density (in $y$) of the transition probabilities of  $B_{t\wedge T_{D^c}}$ under $\p^x$. We define the Green function $G_D$ by 
$$ 
G_D(x, y):= \pi \int_0^\infty p_D(t, x, y) \d t , \qquad x , y  \in D.
$$

\noindent Then
\begin{equation} 
    \label{Gdxy} 
    G_D(x, y) = \log \frac1{|x-y|} +O(1),  \qquad y \to x.
\end{equation}

For any nice point $y \in \partial D$, denote by ${\bf n}_{y}$ the inward normal   at $y$. We have (see Lawler \cite{Lawler}, p.~55) that  
\begin{equation}
    \label{2ndgreen} 
    \lim_{\varepsilon\to0^+} \frac1\varepsilon\, G_D(x, y + \varepsilon {\bf n}_y)= 2 \pi H_D(x, y).
\end{equation}

Following Lawler and Werner \cite{LW04} and Lawler \cite{Lawler}, Section 5.2, we introduce the boundary Poisson kernel: for $y'\neq y$  distinct nice points of $\partial D$,  define 
\begin{equation}\label{hdyy}  H_{D}(y,y') :=  \lim_{\varepsilon\to 0^+}{1\over \varepsilon} H_{D}(y+\varepsilon {\bf n}_{y} ,y') ,  \end{equation}

\noindent For sake of concision, we will use the following notation: if $z,z'$ are distinct points of $D$, we define 
\begin{equation} \label{HDzz'}
H_{D}(z,z') := {1\over 2\pi} G_{D}(z,z').
\end{equation}

\noindent Finally, we define $H_{D}(z,z') := H_{D}(z',z)$ when $z\in \partial D$ and $z'\in D$. The function $H_{D}$ is symmetric and has the following property (see Lawler \cite{Lawler}, Section 5.2): if $D,D'$ are two nice domains,  $z,z'$ nice points of $\overline D$ (recall that  by $y$ nice point of $\overline D$, we mean $y\in D$ or $y$ nice point of $\partial D$), $\Psi:D\to D'$ a conformal transformation %which maps $D$ onto $D'$ 
such that $\Psi(z)$ and $\Psi(z')$ are nice points of $\overline{D'}$, then  
\begin{itemize}
\item if $z$, $z' \in D$, $H_{D'}(\Psi(z),\Psi(z'))=H_{D}(z,z')$,
\item if $z \in D$ and $z'\in \partial D$, $H_{D'}(\Psi(z),\Psi(z')) = |\Psi'(z')|^{-1}\, H_{D}(z,z')$,
\item if $z$, $z' \in \partial D$,  $H_{D'}(\Psi(z),\Psi(z')) = |\Psi'(z)|^{-1}\, |\Psi'(z')|^{-1} \, H_{D}(z,z')$.
\end{itemize}

\noindent Beware that with this definition, $H_{D}$ is not continuous at the boundary of $D$. Indeed, from \eqref{2ndgreen},  we see that for $z\in D$ and $z'$ nice point of $\partial D$,
\begin{equation}\label{HDzz'2}
\lim_{\varepsilon\to 0^+} {1\over \varepsilon} H_{D}(z,z'+\varepsilon {\bf n}_{z'}) =H_{D}(z,z').
\end{equation}

\noindent With this notation, we may give  a unified presentation on  the forthcoming  $h$-transform  \eqref{htransform2}.

%Let ${\mathcal K}$ be the space of   functions $f:   \r_+ \to  \r^2 \cup \{\partial \}$ such that for some $\zeta=\zeta_f \in [0, \infty]$, $f$ is continuous on $[0, \zeta)$ and $f(t)=\partial$ for all $t\ge \zeta$, where $\partial$ denotes some cemetery point.\footnote{When $\zeta_f=0$, $f$ equals identically the cemetery point $\partial$, we denote this function  by ${\mathcal  I}_\partial$.}  

As in Lawler (\cite{Lawler}, Chapter 5), we consider ${\mathcal K}$  the set of all parametrized continuous planar curves $\gamma$  defined on a finite time-interval  $[0, t_\gamma]$ with $t_\gamma \in (0, \infty)$. For any closed subset $A\subset \r^2$, we call $\gamma$ an excursion away from $A$ if
\begin{equation}
    \gamma(0) \in A, 
    \qquad
    \gamma(t_\gamma) \in A,
    \qquad
    \gamma(s) \notin A, \, \forall s\in (0, \, t_\gamma).
    \label{excursion}
\end{equation}

\noindent More generally, for any $0\le s < t< \infty$, we call $(\gamma(u), \, u\in [s, \, t])$ an excursion away from $A$ if $u\mapsto \gamma(u+s)$, $u\in [0, \, t-s]$, is an excursion away from $A$.

The space ${\mathcal K}$ is endowed with the natural filtration of the canonical coordinate process $({\mathfrak e}_t)_{t\ge0}$. We denote by 
$$
T_A:= \inf\{t>0: {\mathfrak e}_t \in A\}\, ,
$$

\noindent the first hitting time of a set $A$ by $({\mathfrak e}_t)$. We will write $(B_t)$ in place of $({\mathfrak e}_t)$ when the underlying measure is a probability measure (such as $\p^z_{D}$ and $\p^{z,z'}_D$).

For any $0\le s < t< \infty$, we denote by $B_{[s,t]}$ the trajectory $u \in [0, t-s]  \mapsto B_{u+s}$, an element in ${\mathcal K}$. We define ${\mathfrak e}_{[s, t]}$ in the same way.

\subsection{Excursion measures inside D}
\label{s:pzd}

Recall that for any Borel set $A\subset \r^2$, 
$$
T_A:=\inf\{t>0\,:\, B_t \in A \} \, .
$$

\noindent Let $D$ be a nice domain. 
For $x\in D$, denote by  $\p^x_{D}$ the probability measure under which $(B_{t},\, t\ge 0)$ is a Brownian motion starting from  $x$ and killed at time $T_{\partial D}$. 
We introduce the probability measure $\p^{z,z'}_D$ (which is the normalized excursion measure denoted by $\mu_{D}^\#(z,z')$ in Lawler \cite{Lawler}, Section 5.2) as follows; the measure is supported by trajectories from $z$ to $z'$ in $D$.

\begin{notation}
 \label{pzd}

 Let $D$ be a nice domain. Let $z$ and $z'$ be distinct nice points of $\overline{D}$. We define the probability measure $\p^{z,z'}_D$ as follows.
 
{\rm (i)} If $z\in D$ and $z'\in \overline D$, we let $\p^{z,z'}_D$ be the law of the Brownian motion starting at $z$,  conditioned to  hit $z'$ before $\partial D$ if $z' \in D$,  and conditioned to exit $D$ at $z'$ if $z' \in \partial D$.  It is given by the  $h$-transform of the Brownian motion: for any  $0< r < | z- z'|$, 
\begin{equation}
\label{htransform2}
 \frac{\d \p^{z,z'}_D}{\d \p^z_{D}} \Big|_{\F_{T_{{\cal C}(z', r)}}}= \frac{H_D(B_{T_{{\cal C}(z', r)}},z')}{H_D(z,z')}1_{\{ T_{{\cal C}(z',r)} < T_{\partial D}\}}  .   
\end{equation}

%\begin{comment}
%\item If $z\in D$ and $z'\in D$, with $z\neq z'$, we let  $\p^{z,z'}_D$ be the law of the Brownian motion starting at $z$ and conditioned at hitting $z'$ before hitting $\partial D$. It can be defined as the $h$-transform of a Brownian motion: For any $r >0$ such that ${\cal B}(z', r) \subset D$,   \begin{equation}\label{htransform}
% \frac{\d \p^{z,z'}_D}{\d \p^z_{D}} \Big|_{\F_{T_{{\cal C}(z', r)}}}= \frac{G_D(z', B_{T_{{\cal C}(z', r)}})}{G_D(z', z)} 1_{\{T_{{\cal C}(z',r)}< T_{\partial D}\}}  .   \end{equation}
%
%
%\item If $z\in D$ and $z'$ a nice boundary point of $\partial D$, we let  $\p^{z,z'}_D$ be the Brownian motion starting at $z$ conditioned to exit $D$ at $z'$. It is given by the  $h$-transform of the Brownian motion:
%\begin{equation}
%\label{htransformH}
% \frac{\d \p^{z,z'}_D}{\d \p^z_{D}} \Big|_{\F_{T_{{\cal C}(z', r)}}}= \frac{H_D(B_{T_{{\cal C}(z', r)}},z')}{H_D(z,z')}1_{\{ T_{{\cal C}(z',r)} < T_{\partial D}\}} , \qquad \mbox{ for all small $r>0$} .   
%\end{equation}
%\end{comment}

{\rm (ii)} If $z$, $z' \in \partial D$, we define $\p^{z,z'}_D$ as the limit of $\p^{y,z'}_D$ as $y\to z$ and $y\in D$ in the sense of \eqref{conv-functionalzz'} below. It is the excursion measure at $z$  conditioned to exit $D$ at $z'$. %When $z=z'$,  we define     $\p^{z', z'}_D$   to be  the Dirac measure at ${\mathcal  I}_\partial$. 

{\rm (iii)} If $z\in \partial D$ and $z'\in D$, we define $\p^{z,z'}_D$ as the limit of $\p^{y,z'}_D$ as $y\to z$ and $y\in D$ in the sense of \eqref{conv-functionalzz'}. It is the excursion measure at $z$ conditioned to hit the interior point $z'\in D$.

\end{notation}

For any $z \in \partial D$ nice point, the limit of $\p^{y,z'}_D$ when $y\to z$ is understood in the following sense: for any  $0< r < | z- z'|$, any ${\bf A}\in \sigma \{B_{T_{D\cap {\cal C}(z,r)+t}}, \, t\ge 0 \}$,
 \begin{equation}\label{conv-functionalzz'}
\p^{z,z'}_D\left( {\bf A} \right)
=  \lim_{\varepsilon\to 0^+  }\p^{z + \varepsilon {\bf n}_z,z'}_D\left( {\bf A} \right).
\end{equation}
 
\noindent See Appendix \ref{A:s1} for a justification of \eqref{conv-functionalzz'}.

 There is a  time-reversal  relationship between $\p^{z,z'}_D$ and $\p^{z', z}_D$ (see Lawler  \cite{Lawler}, Section 5.2)   which can also be checked   by applying the  general  theory on the time-reversal of Markov processes  at cooptional times (see Revuz and Yor \cite{RY}, Theorem VII.4.5):  for any $z, z'$ distinct  nice points  of  $\overline D$, 
\begin{equation}
    \label{timereversal}
    \big((B_{T_{z'}-t}, 0\le t \le T_{z'}) \, \mbox{ under } \p^{z,z'}_D \big) 
    \, \law\, 
    \big( (B_t, 0\le t\le T_{z})  \, \mbox{ under } \p^{z',z}_{D}\big).
\end{equation}

If $D,D'$ are two nice domains,  $z,z'$ nice points of $\overline D$, $\Psi:D\to D'$ a conformal transformation %which maps $D$ onto $D'$ 
such that $\Psi(z)$ and $\Psi(z')$ are nice points of $\overline{D'}$, then the image measure of $\p_D^{z,z'}$ by $\Psi$ is $\p_{D'}^{\Psi(z),\Psi(z')}$.

\medskip

\noindent   We introduce the following notation: for $D_{1}\subset D$   two nice  domains, define   \begin{eqnarray}
 \mathcal J(D, D_1)
 &:=&
 \Big\{z  \in \partial D: \mbox{nice point such that $z\in \partial D_1$ } 
 \nonumber \\
 &&  \mbox{ and  $\exists \, r>0$ such that $D\cap {\cal B}(z, r) = D_1 \cap {\cal B}(z, r)$}
 \Big\}. \label{JDD1}
 \end{eqnarray}

 \begin{lemma} \label{l:hitting-h}
Let $D_{1}\subset D$ be two nice  domains.

(i) For $z'$ nice point of  $\overline D$ and  $z\in D_{1}$   with $z\neq z'$,  for any nice point $y\in \partial D_{1}\backslash \partial D$,
$$
\p^{z,z'}_D(T_{z'}>T_{\partial D_{1}}, B_{T_{\partial D_{1}}} \in \d y)= {H_{D}(y,z') \over H_{D}(z,z')} H_{D_{1}}(z,y) \d y.
$$

(ii) For $z,z'$ distinct points of $D_{1}\cup {\mathcal J}(D,D_{1})$,  we have $$
\p^{z,z'}_D(T_{z'} \le T_{\partial D_{1}})= {H_{D_{1}}(z,z') \over H_{D}(z,z')}.
$$
\end{lemma} 

\noindent {\it Proof}.  (i)     For any nonnegative   measurable  function  $F$ on ${\mathcal K}$, we deduce from the monotone convergence theorem that  
$$
\e^{z,z'}_D(F(B_{[0,T_{\partial D_{1}}]}),T_{z'}>T_{\partial D_{1}}) = \lim_{r\to 0^+} \e^{z,z'}_D(F(B_{[0,T_{\partial D_{1}}]}),T_{{\cal C}(z',r)}>T_{\partial D_{1}}).
$$

\noindent Applying  the change of measures in  \eqref{htransform2} to $T_{{\cal C}(z',r)} \wedge T_{\partial D_{1}}$ in place  of $T_{{\cal C}(z',r)}$,  we have 
\begin{eqnarray*}
 && \e^{z,z'}_D(F(B_{[0,T_{\partial D_{1}}]}),T_{{\cal C}(z',r)}>T_{\partial D_{1}})
 \\
 & =&
 {1\over H_{D}(z,z')} \e^z_{D}\left[  H_{D}(B_{T_{\partial D_{1}}},z') F(B_{[0,T_{\partial D_{1}}]}),  T_{\partial D} \wedge   T_{{\cal C}(z',r)}>T_{\partial D_{1}}\right].
\end{eqnarray*}

\noindent Hence, taking the limit $r\to 0^+$,
\begin{eqnarray}\label{ezxdF}
&& 
	\e^{z,z'}_D\big(F(B_{[0,T_{\partial D_{1}}]}),T_{z'}>T_{\partial D_{1}}\big) \\
&=& 
	{1\over H_{D}(z,z')} \e^z_{D}\left[  H_{D}(B_{T_{\partial D_{1}}},z')F(B_{[0,T_{\partial D_{1}}]}),  T_{\partial D}  >T_{\partial D_{1}}\right], \nonumber
\end{eqnarray}

\noindent  which readily gives (i).

(ii)  First we suppose that %$z\neq z'\in D_{1}$
$z$ and $z'$ are distinct points of $D_1$. By  \eqref{ezxdF},    
\begin{eqnarray*}
	\p^{z,z'}_D(T_{z'}>T_{\partial D_{1}} )  
&=&  
	\frac{\e_D^z\left[  H_{D}(B_{T_{\partial D_{1}}},z'),  T_{\partial D}  >T_{\partial D_{1}}\right]}{H_{D}(z,z')} \\
&=&
	\frac{\e^z\left[  G_{D}(B_{T_{\partial D_{1}}},z')\right]}{G_{D}(z,z')}  .
\end{eqnarray*}

Note that  the function $g$ defined by  $g(y):=G_{D}(y,z')-G_{D_{1}}(y,z')$, $y \in D_1\backslash \{z'\}$,  is harmonic  and bounded. Then $g$ can be continuously extended to $\{z'\}$ so that  $g$ is harmonic on $D_1$.\footnote{We  give here a probabilistic argument on  this known continuous extension.  Let  $z \in D_1 \backslash \{z'\}$   and let $\varepsilon>0$ be small such that ${\cal B}(z, \varepsilon ) \subset D_1$. By the optional stopping theorem for the bounded martingale $g(B_{t\wedge T_{{\cal C}(z', \varepsilon)}})$  at $T_{\partial D_1}$ and by letting $\varepsilon\to 0^+$, we get that $g(z)= \e^z (g(B_{T_{\partial D_1}}))= \int _{\partial D_1} H_{D_1}(z, y) g(y) \d y $ which can be continuously extended to $z'$ and the extension is harmonic on $D_1$.} By the optional stopping theorem for the bounded martingale $g(B_t)$, we get   that  $$\e^z\Big[  G_{D}(B_{T_{\partial D_{1}}},z')\Big]=G_{D}(z,z')-G_{D_{1}}(z,z'),$$

\noindent proving (ii) in this  case. Making $z$ go to $\partial D\cap \partial D_{1}$, using \eqref{2ndgreen} and \eqref{conv-functionalzz'} yields the equality for $z'\in D_{1}$ and $z\in \mathcal J(D,D_{1})$. Using the time-reversal \eqref{timereversal} yields the equality for $z\in D_{1}$ and $z'\in \mathcal J(D,D_{1})$. Making $z$ go to $\partial D \cap \partial D_{1}$ and using \eqref{hdyy} and \eqref{conv-functionalzz'} gives the equality for $z,z' \in  \mathcal J(D,D_{1})$.   \hfill$\Box$

\subsection{Brownian loops}

 Let $x\in D$. We   consider an infinite measure $\nu_D(x,x)$  on the set of Brownian loops  in $D$     that start and end at $x$.   Let us briefly recall the definition     in Lawler (\cite{Lawler}, Chapter 5).     For each $t>0$, the measure $\mu_D(x, x; t)$  of loops  in $D$ of length $t$ is such that  for any    bounded measurable function $\Phi$  on ${\mathcal K}$:  $$
\mu_D(x, x; t)\left[\Phi({\mathfrak e})\right]
=
\lim_{\varepsilon\to0^+}  \frac1{\pi \varepsilon^2} \e^x\left[\Phi(B_{[0, t]}), t < T_{\partial D}, |B_t-x| < \varepsilon\right].$$

\noindent We define
the $\sigma$-finite measure  $\nu_D(x,x)$ by   \begin{equation}\label{def-muDxx}  \nu_D(x,x):= \pi \int_0^\infty \mu_D(x, x; t) \d t .\end{equation}

 With the notation  $\mu_D(x, x)$ of the loop measure   in Lawler (\cite{Lawler}, Chapter 5), we have $\nu_D(x, x) = \pi \mu_D(x, x)$.

Furthermore,    $\nu_D(x, x)$   can also be   viewed as   the excursion measure at $x$ of the Brownian motion conditioned to hit $x$ before hitting ${\partial D}$.  We claim that  
\begin{equation}\label{nuxx}
\nu_{D}(x,x) = \lim_{z \to x}  \log \Big({1 \over |z-x|}\Big) \, \p^{z,x}_{D},
\end{equation}

\noindent 
in the sense  that for  any  $r>0$ such that ${\cal B}(x, r) \subset D$,  for any ${\bf A} \in \sigma \{{\mathfrak e}_{{\cal C}(x,r)+t},   t\ge 0 \} $, \begin{equation}\label{conv-functional4}
\nu_{D}(x, x)\left( {\bf A} \cap \{ T_{{\cal C}(x,r)} < T_{x}\}\right)
= \lim_{z \to x}  \log \Big({1 \over |z-x|}\Big) \, \p^{z,x}_{D}\left({\bf A} \cap \{ T_{{\cal C}(x,r)}< T_{x}\}\right).
\end{equation}

\noindent See Appendix  \ref{A:s2} for a justification of \eqref{conv-functional4}. Equation \eqref{nuxx} implies that the measure $\nu$ is conformally invariant: if $D$ and $D'$ are two nice domains, $\Psi:D\to D'$ a conformal transformation and $x\in D$, the image by $\Psi$ of a loop under $\nu_D(x,x)$  is ``distributed'' as $\nu_{D'}(\Psi(x),\Psi(x))$.

We summarize  some quantitative results  on $\nu_D(x,x)$ in the following Lemma:   

\begin{lemma}\label{l:ito-muxx}
Let $D_{1}\subset D$ be two nice domains and $x\in    D_{1}$. 
\begin{enumerate}[(i)]
\item For any nice point $y \in \partial D_{1}$,
$$
\nu_{D}(x,x)(T_{\partial D_{1}}<T_{x}, {\mathfrak e}_{T_{\partial D_{1}}}\in \d y)=  G_{D}(x,y)H_{D_{1}}(x,y)\d y.
$$
Moreover  for any nonnegative measurable function $F$ on ${\mathcal K}$,
$$
\nu_{D}(x,x) (F({\mathfrak e}_{[0, T_{\partial D_{1}}]}) , T_{\partial D_{1}} <T_{x}) 
=
\e^x \left[G_D(x, B_{T_{\partial D_1}})\, F(B_{[0, T_{\partial D_1}]})\right].
$$
\item   We have 
\begin{eqnarray}\label{CDD1}
C_{D,D_{1}}(x)&:=& \nu_{D}(x,x)(T_{\partial D_{1}}<T_{x})
\nonumber
\\
&=&
 \int_{\partial D_{1}} G_{D}(x,y)H_{D_{1}}(x,y)\d y 
 \nonumber
 \\
 &=&C_{D_{1}}(x)-C_{D}(x), 
\end{eqnarray}
where $C_{S}(x):= - \int_{\partial S} \log(|x-y|)H_{S}(x,y)\d y $ for any nice domain $S$. Moreover, $C_{D,D_{1}}(x)= \log (R_{x,D})-\log (R_{x,D_{1}})$ in the case that $D_{1},D$ are simply connected, where  $R_{x,D}$ (resp. $R_{x,D_{1}}$) is the conformal radius of $D$ (resp. $D_{1}$) seen from $x$.
\item For any nonnegative measurable function $F$ on ${\mathcal K}$,
\begin{eqnarray*}
&&	\nu_{D}(x,x) (F({\mathfrak e}_{[T_{\partial D_{1}}, T_{x}]}) , T_{\partial D_{1}} <T_{x}) \\
&=&  \int_{\partial D_{1}}\e^{y,x}_{D}[F(B_{[0,T_{x}]})]G_{D}(x,y)H_{D_{1}}(x,y)\d y.
\end{eqnarray*}
Consequently, $$\nu_{D}(x,x) \left(F({\mathfrak e}_{[T_{\partial D_{1}}, T_{x}]}) \big| T_{\partial D_{1}} <T_{x},  {\mathfrak e}_{T_{\partial D_1}}=y \right) = \e^{y,x}_{D} \left[  F(B_{[0,T_{x}]})\right].$$
\end{enumerate}
\end{lemma}

{\noindent\it Proof}.   (i) By \eqref{conv-functional4}  and    \eqref{ezxdF},   for any $r>0$ such that ${\cal B}(x, r) \subset D_1$,  \begin{eqnarray*}
&&
\nu_{D}(x,x) (F({\mathfrak e}_{[T_{{\cal C}(x, r)}, T_{\partial D_{1}}]}) , T_{\partial D_{1}} <T_{x}) 
\\
&=&
\lim_{z \to x}  \frac{\log (1 / |z-x|)}{G_{D}(x,z)} \e^z\left[  G_{D}(x,B_{T_{\partial D_{1}}})F(B_{[T_{{\cal C}(x, r)},T_{\partial D_{1}}]})\right]
\\
&=&
\e^x\left[  G_{D}(x,B_{T_{\partial D_{1}}})F(B_{[T_{{\cal C}(x, r)},T_{\partial D_{1}}]})\right],
\end{eqnarray*}

%Lemma \ref{l:hitting-h}, we have 
%$$ \nu_D(x,x)(T_{\partial D_{1}}<T_{x},B_{T_{\partial D_{1}}}\in dy) = \lim_{z \to x}  \log \left({1 \over |z-x|}\right) {  G_{D}(x,y)H_{D_{1}}(x,y)dy \over G_{D}(x,z)} $$

\noindent   by \eqref{Gdxy}.   Under $\p^x$, $B_{T_{\partial D_{1}}}$ is distributed as $H_{D_1}(x, y ) \d y$, which implies the first equality in (i) by taking $F(B_{[T_{{\cal C}(x, r)},T_{\partial D_{1}}]})$ a function of $B_{T_{\partial D_{1}}}$. By considering $F$ a continuous and bounded function on ${\mathcal K}$, and then letting $r\to 0^+$, we get the    second equality in (i).    
%%Or by the monotone class theorem, it is enough to consider $F$ of form $F(B_{\varepsilon+\cdot})$ for any $\varepsilon>0$

(ii) 
The first equality in \eqref{CDD1} readily follows from (i). Let $g(y):=G_{D}(x,y)-G_{D_{1}}(x,y)$ for $y \in D_1\backslash\{x\}$ [$x\in D_1$ being fixed].  We have already observed  that  $g$ can be continuously extended to $x$. Then   $g$ is harmonic  and bounded on $D_{1}$,  and  by the mean property, 
\begin{equation}\label{gx}
g(x)=\int_{\partial D_{1}} g(y)H_{D_{1}}(x,y)\d y=\int_{\partial D_{1}} G_{D}(x,y)H_{D_{1}}(x,y)\d y. 
\end{equation}
The same argument  shows that, for $S=D$ or $S=D_{1}$,
 $$
 \int_{\partial S} \log(|x-y|)H_{S}(x,y)dy = \lim_{y\to x}\Big( \log(|x-y|) + G_{S}(x,y)\Big),
 $$

\noindent which  implies that $C_{D,D_{1}}(x)=C_{D_{1}}(x)-C_{D}(x)$, as stated in \eqref{CDD1}. 

Now we suppose that   $D$  is  simply connected.   It suffices  to  check that 
 $$
C_{D}(x)
=
- \log (R_{x,D}) .
$$ 

  Let $\Phi_{D}$   be a conformal map from $D$ (resp. $D_{1}$) to ${\cal B}(0,1)$ which sends $x$ to $0$ and satisfies that $\Phi_D'(x) >0$.    By the conformal invariance, $G_D(x, y)= \log \frac1{|\Phi_D(y)|}$.  It follows that $$
  \lim_{y\to x}\Big( \log(|x-y|) + G_{D}(x,y)\Big)
  =
  \log \frac1{|\Phi'_D(x)|} 
  =
  \log R_{x, D},$$
  
%%  Notice that  $G_D(x, y)= \lim_{\varepsilon\to0} (\log \frac1\varepsilon) \,\p^x \big( T_{{\cal C}(y, \varepsilon)} < T_{\partial D}\big)$ (see Lawler \cite{Lawler}, p.54). The conformal mapping $\Phi_D$ gives that $\p^x \big( T_{{\cal C}(y, \varepsilon)} < T_{\partial D}\big)= \p^0\big( T_{\partial \Phi_D({\cal B}(y, \varepsilon))} < T_{{\cal C}(0,1)}\big)$. Observe that as $\varepsilon \to 0^+$,  $ \Phi_D({\cal B}(y, \varepsilon) )   \approx  {\cal B}(\Phi_D(y), \varepsilon \, |\Phi_D'(y)|)$, \footnote{By $``\approx"$ we mean that for any small $\delta>0$,   ${\cal B}(\Phi_D(y), \varepsilon \, (|\Phi_D'(y)|- \delta)) \subset   \Phi_D({\cal B}(y, \varepsilon) ) \subset  {\cal B}(\Phi_D(y), \varepsilon \, (|\Phi_D'(y)|+\delta))$.}  we apply Lemma 2.1 in \cite{BBK}  and get that $G_D(x, y)= \log \frac1{|\Phi_D(y)|}$.  The same holds for $G_{D_1}(x, y)$. 

   \noindent   by definition of the conformal radius.

 (iii)  It  comes from a straightforward application of the  strong Markov property together with the first equality in  (i). \hfill $\Box$

\section{A change of measures}\label{s:changeofmeasure}

Let $\alpha\ge 0$ be fixed. Let $D$ be a nice domain. We consider    $D_{1}\subset D$   a nice domain   containing  $x$.  We define for any distinct nice points $y,z\in  \overline D_{1} $,  different from $x$, 
\begin{equation}
    \label{def-xi}
    \xi_{D_{1}}(x,y,z)
    := 
    \frac{2\pi H_{D_{1}}(x,y)H_{D_{1}}(x,z)}{H_{D_{1}}(y,z)} .
\end{equation}

\noindent [Taking $D_1=D$, we get the definition of $ \xi_{D}(x,y,z)$ for any $x\in D$ and distinct nice points $y, z \in  \overline D$,  different from $x$].  %The quantities  $\xi_{D}$ and $\xi_{D_{1}} $ will appear in  the formula of  the change of measures, see the forthcoming Corollary \ref{c:radon}. 

Fix two nice points $z \neq z'$  of $\overline D$.  Let $x\in D$ distinct of $z$ and $z'$. Let us take a close look at the probability measure $\q^{z,z',\alpha}_{x, D}$ defined in the Introduction. Under $\q^{z,z',\alpha}_{x, D}$, the Brownian motion $B$ starts from $z$ and is conditioned to hit $x$. After $T_x$, the trajectory is a concatenation of Brownian loops generated by a Poisson point process $({\mathfrak e}_s)_{s\ge 0}$ with intensity $1_{[0,\alpha]}\d t \times \nu_{D}(x,x)$. The last part of trajectory   is a Brownian motion in $D$ starting from $x$ and conditioned to hit $z'$.

The purpose of this section is to study the absolute continuity of $\q_{x,D}^{z, z', \alpha}$ with respect to $\p^{z,z'}_D$, both restricted to the sigma-algebra generated by the excursions outside a domain $D_1$ containing $x$.

As a first step, we   look at   events of crossings. We introduce some notation. Let ${\mathcal K}$ denote as before the set of all parametrized continuous planar curves $\gamma$  defined on a finite time-interval $[0, t_\gamma]$ with $t_\gamma \in (0, \, \infty)$.

\begin{notation}
 \label{ti}

 Let $D_{2}\subset D_{1}\subset D$ be nice domains with $d(D_{2},\partial D_{1})>0$. Let $\gamma \in {\mathcal K}$.%, and denote  $b:= \gamma(t_\gamma)$.  

{\rm (i)} Define $\s_0 := 0$ and for all $i\ge 0$ (with $\inf \varnothing := \infty$),  
\begin{eqnarray*}
    \t_{i+1} 
 &:=& \inf\{t>\s_{i}\,:\, \gamma(t) \in \partial D_{2}\} ,
    \\
    \s_{i+1} 
 &:=& \inf\{t>\t_{i+1}\,:\,  \gamma(t) \in \partial D_{1} \cup \{\gamma(t_\gamma)\} \}.
\end{eqnarray*}

{\rm (ii)} Let $x\in D_{2}$. For any integer $\ell\ge 1$, denote, if it exists, by ${\bf i}_{\ell}$ the $\ell$-th smallest index $i\ge 1$ such that $\gamma$ hits $x$ during the time interval $[\t_{i},\s_{i}]$ and set ${\bf U}_{\ell}:=\t_{{\bf i}_{\ell}}$, ${\bf V}_{\ell}:=\s_{{\bf i}_{\ell}}$.

\end{notation}

\medskip

Let $D_{2}\subset D_{1}\subset D$ be nice domains with $d(D_{2},\partial D_{1})>0$.  Let $x\in D_{2}$ and $z\neq z'$ be    nice points of $\overline D$ different from $x$,   neither  in  $ \partial D_1 \backslash  {\mathcal J}(D,D_{1})$ nor in $\overline{D_2} $\footnote{In other words, $z, z'  \in    {\mathcal J}(D,D_{1}) \cup (\overline D\backslash (\{x\}\cup \partial D_1\cup\overline{D_2} ))$.  Note that  $z\not \in \partial D_1 \backslash  {\mathcal J}(D,D_{1})$ means that $z $ does not belong to $\partial D_1$ unless $z \in  {\mathcal J}(D,D_{1})$.}. 

\begin{figure}[h]
\centering
\includegraphics[scale=0.66]{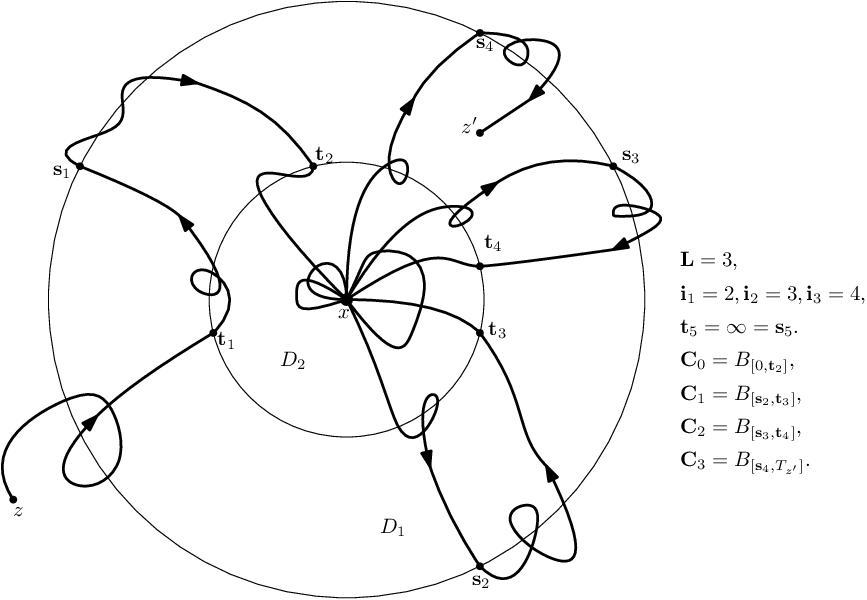}
\end{figure}

Let 
$$
{\bf L} := 1+ \#\{ \hbox{\rm loops from $x$ which hit $\partial D_{1}$} \} .
$$

\noindent By Lemma \ref{l:ito-muxx},
$$
\q^{z,z',\alpha}_{x,D}({\bf L}=L) =\ee^{-\alpha C_{D,D_{1}}(x)} {(\alpha C_{D,D_{1}}(x))^{L-1}\over (L-1)!}.
$$
 
\noindent Notice that ${\bf i}_{\ell}$ is well-defined for any $1\le \ell \le {\bf L}$. 
 %\noindent For any $\ell \le L$, denote by ${\bf i}_{\ell}$ the $\ell$-th smallest index $i\ge 1$ such that the Brownian motion hits $x$ during the time interval $[\t_{i},\s_{i}]$ and ${\bf U}_{\ell}:=\t_{{\bf i}_{\ell}}$, ${\bf V}_{\ell}:=\s_{{\bf i}_{\ell}}$. The first path from $z$ to $x$ is associated with ${\bf U}_{1}$, then each Brownian loop that hits  $\partial D_1$ is associated with a unique ${\bf U}_{\ell}$. 
 The next proposition gives the law of the points $(B_{ {\bf U}_\ell},\, B_{ {\bf V}_\ell} )_{\ell \le {\bf L}}$ under $\q^{z,z',\alpha}_{x,D}$. Notice that we may have 
 %${\bf i}_{1}=0$ in the case where $z\in D_{1}$ and the excursion from $z$ to $x$ does not hit $\partial D_{1}$.  We may also have
  $B_{ {\bf V}_{\bf L}}=z'$ in the case where $z'\in \overline{D_1}$ and the excursion from $x$ to $z'$ does not hit $\partial D_{1}$ (except, of course, possibly at the ending point).
 
\begin{proposition}\label{p:crossing}
Let $D_{2}\subset D_{1}\subset D$ nice domains with $d(D_{2},\partial D_{1})>0$. Let $x\in D_{2}$. Let $z\neq z'$ be nice points of $\overline D$ different from $x$,   neither  in  $\partial D_1 \backslash  {\mathcal J}(D,D_{1})$ nor in $\overline{D_2}$.\footnote{The assumption $z$, $z' \notin \partial D_1 \backslash  {\mathcal J}(D,D_{1})$, which may look somehow uncomfortable, will be automatically satisfied in the applications in the forthcoming sections. The same remark applies to other propositions and corollaries in this section.} Let $L\ge 1$ be  an integer. Let $y_{1},\ldots,y_{L}$ be nice points of $\partial D_{2}$. 
%We can take as well $y_{1}=z$ if $z\in D_{1}$.  
Let $y_{1}',\ldots,y_{L}'$ be nice points of $\partial D_{1}\backslash \partial D$. We can take as well $y_{L}'=z'$ in the case $z'\in D_{1}\cup {\mathcal J}(D,D_{1})$. Let $1\le i_{1}<\ldots <i_{L}$ be integers. Then
\begin{eqnarray*}
&& 
 	\q^{z,z',\alpha}_{x,D}\left( {\bf i}_{j}=i_{j},\, B_{ {\bf U}_j} \in \d y_{j},\, B_{ {\bf V}_j} \in \d y'_{j}, \, 1\le j\le L \mid {\bf L}=L \right) \\
&=&
	{C_{D,D_{1}}(x)^{-(L-1)}\over \xi_{D}(x,z,z')} \times \\
&&
	\prod_{j=1}^L \xi_{D_{1}}(x,y_{j},y_{j}')\,   \p^{z,z'}_{D}\left(  B_{\t_{i_{j}}} \in \d y_{j},\, B_{\s_{i_{j}}} \in \d y'_{j},\, 1\le j\le L ,\, \t_{i_{L}} < T_{z'} \right).
\end{eqnarray*}
In this equation, we mean, 
%in the case $y_{1}=z$, $B_{{\bf U}_{1}}=z$ and $B_{\t_{i_{1}}}=z$ when writing $B_{{\bf U}_{1}}\in \d y_{1}$ and $B_{\t_{i_{1}}} \in \d y_{1}$. We mean, 
in the case  $y'_{L}=z'$, $B_{{\bf V}_{L}}=z'$ and $B_{\s_{i_{L}}}=z'$  when writing $B_{{\bf V}_{L}}\in \d y'_{L}$ and $B_{\s_{i_{L}}} \in \d y_{L}'$.  \end{proposition}

Taking $L=1$ and integrating over $y_{j}$ and $y'_{j}$ yields the following equation which will be used in \eqref{eq:cond}:
\begin{equation}\label{c:crossing}
\xi_{D}(x,z,z') = \e^{z,z'}_{D}\Big[\sum_{i\ge 1} \xi_{D_{1}}(x,B_{\t_{i}},B_{\s_{i}}) 1_{\{\t_{i} < T_{z'}\}}\Big].
\end{equation}

\noindent {\it Proof of Proposition \ref{p:crossing}}.  To deal with the case $z\in \partial D$, we extend the $h$-transform representation \eqref{htransform2} to $z\in \partial D$. To this end, we   abusively write, when $z\in \partial D$, $\p^{z}_{D}$ for $\nu_{D}(z)$ defined in \eqref{mudy} (which is not a finite measure). With this notation, \eqref{htransform2} still holds for  nice point $z$ of $\partial D$: this  comes from \eqref{conv-functionalzz'}, \eqref{htransform2}  (replacing   $z$ there by $z+\varepsilon {\bf n}_{z}$), and \eqref{mudy}.

Denote by $Q$ the probability expression $\q^{z,z',\alpha}_{x,D}(\ldots)$ in the proposition. For simplicity, we suppose first that $y_{L}' \neq z'$. We notice that for each $1\le \ell \le L-1$,  the trajectory between ${\bf V}_{\ell}$ and the next hitting time of $x$ is, conditionally on  $B_{ {\bf V}_{\ell}}=y'$, distributed as $\p^{y',x}_{D}$. Moreover, conditionally on ${\bf L}=L$, $B_{{\bf V}_{\ell}}$ is distributed as   ${\mathfrak e}_{T_{\partial D_{1}}}$ under  $\nu_D(x,x)( \cdot \mid T_{\partial D_{1}}< T_{x})$ when $\ell < L$, and as  $B_{T_{\partial D_{1}\cup\{z'\}}}$ under $\p^{x,z'}_{D}$ when $\ell= L$. 
Using the strong Markov property at the hitting times of $x$ and stopping times  ${\bf V}_{1}, ..., {\bf V}_{L-1}$, the probability expression $Q$ is
\begin{align*} 
& 
	\prod_{j=1}^{L-1} \p^{y'_{j-1},x}_{D}( {\bf i}_{1}=i_{j}-i_{j-1},\,B_{ {\bf U}_{1}} \in \d y_{j}) \, \nu_{D}(x,x)({\mathfrak e}_{T_{\partial D_{1}}} \in \d y_{j}' \mid B_{T_{\partial D_{1}}} < T_{x} ) \\
& 
	\times \, \p^{y'_{L-1},x}_{D}( {\bf i}_{1}=i_{L}-i_{L-1},\, B_{ {\bf U}_{1}} \in \d y_{L})\,  \p_D^{x,z'}( B_{T_{\partial D_1}\cup\{z'\}} \in \d y_L') 
\end{align*}

\noindent where we set $y'_{0}=z$ and $i_{0}=0$. We can write it as
\begin{equation} \label{eq:Q}
Q = \Pi_{1} \Pi_{2} 
\end{equation}

\noindent where 
\begin{eqnarray*}
\Pi_{1} &:=&  \prod_{j=1}^{L} \p^{y'_{j-1},x}_{D}( {\bf i}_{1}=i_{j}-i_{j-1},\,B_{ {\bf U}_{1}} \in \d y_{j}), \\
\Pi_{2}  &:=& \prod_{j=1}^{L-1} \nu_{D}(x,x)({\mathfrak e}_{T_{\partial D_{1}}} \in \d y_{j}' \mid T_{\partial D_{1}} < T_{x} )\, \p_D^{x,z'}( B_{T_{\partial D_1}\cup\{z'\}} \in \d y_L').
\end{eqnarray*}

\noindent By the strong Markov property at the stopping time $\t_{i}$, we have  for any $y'\in\partial D_{1}\cup\{z\}$, % xxx
$y\in\partial D_{2}$, %(or $y=z$ if $z\in D_{1}$),  
$\p^{y',x}_{D}( {\bf i}_1= i,\, B_{ {\bf U}_1} \in \d y)
=
\p^{y',x}_{D}(  B_{ \t_{i}} \in \d y,    \t_i < T_x  ) \,  \p^{y,x}_{D}(T_{x}<T_{\partial D_{1}})
$
and by Lemma \ref{l:hitting-h} (ii) (with $z=y$, $z'=x$ there), $\p^{y,x}_{D}(T_{x}<T_{\partial D_{1}}) =  {  H_{D_{1}}(y,x)\over   H_{D}(y,x) }$. We get that
$$
\p^{y',x}_{D}( {\bf i}_1= i,\, B_{ {\bf U}_1} \in \d y)
=
\p^{y',x}_{D}(  B_{ \t_{i}} \in \d y,  \t_i < T_x  )   \, {  H_{D_{1}}(y,x)\over   H_{D}(y,x) } 
$$

\noindent which is, using the $h$-transform \eqref{htransform2}, 
$$
 \p^{y'}_D(    B_{ \t_{i}} \in \d y, \t_i < T_{\partial D} ) \, {  H_{D_{1}}(y,x)\over   H_{D}(y',x) }. 
$$

\noindent Hence,
\begin{equation} \label{eq:pi1}
\Pi_{1} = \prod_{j=1}^{L}  \p^{y'_{j-1}}_D (  B_{ \t_{i_j-i_{j-1}}} \in \d y_j, \, \t_{i_j-i_{j-1}}< T_{\partial D} ){  H_{D_1}(y_j,x)\over   H_D(y_{j-1}',x) }.
\end{equation}

\noindent  On the other hand, the law of ${\mathfrak e}_{T_{\partial D_{1}}}$ under  $\nu_D(x,x)( \cdot \mid T_{\partial D_{1}}< T_{x})$ has, by Lemma \ref{l:ito-muxx}, density on $\partial D_{1}$
$$
{1\over C_{D,D_{1}}(x)} G_{D}(x,y')H_{D_{1}}(x,y') \d y',
$$ 

\noindent whereas, by Lemma \ref{l:hitting-h} (i) (with $z=x$, $z'=z'$ there), 
\begin{equation} \label{eq:dy'}
\p^{x,z'}_{D}(B_{T_{\partial D_{1}}} \in \d y', T_{\partial D_{1}}<T_{z'} )= {H_{D_{1}}(x,y')H_{D}(y',z')\over H_{D}(x,z')}\d y'.
\end{equation}

\noindent  Recalling $y'_L \neq z'$, we obtain that  
\begin{equation} \label{eq:pi2}
\Pi_{2} =  C_{D,D_{1}}(x)^{-(L-1)}   {H_{D_{1}}(x,y'_{L}) H_{D}(y'_{L},z') \over H_{D}(x,z')} \d y'_{L} \prod_{j=1}^{L-1}   G_{D}(x,y'_{j})H_{D_{1}}(x,y'_{j})\d y'_{j}.
\end{equation}

\noindent  Recall that $G_{D}(y,y')=2\pi H_{D}(y,y')$ when $y,y'\in D$. By \eqref{eq:Q}, \eqref{eq:pi1} and \eqref{eq:pi2}, we get that
\begin{align*}
	Q
=&
	(2\pi)^{L-1}C_{D,D_{1}}(x)^{-(L-1)}  {H_{D}(y'_{L},z')  \over H_{D}(z,x)H_{D}(x,z')} \times \\ 
&	
	\prod_{j=1}^{L} H_{D_{1}}(x,y_{j})H_{D_{1}}(x,y'_{j}) \, \d y'_j 
	\prod_{j=1}^{L}  \p^{y'_{j-1}}_D (  B_{ \t_{i_j-i_{j-1}}} \in \d y_j,\, \t_{i_j-i_{j-1}}< T_{\partial D} ) .
\end{align*}

\noindent In view of the definition of $\xi_{D_{1}}$ in \eqref{def-xi}, we have
\begin{eqnarray*}
&& (2\pi)^{L-1}{H_{D}(y'_{L},z')  \over H_{D}(z,x)H_{D}(x,z')} \prod_{j=1}^{L} H_{D_{1}}(x,y_{j})H_{D_{1}}(x,y'_{j}) \\
&=&
{H_{D}(y'_{L},z')  \over H_{D}(z,z')\xi_{D}(x,z,z')} \prod_{j=1}^{L} \xi_{D_{1}}(x,y_{j},y'_{j}) \prod_{j=1}^L H_{D_{1}}(y_{j},y'_{j}).
\end{eqnarray*}

\noindent Comparing with the statement of the Proposition, we see that it remains to prove that
\begin{eqnarray} \label{eq:Q2}
&& H_{D}(y'_{L},z') \, 
\prod_{j=1}^{L}  \p^{y'_{j-1}}_D (  B_{ \t_{i_j-i_{j-1}}} \in \d y_j,\, \t_{i_j-i_{j-1}}< T_{\partial D} ) \, H_{D_{1}}(y_{j},y'_{j})\,  \d y'_j \\
&=&
  H_{D}(z,z')\, \p^{z,z'}_{D}\left(  B_{\t_{i_{j}}} \in \d y_{j},\, B_{\s_{i_{j}}} \in \d y'_{j},\, 1\le j\le L ,\, \t_{i_{L}} < T_{z'}  \right). \nonumber
\end{eqnarray}

\noindent Consider the right-hand side. By the $h$-transform \eqref{htransform2} (we recall that $y'_{L}\neq z'$ hence $\s_{i_{L}}<T_{z'}$), it is equal to
$$
\p^{z}_{D}\left(  B_{\t_{i_{j}}} \in \d y_{j},\, B_{\s_{i_{j}}} \in \d y'_{j},\, 1\le j\le L,\, \s_{i_{L}} < T_{\partial D}  \right)H_{D}(y_{L}',z').
$$

\noindent Recall that $H_{D_{1}}(y,y')$ for $y\in D_{1}$ and $y'\in \partial D_{1}$ is the density at $y'$ of the harmonic measure of the Brownian motion starting at $y$. It follows that, using the strong Markov property at the hitting times $\s_{i_{j-1}}$, $j\ge 2$, and $\t_{i_{j}}$, $j\ge 1$,
\begin{eqnarray*}
&& \p^z_{D}\left( B_{\t_{i_{j}}} \in \d y_{j}, \, B_{\s_{i_{j}}} \in \d y'_{j},\,1\le j\le L ,\,\s_{i_{L}} < T_{\partial D} \right) \\
&=&
\prod_{j=1}^{L}  \p^{y'_{j-1}}_D (  B_{ \t_{i_j-i_{j-1}}} \in \d y_j,\, \t_{i_j-i_{j-1}}< T_{\partial D} ) \, H_{D_{1}}(y_{j},y'_{j}) \, \d y'_j
\end{eqnarray*}

\noindent which completes the proof in the case  $y'_{L}\neq z'$. In the remaining case $y'_{L}=z'$, the same proof applies by replacing  \eqref{eq:dy'} by 
$$
\p^{x,z'}_{D} (T_{z'} \le T_{\partial D_{1}} ) = {H_{D_{1}}(x,z')\over H_{D}(x,z')},
$$

\noindent which holds by Lemma \ref{l:hitting-h} (ii) (with $z=x$ and $z'=z'$ there), then replacing \eqref{eq:Q2} by
\begin{eqnarray} 
&& 
\prod_{j=1}^{L}  \p^{y'_{j-1}}_D (  B_{ \t_{i_j-i_{j-1}}} \in \d y_j,\, \t_{i_j-i_{j-1}}< T_{\partial D} ) \, H_{D_{1}}(y_{j},y'_{j}) \, \prod_{j=1}^{L-1}\d y'_j  \nonumber \\
&=&
  H_{D}(z,z')\, \p^{z,z'}_{D}\left(  B_{\t_{i_{j}}} \in \d y_{j},\, B_{\s_{i_{j}}} \in \d y'_{j},\, 1\le j\le L ,\, \t_{i_{L}} < T_{z'}  \right)   \label{eq:Q2new}
\end{eqnarray}

\noindent where we write $B_{\s_{i_{L}}} \in \d y'_{L}$ for $B_{\s_{i_{L}}}=z'$.

As in the case $y'_L\neq z'$,  it   remains to check \eqref{eq:Q2new}.  By the strong Markov property at time $\t_{i_{L}}$ and Lemma \ref{l:hitting-h} (ii) (with $z=y_{L}$ and $z'=z'$ there), the right-hand side in \eqref{eq:Q2new}  is 
\begin{eqnarray*}
&& 
	H_{D}(z,z') \times \\
&&	
	\p^{z,z'}_{D}\left(  B_{\t_{i_{j}}} \in \d y_{j},\, B_{\s_{i_{j}}} \in \d y'_{j},\, 1\le j\le L-1 ,\, B_{\t_{i_{L}}}\in \d y_{L}, \,\t_{i_{L}} < T_{z'}  \right) \times \\
&&	
	{H_{D_{1}}(y_{L},z') \over H_{D}(y_{L},z)}
\end{eqnarray*}

\noindent which, in view of the $h$-transform \eqref{htransform2}, is
\begin{eqnarray*}
&&
	\p^{z}_{D}\left(  B_{\t_{i_{j}}} \in \d y_{j},\, B_{\s_{i_{j}}} \in \d y'_{j},\, 1\le j\le L-1 ,\, B_{\t_{i_{L}}}\in \d y_{L}, \,\t_{i_{L}} < T_{\partial D}  \right) \times 
	\\
&&
H_{D_{1}}(y_{L},z'). 
\end{eqnarray*}

\noindent  We use, as before,  the strong Markov property at times $\s_{i_{j-1}},\, j\ge 2$ and $\t_{i_{j}}$, $j\ge 1$ to get \eqref{eq:Q2new}. This completes the proof of Proposition \ref{p:crossing}. \hfill $\Box$

\medskip

Our next aim is to compute the Radon--Nikodym derivative of $\q^{z,z',\alpha}_{x,D}$ with respect to $\p^{z,z'}_D$ when both measures are restricted to an appropriate sigma-algebra denoted by $\F_{D_{1}}^+$ (defined in Notation \ref{notation:tribu2} below). We start by computing the conditional law of $\q^{z,z',\alpha}_{x,D}$ given ${\bf L}$.

Take the notation of Proposition \ref{p:crossing}. Fix an arbitrary $L\ge 1$ as well as arbitrary integers $1 \le i_1<i_2 < \cdots < i_L$. On the event $\{\t_{i_{L}} < T_{z'}\}$, let ${\bf C}_{j}$ be the trajectory from $\s_{i_{j}}$ to $\t_{i_{j+1}}$ for $1\le j\le L-1$, ${\bf C}_{0}$ the trajectory from (time) $0$ to $\t_{i_{1}}$, and ${\bf C}_{L}$ the trajectory from $\s_{i_{L}}$ to $T_{z'}$. [In case $\t_{i_{L}} \ge T_{z'}$, ${\bf C}_{j}$ will play no role, and can be defined as any trajectories.] 

Observe that, under $\q^{z,z',\alpha}_{x,D}$, conditionally on ${\bf L}=L$ and on ${\bf i}_{j}=i_{j}$, $B_{{\bf U}_j} \in \d y_{j}$, $B_{{\bf V}_j}\in \d y'_{j}$, $1\le j\le L$ (in particular, $\t_{i_{L}} < T_{z'}$), the trajectories $({\bf C}_{j},\, 0\le j\le L)$ are independent; for $0\le j\le L-1$, ${\bf C}_{j}$ is distributed as a Brownian excursion in $D$ from $y'_{j}$ to $y_{j+1}$ conditioned on $\t_{i_{j+1}-i_{j}}=T_{y_{j+1}}$   (with $i_{0}:=0$ and $y'_{0}:=z$); ${\bf C}_{L}$ is a Brownian excursion in $D$ from $y'_{L}$ to $z'$. [Beware of the degenerate situation that %${\bf C}_{0}=\{z\}$ if $i_{1}=0$ and that 
${\bf C}_{L}=\{z'\}$ if  $y'_{L}=z'$.] In other words, $({\bf C}_{j},\, 0\le j \le L)$ has the same distribution under $\q^{z,z',\alpha}_{x,D}(\cdot \mid {\bf L}=L, {\bf i}_{j}=i_{j},\,B_{{\bf U}_j} \in \d y_{j}, B_{{\bf V}_j}\in \d y'_{j},\, 1\le j\le L) $ and under $\p^{z,z'}_{D}(\cdot \mid B_{ \t_{i_{j}}} \in \d y_{j}, B_{\s_{i_{j}}}\in \d y'_{j},\, 1\le j\le L,\, \t_{i_{L}} <T_{z'})$. It implies that for any measurable set $A\in \sigma({\bf C}_{j},\, 0\le j \le L)$, 
\begin{eqnarray} \label{QA}
&&
\q^{z,z',\alpha}_{x,D}(A \mid {\bf L}=L, {\bf i}_{j}=i_{j},\,B_{{\bf U}_j} \in \d y_{j}, B_{{\bf V}_j}\in \d y'_{j},\, 1\le j\le L) \\
&=&
\p^{z,z'}_{D}(A \mid B_{ \t_{i_{j}}} \in \d y_{j}, B_{\s_{i_{j}}}\in \d y'_{j},\, 1\le j\le L,\, \t_{i_{L}} <T_{z'}). \nonumber
\end{eqnarray}

\noindent We introduce the following notation. Let ${\mathcal K}$ denote as before the set of all parametrized continuous planar curves $\gamma$  defined on a finite time-interval $[0, t_\gamma]$ with $t_\gamma \in (0, \, \infty)$.  

\begin{notation}
\label{notation:tribu}

 Let $D_{1}\subset D$ be nice domains. Let $\gamma \in {\mathcal K}$.% such that $a:= \gamma(0)$ and $b:= \gamma(t_\gamma)$ are nice points of $\overline D$.  

{\rm (i)} Let ${\mathcal E}_{D_{1}}$ be the set of excursions %inside $D_{1}$  to which we add: if $a\in D_{1}$ the excursion from $a$ to $\partial D_{1}$  (if it exists); if $b \in D_{1}$, the excursion from $\partial D_{1}$ to $b$ (if it exists).
away from $D_1^c \cup \{ \gamma(0)\} \cup \{ \gamma(t_\gamma)\}$ in the sense of \eqref{excursion}; an element of ${\mathcal E}_{D_{1}}$ is called an excursion inside $D_1$.

{\rm (ii)} For $\eee\in {\mathcal E}_{D_{1}}$, call $\eee_{g}$ and $\eee_{d}$ its starting and ending points. 

{\rm (iii)} An excursion away from $\overline{D_1} \cup \{ \gamma(0)\} \cup \{ \gamma(t_\gamma)\}$ is called an excursion outside $D_1$.

%{\rm (iii)} Let $\F_{D_{1}}^0$ be the sigma-algebra generated by $(\eee_{g},\eee_{d})$ for $\eee\in{\mathcal E}_{D_{1}}$ with the order of their appearances.

\end{notation}

\begin{notation}
\label{notation:tribu2}

 Let $D_{1}\subset D$ be nice domains. Let $\F_{D_{1}}^+$ be the sigma-algebra generated by the excursions of Brownian motion outside $D_{1}$ %(including the ones, if they exist, from the initial position (if not lying in $\overline{D_1}$) to $\partial D_{1}$ and from $\partial D_{1}$ to the ending position (if not lying in $\overline{D_1}$)) and by the starting and ending points of excursions inside $D_1$ 
together
%$(\eee_{g},\eee_{d})$ for $\eee\in{\mathcal E}_{D_{1}}$ 
with the order of their appearances.\footnote{
For any $\gamma \in \mathcal{K}$, we can associate  a function $\phi(\gamma)= ((\eee_{g},\eee_{d})_{D_{1}}, \mathcal{R})$ where $(\eee_{g},\eee_{d})_{D_{1}}$ are starting and ending points of excursions in $D_{1}$ of $\gamma$ and $\mathcal{R}$ is the order relation defined as $\mathcal{R}(\eee_{g},\eee'_{g})=0$ if $\eee_{g}$ is visited before $\eee'_{g}$, and $1$ otherwise, for any $\eee_{g}$, $\eee'_{g}$ starting points of excursions in ${\mathcal E}_{D_{1}}$. Then the sigma-algebra generated by $(\eee_{g},\eee_{d})$ for $\eee\in{\mathcal E}_{D_{1}}$ with the order of their appearances, is the sigma-algebra generated by $\phi$.}

\end{notation}

Observe that for any $A\in \F_{D_{1}}^+$, $A \cap \{ \t_{i_{L}} <T_{z'}\} \in \sigma({\bf C}_{j},0\le j \le L)$. Hence, by Proposition \ref{p:crossing}, we have for any $A\in\F_{D_{1}}^+$,
\begin{eqnarray}
 && \q^{z,z',\alpha}_{x,D}\left(A, {\bf i}_{j}=i_{j},\, B_{ {\bf U}_j} \in \d y_{j},\, B_{ {\bf V}_j} \in \d y'_{j}, \, 1\le j\le L \mid {\bf L}=L \right)   \nonumber\\
&=&
{C_{D,D_{1}}(x)^{-(L-1)}\over \xi_{D}(x,z,z')} \prod_{j=1}^L \xi_{D_1}(x,y_{j},y_{j}')  \,
\times \nonumber
\\
&&  \p^{z,z'}_{D}\left(A,  B_{\t_{i_{j}}} \in \d y_{j},\, B_{\s_{i_{j}}} \in \d y'_{j},\, 1\le j\le L,\, \t_{i_{L}} < T_{z'}  \right), \label{cafe1}
\end{eqnarray}

\noindent which gives that
\begin{eqnarray}
    \label{eq:radon}
 &&\q_{x,D}^{z,z',\alpha}(A\mid {\bf L}=L) 
    \\
 &=& {  C_{D,D_{1}}(x)^{-(L-1)} \over \xi_{D}(x,z,z')} \e^{z,z'}_D\Big[ 1_{A} \sum_{1 \le i_{1}<\ldots < i_{L}} \prod_{j=1}^{L}\xi_{D_{1}}(x, B_{\t_{i_{j}}}, B_{\s_{i_{j}}}) 1_{\{\t_{i_{L}} < T_{z'}\}} \Big].
    \nonumber
\end{eqnarray}

%\noindent [The indicator function $1_{\{i_{1}=0,z\notin D_{1}\}^c}$  on the right-hand side of \eqref{eq:radon} is to ensure that we consider the case $i_{1}=0$ only when $z\in D_{1}$; see Proposition \ref{p:crossing}.]

Contrarily to the left-hand side in \eqref{eq:radon}, the right-hand side depends on $D_{2}$ via the stopping times. To get rid of the dependence on $D_{2}$, we want to compute
$$
\e^{z,z'}_D\Big[ \sum_{1 \le i_{1}<\ldots < i_{L} } \prod_{j=1}^{L}\xi_{D_{1}}(x, B_{\t_{i_{j}}}, B_{\s_{i_{j}}}) 1_{\{\t_{i_{L}} < T_{z'}\}} \mid \F_{D_{1}}^+\Big].
$$

\noindent Observe that any $(\eee_{g},\eee_{d})$ for $\eee\in {\mathcal E}_{D_{1}}$ is measurable with respect to %$\F_{D_{1}}^0$ by definition, hence also to 
$\F_{D_{1}}^+$. Conditionally on $\F_{D_{1}}^+$, the excursions between $\eee_{g}$ and $\eee_{d}$ for $\eee\in {\mathcal E}_{D_{1}}$ are independent Brownian excursions inside $D_{1}$ with law $\p^{\eee_{g},\eee_{d}}_{D_{1}}$. %\footnote{Recall that for $u=z$ or $z'$, if $u\in \partial D_{1}$, then $u\in {\mathcal J}(D,D_{1})$.} 
Notice that each excursion inside $D_{1}$ can be associated to at most one $\t_{i}$: for any excursion $\eee\in {\mathcal E}_{D_{1}}$, write (if it exists) $\t(\eee)$ for the time $\t_{i}$ associated, which is the hitting time of $\partial D_{2}$ by the excursion. %except for the excursion from $z$ to $\partial D_{1}$ in the case $z\in D_{1}$, for which $\t(\eee)=0$. 
Set $\t(\eee)=\infty$ otherwise. Similarly, we write 
$$
\s(\eee)
:=
\inf \{ t>\t(\eee): \, \eee_t \in \partial D_1 \cup \{ z'\} \} \, .
$$

\noindent %which stands for the return time to $\partial D_{1}$, except for the excursion from $\partial D_{1}$ to $z'$ in the case $z'\in D_{2}$ for which $\s(\eee)$ is the hitting time of $z'$. 
Consequently, 
\begin{eqnarray*}
&&
	\sum_{1 \le i_{1}<\ldots < i_{L} } \prod_{j=1}^{L}\xi_{D_{1}}(x, B_{\t_{i_{j}}}, B_{\s_{i_{j}}}) 
%1_{\{i_{1}=0,z\notin D_{1}\}^c}
1_{\{\t_{i_{L}} < T_{z'}\}} \\
&=&
	\sum_{ {\underline \eee}_{D_{1}}^L} \prod_{j=1}^L \xi_{D_{1}}(x, B_{\t(\eee^j)}, B_{\s(\eee^j)})1_{\{\t(\eee^j)<\infty\}} ,  
\end{eqnarray*}

\noindent where we write $\sum_{{\underline \eee}_{D_{1}}^L}$ as a short way for sum over ordered (distinct) excursions $(\eee^1,\ldots,\eee^L) \in ({\mathcal E}_{D_{1}})^L$.\footnote{The  excursions $\eee^1,\ldots,\eee^L$ are naturally ordered in terms of their appearances, in particular any set of $L$ distinct excursions appears only once in the sum.} Apply \eqref{c:crossing} to $D_{1}=D$, $z=\eee_{g}$, and $z'=\eee_{d}$. Notice that in this case, the sum $\sum_{i\ge 1} \ldots$ in \eqref{c:crossing} has at most one term.  This  yields that for any $\eee \in {\mathcal E}_{D_{1}}$,  
\begin{eqnarray*}
&& 
	\e^{z,z'}_{D}[ \xi_{D_{1}}(x, B_{\t(\eee)}, B_{\s(\eee)}) 1_{\{\t(\eee)<\infty\}}
 \mid \F_{D_{1}}^+] \\
&=&
\e^{\eee_{g},\eee_{d}}_{D_{1}}\Big[\sum_{i\ge 1} \xi_{D_{1}}(x,B_{\t_{i}},B_{\s_{i}})  \, 1_{\{\t_{i} < T_{\eee_{d}}\}}\Big]\\
 &=&
   \xi_{D_{1}}(x,\eee_{g},\eee_{d})
.
\end{eqnarray*}

\noindent Hence,
\begin{eqnarray}
\label{eq:cond}
&&
	\e^{z,z'}_D\Big[ \sum_{1\le i_{1}<\ldots < i_{L} } \prod_{j=1}^{L}\xi_{D_{1}}(x, B_{\t_{i_{j}}}, B_{\s_{i_{j}}}) 1_{\{\t_{i_{L}}< T_{z'}\}}\mid \F_{D_{1}}^+\Big] \\
&=&
	\sum_{{\underline \eee}_{D_{1}}^L} \prod_{j=1}^L \xi_{D_{1}}(x,\eee_{g}^j,\eee_{d}^j). \nonumber
\end{eqnarray}

\noindent For future use, we observe that the same argument also gives (noting that $\s_{i_L}< T_{z'}$ ensures $\t_{i_{L}}< T_{z'}$)
\begin{eqnarray}
 \label{eq:cond2}
 &&\e^{z,z'}_D\Big[ \sum_{1\le i_{1}<\ldots < i_{L} } \prod_{j=1}^{L}\xi_{D_{1}}(x, B_{\t_{i_{j}}}, B_{\s_{i_{j}}}) 1_{\{T_{\partial D_1} < \t_{i_1}\}} 1_{\{\s_{i_L}< T_{z'}\}}\mid \F_{D_{1}}^+\Big]
    \\
 &=& \sum_{{\underline \eee}_{D_{1}}^L} \prod_{j=1}^L \xi_{D_{1}}(x,\eee_{g}^j,\eee_{d}^j) \, 1_{ \{\eee_{g}^1\neq z, \, \eee_{d}^L \ne z'\}} .\nonumber 
\end{eqnarray}

\noindent Going back to \eqref{eq:radon},  the following proposition is already proved:

\begin{proposition}\label{p:radon}
Let $D_{1}\subset D$ be nice domains. Let $x\in D_{1}$. Suppose that $z$ and $z'$ are distinct  nice points of $\overline D$, different from $x$,   and not in $\partial D_{1} \backslash {\mathcal J}(D,D_{1})$. For any $A\in \F_{D_{1}}^+$, 
$$
\q_{x,D}^{z,z',\alpha}(A\mid {\bf L}=L) = {  C_{D,D_{1}}(x)^{-(L-1)} \over \xi_{D}(x,z,z')} \e^{z,z'}_D\Big[ 1_{A} \sum_{  \underline{\eee}^L_{D_{1}}} \prod_{j=1}^{L}\xi_{D_{1}}(x,{\mathfrak e}_{g}^j, {\mathfrak e}_{d}^{j}) \Big],
$$
where as before, $\sum_{  \underline{\eee}^L_{D_{1}}} $ means that the sum runs over all ordered (distinct) excursions $(\eee^1,\ldots,\eee^L)   \in ({\mathcal E}_{D_{1}})^L $. In particular, we have \begin{equation}\label{edzz'} 
\e^{z,z'}_D\Big[  \sum_{  \underline{\eee}^L_{D_{1}}} \prod_{j=1}^{L}\xi_{D_{1}}(x,{\mathfrak e}_{g}^j, {\mathfrak e}_{d}^{j}) \Big]
=
\xi_{D}(x,z,z')\, C_{D,D_{1}}(x)^{L-1} .
\end{equation}
\end{proposition}

\medskip

Recall that $C_{D,D_{1}}(x)=C_{D_{1}}(x)-C_{D}(x)$ by \eqref{CDD1}. Since ${\bf L}-1$ is by construction a Poisson random variable with parameter $\alpha\,  \nu_D(x,x)(T_{\partial D_{1}}<T_{x}) = \alpha \, C_{D, D_1}(x)$, we get the following corollary. 

\begin{corollary}\label{c:radon}
With the notation and assumptions of Proposition \ref{p:radon},
$$
\frac{\d \q^{z,z',\alpha}_{x,D}}{\d \p^{z,z'}_D}\, \big|_{ \F_{D_{1}}^+} = {M_{D_{1}}(x,\alpha)   \over\ee^{-\alpha C_{D}(x)}\xi_{D}(x,z,z')}, $$
where  \begin{equation} 
M_{D_{1}}(x,\alpha):= \ee^{-\alpha C_{D_{1}}(x)}\sum_{L \ge 1} {\alpha^{L-1}\over (L-1)!}\sum_{{\underline{\eee}_{D_{1}}^L}   } \prod_{j=1}^{L}\xi_{D_{1}}(x,{\mathfrak e}_{g}^j, {\mathfrak e}_{d}^j).  
\label{MDalpha}
\end{equation}
\end{corollary}

\medskip
 By definition,  $M_D(x, \alpha)=\ee^{-\alpha C_{D}(x)}\xi_{D}(x,z,z')$ under $\p^{z, z'}_D$.  
 \medskip
 
%  \medskip
%
%\begin{remark}\label{r:radon}
%We also deduce that for any $A\in \F_{D_{1}}^0$ and $B\in \F_{D_{1}}^+$, with the notation of the corollary,
%\begin{eqnarray*}
%\q_{x,D}^{z,z',\alpha}(A\cap B) 
%&=&
% \e^{z,z'}_D\Big[1_{A\cap B} {M_{D_{1}}(x,\alpha) \over M_{D}(x,\alpha)} \Big]
% \\
% &=&\e^{z,z'}_D\Big[ 1_{A} {M_{D_{1}}(x,\alpha) \over M_{D}(x,\alpha)}\p^{z,z'}_D(B\mid \F_{D_{1}}^0 )\Big] 
% \\
% &=& \e_{\q^{z,z',\alpha}_{x,D}}\Big[ 1_{A} \p^{z,z'}_D(B\mid \F_{D_{1}}^0 )\Big].
%\end{eqnarray*}
%Therefore, conditionally on $\F_{D_{1}}^0$, the law of $\p^{z,z'}_D$ and $\q^{z,z', \alpha}_{x,D}$ on $\F_{D_{1}}^+$ are the same.
%\end{remark}
%
\begin{corollary}\label{r:martingale}
With the same notation and assumptions of Proposition \ref{p:radon}. Let $D_{2}\subset D_{1}$ be another nice domain such that $x\in D_2$ and that $z$, $z'$ do not lie in $\partial D_2 \backslash {\mathcal J}(D,D_2)$.\footnote{It is elementary to check that if $z\notin \partial D_1 \backslash {\mathcal J}(D,D_1)$, then saying $z\notin \partial D_2 \backslash {\mathcal J}(D,D_2)$ and saying $z\notin \partial D_2 \backslash {\mathcal J}(D_1,D_2)$ are equivalent.}  Then
$$
\e^{z,z'}_D[ M_{D_{2}}(x,\alpha) |  \F_{D_{1}}^+ ] = M_{D_{1}}(x,\alpha) . 
$$
\end{corollary}

\noindent {\it Proof of Corollary \ref{r:martingale}}. Let $C\in \F_{D_{1}}^+$. Then $C\in  \F_{D_{2}}^+$, hence in view of Corollary \ref{c:radon}, we have 
$$
\e^{z,z'}_D\left[1_{C} {M_{D_{2}}(x,\alpha) \over M_{D}(x,\alpha)}\right]=\q_{x,D}^{z,z',\alpha}(C)=\e^{z,z'}_D\left[1_{C} {M_{D_1}(x,\alpha)\over M_{D}(x,\alpha)}  \right] 
$$ which implies Corollary \ref{r:martingale}. 
\hfill$\Box$

\medskip
 
 The rest of this section is devoted to Proposition \ref{p:martingale2}, which  controls the conditional expectation of  $M_{D_{2}}(\cdot,\alpha) $ under  $\q_{x,D}^{z,z',\alpha}$,  and   is  also  the main technical tool in the proof of the forthcoming Proposition \ref{p:Mfini}.  At first we  compute  the expectation of  
$ \sum_{  \underline{\eee}_{D_{1}}^L } \prod_{j=1}^{L}\xi_{D_{1}}(x,{\mathfrak e}_{g}^j, {\mathfrak e}_{d}^{j}) 1_{ \{\eee_{g}^1\neq z,\eee_{d}^L \ne z'\}}$ under different measures:

\begin{lemma}\label{l:equalities}
 Let $x\in D_{1}$.  Suppose that  $z,z'$ are   distinct nice points of $\overline D$, different from $x$,  and    do not belong to $\partial D_{1} \backslash {\mathcal J}(D,D_{1})$. We take the convention that $H_{D_{1}}(x,y)=0$ if $y\notin  \overline D_{1}$. We have the following equalities: For all $L\ge 1$,  
 \begin{align}
 \begin{split} \label{eq:radonzz'} 
 	\frac{C_{D,D_{1}}(x)^{-(L-1)}}{\xi_{D}(x,z,z') } \e^{z,z'}_{D} \Big[\sum_{  \underline{\eee}_{D_{1}}^L }  \prod_{j=1}^{L} 
	& \xi_{D_{1}}(x,{\mathfrak e}_{g}	^j, {\mathfrak e}_{d}^{j}) 1_{ \{\eee_{g}^1\neq z,\eee_{d}^L \ne z'\}}\Big]   \\
=
	\Big(1- 
	& {H_{D_{1}}(x,z)\over H_{D}(x,z)}\Big)\Big(1- {H_{D_{1}}(x,z')\over H_{D}(x,z')} \Big)  , 
\end{split}  
\\
\label{eq:radonzx}
	\e^{z,x}_{D} \Big[\sum_{  {\underline{\eee}^L_{D_{1}}}} \prod_{j=1}^{L}\xi_{D_{1}}(x,{\mathfrak e}_{g}^j, {\mathfrak e}_{d}^{j}) 1_{\{\eee_{g}^1\neq z, \eee_{d}^L\neq x\}}\Big] 
& = 
	\Big(1- {H_{D_{1}}(x,z)\over H_{D}(x,z)} \Big)C_{D,D_{1}}(x)^{L} , \\
 \label{eq:radonxx}
 	\nu_{D}(x,x) \Big[\sum_{ {\underline{\eee}_{D_{1}}^L}} \prod_{j=1}^{L}\xi_{D_{1}}(x,{\mathfrak e}_{g}^j, {\mathfrak e}_{d}^{j}) 1_{\{\eee_{g}^1\neq x, 	\eee_{d}^L\neq x\}}\Big] 
& =
 	C_{D,D_{1}}(x)^{L+1} ,
\end{align}
where as before, $\sum_{  \underline{\eee}^L_{D_{1}}} $ means that the sum runs over all ordered (distinct) excursions $(\eee^1,\ldots,\eee^L)   \in ({\mathcal E}_{D_{1}})^L$; if furthermore, $z\in D$, then
\begin{eqnarray} 
 &&C_{D,D_{1}}(x)^{-(L-1)}  \nu_{D}(z,z)\Big[\sum_{ {\underline{\eee}_{D_{1}}^L}} \prod_{j=1}^{L}\xi_{D_{1}}(x,{\mathfrak e}_{g}^j, {\mathfrak e}_{d}^{j}) 1_{\{\eee_{g}^1\neq z, \eee_{d}^L\neq z\}}\Big] 
    \nonumber
    \\
 &=& (2\pi )^2 \big(H_{D}(x,z)- H_{D_{1}}(x,z)\big)^2 .
    \label{eq:radonzz}
\end{eqnarray}

% \begin{align}
% \begin{split} \label{eq:radonzz'} 
%\Big(1- {H_{D_{1}}(x,z)\over H_{D}(x,z)}\Big)\Big(1- {H_{D_{1}}(x,z')\over H_{D}(x,z')} \Big)  
%&= 
% {  C_{D,D_{1}}(x)^{-(L-1)} \over \xi_{D}(x,z,z') } \e^{z,z'}_{D} \Big[\sum_{  \underline{\eee}_{D_{1}}^L } \prod_{j=1}^{L}\xi_{D_{1}}(x,{\mathfrak e}_{g}^j, {\mathfrak e}_{d}^{j}) 1_{ \{\eee_{g}^1\neq z,\eee_{d}^L \ne z'\}}\Big] , 
%\end{split}  
%\\
%\label{eq:radonzx}
%\Big(1- {H_{D_{1}}(x,z)\over H_{D}(x,z)} \Big)C_{D,D_{1}}(x)^{L} 
%&=     \e^{z,x}_{D} \Big[\sum_{  {\underline{\eee}^L_{D_{1}}}} \prod_{j=1}^{L}\xi_{D_{1}}(x,{\mathfrak e}_{g}^j, {\mathfrak e}_{d}^{j}) 1_{\{\eee_{g}^1\neq z, \eee_{d}^L\neq x\}}\Big], \\
% \label{eq:radonxx}
%C_{D,D_{1}}(x)^{L+1} &=      \nu_{D}(x,x) \Big[\sum_{ {\underline{\eee}_{D_{1}}^L}} \prod_{j=1}^{L}\xi_{D_{1}}(x,{\mathfrak e}_{g}^j, {\mathfrak e}_{d}^{j}) 1_{\{\eee_{g}^1\neq x, \eee_{d}^L\neq x\}}\Big] , \\
%  \begin{split} 
% \label{eq:radonzz} 
%(2\pi )^2 \Big(H_{D}(x,z)- H_{D_{1}}(x,z)\Big)^2 &=   C_{D,D_{1}}(x)^{-(L-1)}  \nu_{D}(z,z)\Big[\sum_{ {\underline{\eee}_{D_{1}}^L}} \prod_{j=1}^{L}\xi_{D_{1}}(x,{\mathfrak e}_{g}^j, {\mathfrak e}_{d}^{j}) 1_{\{\eee_{g}^1\neq z, \eee_{d}^L\neq z\}}\Big] . 
%\end{split}
%\end{align}
 \end{lemma}
 
 \noindent {\it Proof of Lemma \ref{l:equalities}}.  Let  $A$ be the event that, under $\q^{z,z',\alpha}_{x,D}$,   the first part of the trajectory    from $z$ to $x$  hits $\partial D_{1}$,  and  so does the last part from $x$ to $z'$.   
 
We have by Lemma \ref{l:hitting-h}  (ii), 
\begin{eqnarray}
\q^{z,z',\alpha}_{x,D}\big( A \mid {\bf L}= L \big)
&=&
\p^{z, x}_D( T_{\partial D_{1}} < T_x) \, \p^{x, z'}_D( T_{\partial D_{1}} < T_{z'})
\nonumber
\\
&=&
\Big(1- {H_{D_{1}}(x,z)\over H_{D}(x,z)}\Big)\Big(1- {H_{D_{1}}(x,z')\over H_{D}(x,z')} \Big).
\label{cafe1001}
\end{eqnarray}

Let $D_2$ be a nice domain such that $x\in D_2$, $\overline{D_2} \subset D_1$ and that $z$, $z' \notin D_2$. Observe that for all $1\le i_1<\ldots <i_L$, $A \cap \{{\bf L}=L, \, {\bf i}_{1}=i_{1},\,  {\bf i}_{L}=i_{L}\}= A(i_1, i_L) \cap \{{\bf L}=L, \, {\bf i}_{1}=i_{1},\,  {\bf i}_{L}=i_{L}\}$, with $A(i_1, i_L):= \{T_{\partial D_1} < {\bf t}_{i_1}, \, B_{{\bf s}_{i_L}} \neq z'\} \in \F_{D_1}^+$. Then
\begin{eqnarray*}
 &&\q^{z,z',\alpha}_{x,D}\big( A \mid {\bf L}= L \big)
    \\
 &=& \sum_{1\le i_1<\ldots <i_L} \q^{z,z',\alpha}_{x,D}\big( A(i_1, i_L), \, {\bf i}_j=i_j, 1\le j\le L \,  \mid {\bf L}= L \big)
    \\
 &=& \frac{C_{D,D_{1}}(x)^{-(L-1)}}{\xi_{D}(x,z,z')} \e^{z,z'}_D \Big( \sum_{1\le i_1<\ldots <i_L} 1_{\{ T_{\partial D_1} < {\bf t}_{i_1}, \, B_{{\bf s}_{i_L}} \neq z'\} } \prod_{j=1}^L \xi_{D_1}(x,y_j,y_j') \Big) ,
\end{eqnarray*}

\noindent the last identity being a consequence of \eqref{cafe1}. Applying \eqref{eq:cond2} gives
\begin{eqnarray*}
&&
	\q^{z,z',\alpha}_{x,D}\big( A \mid {\bf L}= L \big) \\
&=&
	\frac{C_{D,D_{1}}(x)^{-(L-1)}}{\xi_{D}(x,z,z') } \e^{z,z'}_{D} \Big[\sum_{  \underline{\eee}_{D_{1}}^L } \prod_{j=1}^{L}\xi_{D_{1}}(x,{\mathfrak e}_{g}^j, {\mathfrak e}_{d}^{j}) 1_{ \{\eee_{g}^1\neq z,\eee_{d}^L \ne z'\}}\Big],
\end{eqnarray*}

\noindent which is the same expression as on the left-hand side of \eqref{eq:radonzz'}. Together with \eqref{cafe1001}, this yields \eqref{eq:radonzx}.  

Recall the definition of $C_{D,D_{1}}(x)$ in \eqref{CDD1} and that $\lim_{z'\to x}  \big( G_{D}(x,z')-G_{D_{1}}(x,z')\big)=C_{D,D_{1}}(x)$ as proved in \eqref{gx}.  It follows that  $$
 \xi_D(x, z, z')  \, \Big(1- {H_{D_{1}}(x,z')\over H_{D}(x,z')} \Big)  \, \to \, C_{D, D_1}(x), \qquad \mbox{as } z'\to x. $$

\noindent Then making $z'\to x$ in   \eqref{eq:radonzz'}  implies \eqref{eq:radonzx}. 

Notice that as $z\to x$, $ 1- {H_{D_{1}}(x,z)\over H_{D}(x,z)}  = \frac{G_D(x, z)-G_{D_1}(x, z)}{G_D(x, z)} \sim \frac{C_{D, D_1}(x)}{\log 1/|x-z|}$ by \eqref{Gdxy}.  Using    \eqref{nuxx}, making $z\to x$ in \eqref{eq:radonzx} gives \eqref{eq:radonxx}.  
 
 Recall that $H_{D}(y,y')=G_{D}(y,y')/2\pi$ when $y$, $y'$ are in $D$.  When $z\in D\backslash \partial D_{1}$, making $z'\to z\neq x$ in \eqref{eq:radonzz'}, and using \eqref{Gdxy} and \eqref{nuxx}  give \eqref{eq:radonzz}. \hfill $\Box$

\medskip

Let us go back to the probability $\q^{z,z',\alpha}_{x,D}$ for $z\neq z'$ which are different from $x$.  For any nice domains $ S\subset D_{1}\subset D$, such that $x\in S$, we introduce ${\mathcal E}_{x,D_{1},S}$, the set of excursions away from $D_1^c \cup \{z\} \cup \{ z'\} \cup \{x\}$ (in the sense of \eqref{excursion}) but excluding all loops at $x$ that lie in $S$.\footnote{In words, ${\mathcal E}_{x,D_{1},S}$ stands for the set of excursions in ${\mathcal E}_{D_{1}}$ that do {\it not} hit $x$, to which we also add: (a) the excursions from $\partial D_{1}$ to $x$ (if they exist); (b) the excursions from $x$ to $\partial D_{1}$  (if they exist); (c) loops at $x$ which lie in $D_{1}$ and hit $\partial S$; (d) if $z\in D_{1}$, the path from $z$ to $x$ if it does not hit $\partial D_1$; (e) if $z'\in D_{1}$, the path from $x$ to $z'$ if it does not hit $\partial D_1$.}

For any nice domain set $S\subset D_{1}$ which contains $x$,  any nice domain $D_2 \subset D_{1}$, and any $u\in D_{2} \backslash \{ x\}$, let 
\begin{equation}\label{defMxDD}
M_{x,D_{1},S,D_{2}}(u,\alpha) :=  \ee^{-\alpha C_{D_{1}}(u)} \sum_{L=1}^\infty {\alpha^{L-1}\over (L-1)!}\, \sum_{ {\underline{\eee}^L_{x,D_1,S}}} \prod_{j=1}^{L}\xi_{x,D_{1},S,D_{2}}(u,{\mathfrak e}_{g}^j, {\mathfrak e}_{d}^j)
\end{equation}

\noindent where the sum $\sum_{ {\underline{\eee}^L_{x,D_1,S}}} $ runs over all ordered (distinct) excursions $(\eee^1,\ldots,\eee^L)$ in $({\mathcal E}_{x,D_1,S})^L$ and 
\begin{align} \nonumber
	\xi_{x,D_{1},S,D_{2}}(u, & \eee_{g},  \eee_{d}) := \\
&
	\begin{cases}
         \xi_{D_{1}}(u,\eee_{g},\eee_{d}), 
         & \mbox{if } \eee_{g}, \, \eee_{d}\in \partial D_{1} ,
       	 \\
          \frac{(2\pi)^2 ( H_{D_1}(u, x)  - H_{D_2}(u,x))^2}{C_{D_1, S}(x)}, 
          & \mbox{if }  \eee_{g}=\eee_{d}=x,  
          \\ 
          \frac{2\pi (H_{D_{1}}(u,\eee_{g})-H_{D_{2}}(u,\eee_{g}))(H_{D_{1}}(u,\eee_{d})-H_{D_{2}}(u,\eee_{d}))}{H_{D_{1}}(\eee_{g},\eee_{d})},
          &\mbox{otherwise}.
    \end{cases}
    \label{xi}  
\end{align}

%\noindent [In particular, $\xi_{x,D_{1},S,D_{2}}(u,\eee_{g},\eee_{d})$ coincides with $\xi_{D_{1}}(u,\eee_{g},\eee_{d})$ if both $\eee_{g}$ and $\eee_{d}$ are on $\partial D_{1}$.]

%$\xi_{x,D_{1},S,D_{2}}(u,\eee_{g},\eee_{d})$ is defined as   
%\begin{itemize}
%\item $\xi_{D_{1}}(u,\eee_{g},\eee_{d})$ if $\eee_{g}$ and $\eee_{d}$ are on $\partial D_{1}$;
%\item $2\pi (H_{D_{1}}(u,\eee_{g})-H_{D_{2}}(u,\eee_{g}))(H_{D_{1}}(u,x)-H_{D_{2}}(u,x))/H_{D_{1}}(\eee_{g},x)$ if $\eee_{d}=x$ and $\eee_{g}\in \partial D_1 \cup \{ z\}$ (in case $\eee_{g}=z$, we necessarily have $z\in \overline{D_1}$);
%\item $2 \pi (H_{D_{1}}(u,\eee_{d})-H_{D_{2}}(u,\eee_{d}))(H_{D_{1}}(u,x)-H_{D_{2}}(u,x))/H_{D_{1}}(x,\eee_{d})$ if $\eee_{g}=x$ and $\eee_{d} \in \partial D_1\cup \{ z'\}$ (in case $\eee_d=z'$, we necessarily have $z'\in \overline{D_1}$);
%\item $(2\pi)^2  \, \big(H_{D_{1}}(u,x)- H_{D_{2}}(u,x)\big)^2  / C_{D_{1},S}(x)$ if $\eee_{g}=\eee_{d}=x$.
%\end{itemize}

%For future use, we make the simple observation that by definition, if $H_{D_{2}}(u,\, y) =0$ for all $y\in \partial D_1\cup \{z\} \cup \{ z'\} \cup \{ x\}$, then
%\begin{equation}
%    \xi_{x,D_{1},S,D_{2}}(u,\eee_{g},\eee_{d})
%    = 
%    \begin{cases}
%       (2\pi)^2 ( H_{D_1}(u, x) )^2/C_{D_1, S}(x), \qquad & \mbox{if }  \eee_{g}=\eee_{d}=x,  
%       \\
%       \xi_{D_1}(u, \eee_{g},\eee_{d}), \qquad &\mbox{otherwise}.
%    \end{cases}
%    \label{xi}
%\end{equation}

Considering \eqref{defMxDD}, and the inequality $\sum_{j\ge 0} {\lambda^j \over (j!)^2} \le \ee^{2\sqrt{\lambda}}$ for any $\lambda\ge 0$, we have the following bound:
\begin{equation}\label{boundMxDD}
M_{x,D_{1},S,D_{2}}(u,\alpha) \le \varUpsilon_{x,D_{1},S,D_{2}}(u) \,  \ee^{-\alpha C_{D_{1}}(u)} \, \ee^{2\sqrt{\alpha \varUpsilon_{x,D_{1},S,D_{2}}(u)}}
\end{equation}

\noindent where
\begin{equation}\label{defSigma}
\varUpsilon_{x,D_{1},S,D_{2}}(u) := \sum_{\eee \in {\mathcal E}_{x,D_{1},S}} \xi_{x,D_{1},S,D_{2}}(u,\eee_{g},\eee_{d}).
\end{equation}

\begin{proposition}\label{p:martingale2}
 Let $S \subset D_{1}\subset D$ be nice domains  such that $x\in S$. Let $D_2\subset D_1$ be a nice domain satisfying $\partial D_2 \subset S^c$. Suppose that $z$ and $z'$ are distinct nice points of $\overline D$, different from $x$, which belong neither to $\partial D_{1} \backslash {\mathcal J}(D,D_{1})$ nor to $\partial D_2 \backslash {\mathcal J}(D_1,D_2)$. Let $u\in D_2 \backslash \{ x\}$. Define
 $$
 \widetilde M_{x,D_{2}}(u,\alpha) :=\ee^{-\alpha C_{D_{2}}(u)}\sum_{L=1}^\infty {\alpha^{L-1}\over (L-1)!}\sum_{\underline{\eee}_{D_{2}}^L} \prod_{j=1}^{L}\xi_{D_{2}}(u,{\mathfrak e}_{g}^j, {\mathfrak e}_{d}^j) 1_{\{ x\notin \eee^j\}}.
 $$
 We have, with the notation of \eqref{defMxDD} and \eqref{defSigma},
 \begin{eqnarray}
  &&\e_{\q^{z,z',\alpha}_{x,D}}[ \widetilde M_{x,D_{2}}(u,\alpha) \, | \, (\eee_{g},\eee_{d})_{\eee\in{\mathcal E}_{x,D_{1},S}}] 
     \nonumber
     \\
  &=&  M_{x,D_1,S,D_{2}}(u,\alpha)
     \label{cafe108}
     \\
  &\le&\varUpsilon_{x,D_{1},S,D_{2}}(u) \, \ee^{-\alpha C_{D_{1}}(u)} \, \ee^{2\sqrt{\alpha \varUpsilon_{x,D_{1},S,D_{2}}(u)}}.
     \label{cafe109}
 \end{eqnarray}

\end{proposition}
 
\noindent {\it Proof}. %Conditionally on $(\eee_{g},\eee_{d})_{\eee\in{\mathcal E}_{x,D_{1},S}}$, the excursions of ${\mathcal E}_{x,D_{1},S}$ are independent Brownian excursions. 
%We decompose the sum $\widetilde M_{x,D_{2}}(u,\alpha)$ with respect to the excursions of ${\mathcal E}_{x,D_{1},S}$. 
Any element of ${\mathcal E}_{D_2}$ that does not hit $x$, necessarily hits $\partial D_2$, and lies (except, possibly, for the starting or the ending point) in $D_1 \backslash \{ x, \, z, \, z'\}$, and is thus contained in an excursion, say $\eee$, away from $D_1^c \cup \{x, \, z, \, z'\}$; $\eee$ cannot be a loop at $x$ lying in $S$ because $\partial D_2 \subset S^c$. In other words, $\eee \in {\mathcal E}_{x,D_1,S}$. It follows that
\begin{eqnarray*}
&& 
	\widetilde M_{x,D_{2}}(u,\alpha)  \\
&=&
	\ee^{-\alpha C_{D_{2}}(u)}\sum_{L=1}^\infty {\alpha^{L-1} \over (L-1)!} \sum_{K=1}^L \sum_{ {\underline{\eee}^K_{x,D_{1},S}}} \sum_{\ell_1+\ldots+	\ell_K=L, \; \ell_j \ge 1, \, \forall j\le K} \, \prod_{j=1}^K \varUpsilon(\eee^j,\ell_j) , 
\end{eqnarray*}

\noindent where   as before, the sum $\sum_{ {\underline{\eee}^K_{x,D_1,S}}} $ runs over ordered (distinct) all excursions $(\eee^1,\ldots,\eee^K) \in ({\mathcal E}_{x,D_1,S})^K$,  and  $\varUpsilon(\eee^j,\ell_{j})$ is the sum, over all $\ell_{j}$ distinct excursions $\eee^{j,1},\ldots,\eee^{j,\ell_{j}}$ in ${\mathcal E}_{D_{2}}$ which do not hit $x$ and which is contained in the excursion $\eee^j$, of the products $\prod_{i=1}^{\ell_{j}}\xi_{D_{2}}(u,\eee^{j,i}_{g},\eee^{j,i}_{d})$.  

Notice that conditionally on $ (\eee_{g},\eee_{d})_{\eee\in{\mathcal E}_{x,D_{1},S}} \in (\C^2)^{{\mathcal E}_{x,D_{1},S}}$, the family $\big(\varUpsilon(\eee^j,\ell_{j})\big)_{1\le j \le K}$ for distinct $\eee^1,\ldots,\eee^K \in {\mathcal E}_{x,D_1,S}$, is independent. We claim that 
\begin{equation}
    \label{2401}
    \e_{\q^{z,z',\alpha}_{x,D}} \Big[  \varUpsilon(\eee^j,\ell_{j}) \, \big|\,   (\eee_{g},\eee_{d})_{\eee\in{\mathcal E}_{x,D_{1},S}}\Big] 
    =
    C_{D_{1},D_{2}}(u)^{\ell_{j}-1} \,  \xi_{x,D_{1},S,D_{2}}(u,\eee_{g}^j,\eee_{d}^j).
\end{equation}

Indeed by definition of ${\mathcal E}_{x,D_1,S}$,  the excursion $\eee^j$ can be an excursion from $\partial D_{1}$ to $\partial D_{1}$  which does not hit $x$, or from $\partial D_{1}$ to $x$ or from $x$ to $\partial D_{1}$, or a loop at $x$ in $D_{1}$ conditioned to hit $\partial S$.  Let us check \eqref{2401} for each of these four cases. Write $Q_{\eqref{2401}}$  for the conditional expectation term on the left-hand side of \eqref{2401}. 

If $\eee^j$ is  an excursion from $\partial D_{1}$ to $\partial D_{1}$ which does not hit $x$, then $\eee^j_g \in \partial D_1$,  $\eee^j_d \in \partial D_1$ and none of $\eee^{j,1},\ldots,\eee^{j,\ell_{j}}$   hits $x$. By the Markov property,
$$
	Q_{\eqref{2401}}
= 
	\e^{\eee^j_g, \eee^j_d}_{D_1}\Big[   \sum_{ {\underline{\eee}_{D_{2}}^{\ell_j}}} \prod_{j=1}^{\ell_j}\xi_{D_{2}}(u,{\mathfrak e}_{g}^j, {\mathfrak e}_{d}^{j})     
	\Big] 
=
	C_{D_{1},D_{2}}(u)^{\ell_{j}-1} \,  \xi_{D_{1}}(u,\eee_{g}^j,\eee_{d}^j),
$$

\noindent  by using  \eqref{edzz'} for the second equality.  This implies \eqref{2401} as $ \xi_{x,D_{1},S,D_{2}}(u,\eee_{g}^j,\eee_{d}^j)=\xi_{D_{1}}(u,\eee_{g}^j,\eee_{d}^j)$ in this case.

If $\eee^j$ is an excursion from $\partial D_{1}$ to $x$, then $\eee^j_g \in \partial D_1$ and $\eee^j_d=x$. By the Markov property and  \eqref{eq:radonzz'},
$$
Q_{\eqref{2401}}
=
C_{D_{1},D_{2}}(u)^{\ell_{j}-1} \,  \xi_{D_{1}}(u,\eee_{g}^j,\eee_{d}^j)  \Big( 1- \frac{H_{D_2}(u, \eee_{g}^j)}{H_{D_1}(u, \eee_{g}^j)}\Big)  \Big( 1- \frac{H_{D_2}(u, \eee_{d}^j)}{H_{D_1}(u, \eee_{d}^j)}\Big),
$$

\noindent which yields \eqref{2401}  by the definition of $\xi_{x,D_{1},S,D_{2}}(u,\eee_{g}^j,\eee_{d}^j)$ (recalling that $\eee^j_g \neq x$ and $\eee^j_d=x$ in this case). The case when $\eee^j$ is an excursion from  $x$ to $\partial D_{1}$ follows from the same way.

Finally for the case when $\eee^j$ is a loop at $x$ in $D_{1}$ conditioned to hit $\partial S$, $\eee_{g}^j=\eee_{d}^j=x$, so
\begin{eqnarray*}
    Q_{\eqref{2401}}
 &=& \nu_{D_1}(x,x) \Big[   \sum_{ {\underline{\eee}_{D_{2}}^{\ell_j}}} \prod_{j=1}^{\ell_j}\xi_{D_{2}}(u,{\mathfrak e}_{g}^j, {\mathfrak e}_{d}^{j}) 1_{\{\eee_{g}^1\neq x, 	\eee_{d}^{\ell_j}\neq x\}}   \, \Big| \, T_{\partial S} < T_x \Big] 
    \\
 &=& \frac{\nu_{D_1}(x,x) \Big[\sum_{ {\underline{\eee}_{D_{2}}^{\ell_j}}} \prod_{j=1}^{\ell_j}\xi_{D_{2}}(u,{\mathfrak e}_{g}^j, {\mathfrak e}_{d}^{j})1_{\{\eee_{g}^1\neq x, 	\eee_{d}^{\ell_j}\neq x\}} \Big]}{\nu_{D_1}(x,x) (T_{\partial S} < T_x)} \, .
\end{eqnarray*}

\noindent Recall from \eqref{CDD1} that $\nu_{D_{1}}(x,x)(T_{\partial S}<T_{x}) = C_{D_{1},S}(x)$. Using \eqref{eq:radonzz} completes the proof of \eqref{2401}. 

By \eqref{2401},   $ \e_{\q^{z,z',\alpha}_{x,D}}[ \widetilde M_{x,D_{2}}(u,\alpha) \, | \, (\eee_{g},\eee_{d})_{\eee\in{\mathcal E}_{x,D_{1},S}} ] $ is given by
\begin{align*}
 	\ee^{-\alpha C_{D_{2}}(u)}  \sum_{L=1}^\infty {\alpha^{L-1} \over (L-1)!} 
	& \sum_{K= 1}^L  C_{D_{1},D_{2}}(u)^{L-K} \times   \\ 
&	
	\Big\{  \sum_{ \underline{\eee}_{x,D_{1},S}^K}  \sum_{ \ell_1+\ldots+	\ell_K=L, \; \ell_j \ge 1, \, \forall j\le K}  \prod_{j=1}^K \xi_{x,D_{1},S,D_{2}}(u,\eee_{g}^j,\eee_{d}^j) \Big\}.
\end{align*}

\noindent The sum $ \sum_{L=1}^\infty \ldots$ is also
\begin{eqnarray*}
 &&
 	\sum_{L=1}^\infty {\alpha^{L-1} \over (L-1)!} \sum_{K=1}^L C_{D_{1},D_{2}}(u)^{L-K} 
 	 \sum_{   \underline{\eee}_{x,D_{1},S}^K } \binom{L-1}{K-1} \prod_{j=1}^K \xi_{x,D_{1},S,D_{2}}(u,\eee_{g}	^j,\eee_{d}^j)
    \\
 &=& \sum_{K=1}^\infty \sum_{ \underline{\eee}_{x,D_{1},S}^K}   \prod_{j=1}^K \xi_{x,D_{1},S,D_{2}}(u,\eee_{g}^j,\eee_{d}^j) \sum_{L=K}^\infty {\alpha^{L-1} \over (K-1)! (L-K)!} C_{D_{1},D_{2}}(u)^{L-K} . 
\end{eqnarray*}

\noindent Since
$$
\sum_{L=K}^\infty {\alpha^{L-1} \over (K-1)! (L-K)!} C_{D_{1},D_2}(u)^{L-K} = {\alpha^{K-1} \over (K-1)!}\ee^{\alpha C_{D_{1},D_2}(u)} ,
$$

\noindent and $C_{D_{1},D_2}(u)=C_{D_{2}}(u)-C_{D_{1}}(u)$, we get \eqref{cafe108}, which, in turn (in view of \eqref{boundMxDD}), yields \eqref{cafe109}.\hfill $\Box$

\section{Construction of the measure ${\mathcal M}^\alpha_\infty$}\label{s:construction}

Let $D$ be a simply connected nice domain. Let $z$ and $z'$ be distinct nice points of $\overline D$. For $n\ge 0$, we let $\D_{n}$ be the set of the connected components of $D$ minus the grid of mesh size $2^{-n}$. We choose $(\D_n, \, n\ge 0)$ such that $z$ and $z'$ never lie in any of the grid. We will abusively call an element of $\D_{n}$ a square (it is not necessarily a square near the boundary of $D$). We say that $x$ is {\it suitable} if $x$ is in some square $D_{n}\in \D_{n}$ at any level $n$, and is different from $z$ and $z'$. The set of suitable points has    full Lebesgue measure.

For any suitable $x$, we let $D_{n}^{(x)}$ be the square in $\D_{n}$ that contains $x$. We let $\F_{\D_{n}}$
%$\F_{\D_{n}}:=\bigcap_{D_{n}\in\D_{n}} \F_{D_{n}}^+$ 
be the sigma-algebra generated by the starting and return points  of all excursions inside some square $D_{n}$ together with the order of their appearances. %Notice that 
Then $(\F_{\D_{n}})_{n\ge 0}$ is a filtration, and $\sigma (\cup_{n=0}^\infty \F_{\D_{n}})$ %is the sigma-algebra generated by the whole sigma-algebra ${\mathscr F}$
coincides with $\sigma(B_t, \, t\in [0, \, T_{z'}]).$\footnote{
To identify the two sigma-algebras, note that by continuity, it suffices to show that for any disc ${\cal B} (x, \, r) \subset D$ (where $r>0$ is rational, and $x$ is with rational coordinates), the duration and the exiting position of Brownian motion starting at any $y\in {\cal B} (x, \, r)$ and killed upon exiting from ${\cal B} (x, \, r)$ are measurable with respect to $\sigma (\cup_{n=0}^\infty \F_{\D_{n}})$. This, however, is quite straightforward because by continuity, the number of crossings and their positions between two concentric circles (hence the local time on any circle, hence the duration by integration over local time) are measurable with respect to $\sigma (\cup_{n=0}^\infty \F_{\D_{n}})$.
} We define for any Borel set $A$,
\begin{equation}\label{measuredna}
\mathcal{M}_{\D_{n}}^\alpha(A) := \int_{A} {M}_{D_{n}^{(x)}}(x, \,\alpha) \,\d x,
\end{equation}

\noindent where  ${M}_{D_{n}^{(x)}}(x,\alpha)$ is  defined in \eqref{MDalpha} and ${M}_{D_{n}^{(x)}}(x,\alpha):=0$ if $x$ is not suitable.

%\begin{proposition}
%Let $A$ be a Borel set in $D$. The process $(\mathcal{M}_{\D_{n}}^\alpha(A))_{n\ge 0}$ is a $(\F_{\D_{n}})$-martingale under $\p^{z,z'}_{D}$. %More generally, for any   Borel function $f:  D \to \r_+$,  the process $\big( \int  {M}_{D,D_{n}^{(x)}}(x,\alpha) f(x)  \d x \big)_{n\ge 0}$  is a $(\F_{\D_{n}})$-martingale under $\p^z$.
%\end{proposition}
%{\it Proof}.  Note that ${\mathcal M}_{\D_{n}}^\alpha(A)$ is $\F_{\D_{n}}$-adapted.    Use Remark \ref{r:martingale} and integrate over $A$ yield   the Proposition. \hfill $\Box$  

Note that ${\mathcal M}_{\D_{n}}^\alpha(A)$ is $\F_{\D_{n}}$-adapted. For any suitable $x$, $\F_{\D_{n}} \subset \F_{D_n^{(x)}}^+$ (defined in Notation \ref{notation:tribu2} for the latter), so we can use Corollary  \ref{r:martingale} and integrate over $A$ to yield that  $(\mathcal{M}_{\D_{n}}^\alpha(A))_{n\ge 0}$ is an $(\F_{\D_{n}})$-martingale under $\p^{z,z'}_{D}$. Consequently, the following limit exists:
\begin{equation}
    \label{MinfiniA}
    {\mathcal M}_{\infty}^\alpha(A) 
    :=
    \lim_{n\to+\infty} {\mathcal M}_{\D_{n}}^\alpha(A)
\in [0, \, \infty), \qquad  \mbox{$\p^{z,z'}_D$-a.s.} 
\end{equation}

\begin{theorem}
 \label{t:conv-p}

 Let $D$ be a simply connected nice domain, and $z$ and $z'$ be distinct nice points of $\overline D$. Fix $0\le \alpha <2$ and a Borel set $A$   of $D$. Under $\p^{z,z'}_{D}$, the martingale $({\mathcal M}_{\D_{n}}^\alpha(A))_{n \ge 0}$ converges in $L^1$ to ${\mathcal M}_{\infty}^\alpha(A)$.

\end{theorem}

In order to prove Theorem \ref{t:conv-p}, we first remark that by Corollary \ref{c:radon}, the (finite) measure $\int_{A} \q_{x,D}^{z,z',\alpha}(\bullet)\, M_D(x, \, \alpha) \, \d x$ has Radon--Nikodym derivative ${\mathcal M}_{\D_{n}}^\alpha(A)$ with respect to $\p^{z,z'}_{D}$ on $\F_{\D_{n}}$. We assume without loss of generality that $A$ has a positive Lebesgue measure.  Let
\begin{equation}
    \label{QAzz'}
    Q_{A}^{z, z', \alpha}(\bullet)
    := 
    \frac1{ \int_{A}   M_D(x, \alpha) \d x}  \int_A \q_{x,D}^{z,z',\alpha} (\bullet)   \, M_D(x, \alpha) \d x
\end{equation} 

\noindent be the  normalized probability measure defined on $\sigma (\cup_{n=0}^\infty \F_{\D_{n}})$. 

By an elementary fact (see Durrett \cite{Durrett}, Theorem 5.3.3),  ${\mathcal M}_{\D_{n}}^\alpha(A)$ converges in $L^1(\p^{z,z'}_{D})$  if and only if  ${\mathcal M}_{\infty}^\alpha(A)< \infty$,  $Q_{A}^{z, z', \alpha}$-a.s.,  
 where ${\mathcal M}_{\infty}^\alpha(A)$ is the $Q_{A}^{z, z', \alpha}$-a.s.\ limit of ${\mathcal M}_{\D_{n}}^\alpha(A)$, which exists since $1/{\mathcal M}_{\D_{n}}^\alpha(A)$ is a nonnegative supermartingale under $Q_{A}^{z, z', \alpha}$.

%Hence  Theorem \ref{t:conv-p} follows once  we  have proved   that $$ \q_{x,D}^{z,z',\alpha}({\mathcal M}_{\infty}^\alpha(D)<\infty)=1, \qquad \mbox{Lebesgue-a.e. }  x\in D. $$

  Then Theorem \ref{t:conv-p}    follows from the next proposition.
 
 \begin{proposition}
     \label{p:Mfini}
     
     Let $D$ be a simply connected nice domain, and $z$ and $z'$ be distinct nice points of $\overline D$. Fix $0\le \alpha <2$. Then 
 \begin{equation}
     \label{eq:Mfini}
     \q_{x,D}^{z,z',\alpha}({\mathcal M}_{\infty}^\alpha(\mathbb{R}^2)<\infty)
     =1,
     \qquad \mbox{Lebesgue-{\rm a.e.} }  x\in D .
 \end{equation}
  %Furthermore, for any $0<\varrho<2-\alpha$, 
 %\begin{equation}\label{eq:Mfini2}
 % \int_D {{\mathcal M}_{\infty}^\alpha(\d u)\over |u-x|^\varrho} <\infty, \qquad \mbox{$\q_{x,D}^{z,z',\alpha}$-a.s.}
%  \end{equation}
 \end{proposition}

\medskip

As a consequence of Theorem \ref{t:conv-p}, we have

\begin{corollary}
 \label{c:Minfini} 
 
 Under the assumption of Proposition \ref{p:Mfini}, we may define a random finite measure ${\mathfrak m}$ on the Borel sets   such that $\p^{z, z'}_D$-\hbox{\rm a.s.}, ${\mathfrak m}$ is the weak limit of ${\mathcal M}_{\D_{n}}^\alpha$. Moreover, for any rectangle $A$, ${\mathfrak m}(A)= \lim_{n\to\infty} {\mathcal M}_{\D_{n}}^\alpha(A)$, $\p^{z, z'}_D$-\hbox{\rm a.s.}

\end{corollary}

\noindent {\it Proof of Corollary \ref{c:Minfini}.} The argument is routine; we give the details for the sake of completeness. 

First, note that the sequence $({\mathcal M}_{\D_{n}}^\alpha)_{n\ge 0}$ is tight, all the measures being supported in the compact set $\overline{D}$; so we can extract a (random) subsequence, say $(n(k), \, k\ge 1)$, along which ${\mathcal M}_{\D_{n(k)}}^\alpha$ converges weakly to some finite random measure ${\mathfrak m}$. 

We may define $M_\infty^\alpha(A)$ such that  $\p^{z,z'}_{D}$-a.s.,  \eqref{MinfiniA} simultaneously holds for all rectangles $A$ with rational coordinates. For any $\varepsilon>0$ and rectangle $A$, let $A_{+, \varepsilon}$  be an open rectangle and $A_{-, \varepsilon}$ be a closed rectangle both with rational coordinates such that $A_{-, \varepsilon} \subset {\mathring A} \subset {\overline A} \subset A_{+, \varepsilon}$ and $\int_{  A_{+, \varepsilon}\backslash   A_{-, \varepsilon}} M_D(x, \, \alpha) \, \d x \le \varepsilon$. By the $L^1$-convergence in Theorem \ref{t:conv-p}, for any $\varepsilon>0$, $\e_D^{z, z'} \big( M_\infty^\alpha( A_{+, \varepsilon})-M_\infty^\alpha( A_{-, \varepsilon})\big)= \int_{  A_{+, \varepsilon}\backslash   A_{-, \varepsilon}} M_D(x, \alpha) \d x \le \varepsilon$. It implies that for any rectangle $A$, 
\begin{equation}\label{portmanteau}
\inf_{\varepsilon >0 } ( M_\infty^\alpha( A_{+, \varepsilon})-M_\infty^\alpha( A_{-, \varepsilon})) = 0 \qquad \p^{z, z'}_D\mbox{-a.s.}
\end{equation}

\noindent Consider the event $\mathscr{A}$ on which  \eqref{MinfiniA} and \eqref{portmanteau}  hold for all rectangles $A$ with rational coordinates.  By the   Portmanteau theorem, for any $\varepsilon>0$, we necessarily have ${\mathfrak m}( {\overline A} ) \le  {\mathfrak m}( A_{+, \varepsilon}) \le \liminf_{k\to\infty} {\mathcal M}_{\D_{n(k)}}^\alpha( A_{+, \varepsilon})= M_\infty^\alpha( A_{+, \varepsilon})$, and similarly, ${\mathfrak m}({\mathring A} ) \ge M_\infty^\alpha( A_{-, \varepsilon})$. Consequently, $$ {\mathfrak m}(\partial A) \le \inf_{\varepsilon >0 } ( M_\infty^\alpha( A_{+, \varepsilon})-M_\infty^\alpha( A_{-, \varepsilon}))=0$$ for any rectangle $A$ with rational coordinates on the event $\mathscr{A}$. The Portmanteau theorem implies that $ {\mathfrak m}(A)= \lim_{k\to\infty} {\mathcal M}_{\D_{n(k)}}^\alpha(A) = M_\infty^\alpha(A)$ for any such $A$, which by the monotone class theorem yields the uniqueness of the limit measure ${\mathfrak m}$ on $\mathscr{A}$ and proves the $\p^{z, z'}_D$-a.s.\ weak convergence of the sequence $({\mathcal M}_{\D_{n}}^\alpha)_{n\ge 0}$ to ${\mathfrak m}$. Finally, for  any rectangle $A$ (regardless of the rationality of the coordinates of $A$), equation \eqref{portmanteau} holds $\p^{z, z'}_D$-a.s., hence the same reasoning shows that ${\mathfrak m}(A)= \lim_{n\to\infty} {\mathcal M}_{\D_{n}}^\alpha(A)$, $\p^{z, z'}_D$-a.s.\hfill$\Box$

\begin{definition}[Definition of  $M_\infty^\alpha$] In the sequel, by a slight abuse of notation, we still denote  by $M_\infty^\alpha$ the finite random measure ${\mathfrak m}$  in Corollary \ref{c:Minfini}. 
\end{definition}

Theorem \ref{t:conv-p}  implies  in particular that $\e^{z,z'}_{D} ({\mathcal M}_{\infty}^\alpha(D))>0$. Furthermore,  $\p^{z,z'}_{D}$-a.s.,  ${\mathcal M}_{\infty}^\alpha$ is not  trivial, see Proposition \ref{p:loizeroun}.

The rest of this section is devoted to the proof of Proposition \ref{p:Mfini}. Firstly we present in Section \ref{sub:pre}  some preliminary estimates under $\q^{z,z',\alpha}_{x,D}$ by means of elementary properties of $\xi_{D}$;  then we give  the proof in Section \ref{sub:proof}. 

\subsection{Preliminary estimates}\label{sub:pre}

Let $x\in D$. Let $y$, $z\in  \overline D$ be distinct nice points, different from $x$. Recall from \eqref{def-xi} that $\xi_{D}(x,y,z) = \frac{2\pi H_D(x,y) H_D (x,z)}{H_D(y,z)}$. By properties of harmonic functions under conformal transformations, we see that $\xi_{D}$ is invariant under conformal transformations: if $D$ and $D'$ are two nice domains,  $z \neq z'$ nice points of $\overline D$, $\Psi:D\to D'$ a conformal transformation %which maps $D$ onto $D'$ 
such that $\Psi(z)$ and $\Psi(z')$ are nice points of $\overline{D'}$, then
$$
\xi_{D'}(\Psi(x),\Psi(z),\Psi(z')) = \xi_{D}(x,z,z').
$$

\begin{lemma}\label{l-xibound} Let  $D$ be a simply connected nice domain and $z' \in \partial D$ be a nice point.  

(i)  For any $x\in D$,
 and  any      nice point $z\in \partial D$ different from $z'$,  we have  \begin{equation}\label{xibound}
\xi_{D}(x,z,z')\le 2.
\end{equation}

(ii)  For any $\delta \in (0, 1)$, there exists some positive constant $c_4=c_4(\delta, D)$ such that for  any  $x, z \in D$ with $|x-z | > \delta$, we have \begin{equation}\label{xibound-2}
 \xi_{D}(x,z,z')\le c_4.
 \end{equation}
\end{lemma}

{\noindent\it Proof of Lemma \ref{l-xibound}.}  Observe that $\xi_{D}(x,z,z')=\xi_{{\cal B}(0,1)}(0,\Psi(z),\Psi(z'))$ for $\Psi$ a conformal map from $D$ to ${\cal B}(0,1)$ which maps $x$ to $0$.   It is well-known (see Lawler \cite{Lawler}, Chapter 2, formula (2.6)) that   for any $|a|<1$ and $|b|=1$,   \begin{equation}\label{lawler2.2.4} H_{{\cal B}(0,1)}(a,b)= { 1-|a|^2\over 2 \pi \,  |a-b|^2}. \end{equation}
 
 It follows that for $|b'|=|b|=1$, $H_{{\cal B}(0,1)}(b',b)=   \frac{1}{ \pi \, |b-b'|^2}$ (see Lawler \cite{Lawler}, Chapter 5, Example 5.6).

 If $z \in \partial D$, then $|\Psi(z)|=1$ and      $\xi_{{\cal B}(0,1)}(0,\Psi(z),\Psi(z'))= \frac12 |\Psi(z)- \Psi(z')|^2  \le 2$  yielding \eqref{xibound}.

It remains to treat  the case $z \in D$ with $|z-x|>\delta$.  Write $a:= \Psi(z)$ for notational brevity.  Using the elementary    fact that for any $|a| <1$, $H_{{\cal B}(0,1)}(0,a)=\frac1{2\pi} \log \frac1{|a|}$,  we have  \begin{equation}\label{xi-bound-inter}\xi_{{\cal B}(0,1)}(0,\Psi(z),\Psi(z'))
=
  \big(\log \frac1{|a|} \big) \frac{|a- \Psi(z')|^2}{1- |a|^2}
 \le
   \frac{4\log \frac1{|a|}}{1-|a|} , 
  \end{equation}
  
\noindent where  we used $|a-\Psi(z')|\le 2$ by the triangular inequality and $1-|a|^2\ge 1-|a|$. %by using the triangular inequality: $|a- \Psi(z')| \ge | \Psi(z')|- | a|= 1-| a|$.    

Now we    estimate $|a|$ from below.     Observe that  $|z-x|= | \Psi^{-1}(a)- \Psi^{-1}(0)| \le |a| \, \sup_{0 \le t \le 1} | (\Psi^{-1})'(at)| .$ 
Applying Corollary 3.19 (Lawler \cite{Lawler}, Chapter 3) to the conformal transformation $\Psi^{-1}$ which maps  ${\cal B}(0,1)$ to $D$, we get that  for any $0 \le t \le 1$, $| (\Psi^{-1})'(at)|  \le  4 \frac{d(\Psi^{-1}(at), \partial D)}{1- |a t|} \le 4 \frac{\mathrm{diam}(D)}{1-|a|}$, where $\mathrm{diam}(D)$ denotes the diameter of $D$. Then   $|z-x| \le 4\, \mathrm{diam}(D) \frac{|a|}{1-|a|}$.  Since $|z-x| > \delta$ by assumption, we deduce that $|a| \ge \eta$, with $\eta:=\frac\delta{4\mathrm{diam}(D) +\delta}$.  Going back to  \eqref{xi-bound-inter}, we get that $\xi_{D}(x,z,z')=\xi_{{\cal B}(0,1)}(0,\Psi(z),\Psi(z')) \le    \sup_{\eta \le r <1} \frac{4\log\frac1{r}}{1-r}=:c_4,$ proving  \eqref{xibound-2}.\hfill$\Box$

\medskip 
 
{F}rom \eqref{lawler2.2.4}, we have the following lemma.

\begin{lemma}
 \label{l:poissonkernel}
 
 For each $\varepsilon>0$, there exists $\eta>0$ such that $| {H_{{\cal B}(0,1)} (a,b) \over H_{{\cal B}(0,1)} (0,b)} -1|<\varepsilon$ for any $|b|=1$ and $|a|<\eta$.

\end{lemma}

In view of the proof of Proposition \ref{p:Mfini}, we now study loops from $x$ to $x$ under $\q_{x,D}^{z,z',\alpha}$. For any $r>0$, we denote by
\begin{eqnarray}
    \label{Nxrbeta}
    N(x,r) 
 &:=& \#\{\mbox{loops from $x$ to $x$ which hit ${\cal C}(x,r)$}\} ,
    \\
    \label{Sigmaxr}
    \varUpsilon(x,r) 
 &:=& \sum_{\eee \in {\mathcal E}_{{\cal B}(x,r)}} \xi_{{\cal B}(x,r)}(x,\eee_{g},\eee_{d}) \, 1_{\{x\notin \eee\}}, 
\end{eqnarray}
 
 \noindent where ${\mathcal E}_{{\cal B}(x,\, r)}$ is the set of excursions inside ${\cal B}(x,\, r)$ as defined in Notation \ref{notation:tribu}, and $\xi$   is defined in \eqref{def-xi}.
 
Let $z\neq z'$ be distinct nice points of $\overline{D}$, different from $x$. Let $r_0\in (0, \, 1)$ be such that ${\cal B}(x,\, r_0) \subset D$ and $z$, $z'\notin \overline{{\cal B}(x,\, r_0)}$. Let $r_k:= \frac{r_0}{2^k}$ for $k\ge 0$. 
% \noindent where, for $\eee\in {\mathcal E}_{x,{\cal B}(x,r),{\cal C}(x,r')}$ $\xi_{r,r'}(x,\eee_{g},\eee_{d})$ is equal to
% \begin{itemize}
% \item $\xi_{{\cal B}(x,r)}$ if $\eee_{g}$ and $\eee_{d}$ belong to ${\cal C}(x,r)$;
% \item $C_{{\cal B}(x,r),{\cal B}(x,r')}(x) \big(1- {H_{D_{1}}(\eee_{g},x) \over H_{D}(\eee_{g},x)} \big)$ if $\eee_{d}=x$
% \item $C_{{\cal B}(x,r),{\cal B}(x,r')}(x) \big(1- {H_{D_{1}}(x,\eee_{d}) \over H_{D}(x,\eee_{d})} \big)$ if $\eee_{g}=x$
% \item $C_{{\cal B}(x,r),{\cal B}(x,r')}(x)$ if $\eee_{g}=\eee_{d}=x$.
% \end{itemize}
 
\begin{lemma}
 \label{l:spine}  
 
  Let $\alpha \ge 0$. As $k\to \infty$,

{\rm (i)} $N(x,r_k) \sim \alpha  \log \frac{1}{r_k}$, $\quad \q_{x,D}^{z,z',\alpha} \hbox{\rm -a.s.}$;

{\rm (ii)} $\varUpsilon(x,r_{k}) \sim \alpha (\log \frac{1}{r_k})^2$, $\quad \q_{x,D}^{z,z',\alpha} \hbox{\rm -a.s.}$

% \item $N_{x}^{(\beta)}(r_{k+1})-N_{x}^{(\beta)}(r_{k}) = o_{k}(1) \log(1/r_{k}) (\log(\log(1/r_{k})))^{-1/2}$ as $k\to +\infty$;
 %\item We have that $\sum_{(g,d)_{{\cal C}(x,r_{k})}} \xi_{{\cal C}(x,r_{k})}(x,{\mathfrak e}_{g},{\mathfrak e}_{d})  \sim \alpha  \,   2\pi  (\log r_{k})^2$ as $k\to+\infty$;

\end{lemma}

 {\noindent\it Proof of Lemma \ref{l:spine}}. (i) 
We start by mentioning a simple fact on the concentration of a Poisson variable: Let $N_\lambda$ be a Poisson random variable with parameter $\lambda >0$.   By the standard  large deviation principle,  for any $\varepsilon>0$, there exist some $\delta>0$ and $\lambda_0>0$ such that for all $\lambda>\lambda_0$, 
\begin{equation}
    \label{deviation-poisson} 
    \p\Big( | N_\lambda - \lambda | \ge \varepsilon \lambda \Big) \le \ee^{- \delta \, \lambda}.
\end{equation}

Under $\q^{z,z',\alpha}_{x, D}$, the number of loops from $x$ to $x$ that hit  ${\cal C}(x,r_{k})$ is a Poisson variable of parameter $\alpha\, \nu_D(x,x)(T_{{\cal C}(x,r_{k})}< T_x)$ (see Section \ref{subs:Q}). By Lemma \ref{l:ito-muxx} (ii), $ \nu_D(x,x)(T_{{\cal C}(x,r_{k})}< T_x) = \int_{{\cal C}(x,r_{k})} G_D(x, y) H_{{\cal B}(x,r_{k})}(x, y) \d y  \sim   \log (1/r_k)$ as $k \to \infty$.  By \eqref{deviation-poisson} and the Borel--Cantelli lemma, $ N(x,r_k) \sim \alpha  \log(1/r_{k})$ as $k \to \infty$,  $\q_{x,D}^{z,z',\alpha}$-a.s. This proves (i). 

(ii) For each $r>0$, the contribution to the right-hand side of \eqref{Sigmaxr} comes from the paths of the excursion from $z$ to $x$, of the $N(x,r)$ loops which hit ${\cal C}(x,r)$, and of the excursion from $x$ to $z'$. Consequently for any $k\ge1$,
$$
\varUpsilon(x,r_k)
=
\varUpsilon_{begin} + \sum_{j=1}^{N(x, r_k)}  \varUpsilon_j + \varUpsilon_{end} \, ,
$$

\noindent  where under $\q_{x,D}^{z,z',\alpha}$,  $\varUpsilon_{begin}$, $\varUpsilon_{end}$,    $ \varUpsilon_1, \varUpsilon_2, ... , $   are mutually independent (and independent of   $N(x, r_k)$)   such that 
$\varUpsilon_{begin}$ and $\varUpsilon_{end}$ are distributed as $\varUpsilon(x, r_k)$ under $\p_D^{z, x}$ and $\p_D^{x, z'}$ respectively, and  for any  $j\ge1$, $\varUpsilon_j$ has the same distribution as that  of  $\varUpsilon(x, r_k)$ under  $\nu_D(x,x) \big(\bullet \, |\, T_{{\cal C}(x,r_k)} < T_x\big)$. [For the sake of presentation, we have introduced $\varUpsilon_j$ for all $j\ge 1$; we have also omitted the dependence on $k$ in the notation $\varUpsilon_{begin}, \varUpsilon_{end}$ and $\varUpsilon_j$, $j\ge 1$.]

We claim that \begin{equation}\label{varUpsilon-moment1}
\q_{x,D}^{z,z',\alpha}(\varUpsilon_{begin})= \q_{x,D}^{z,z',\alpha}(\varUpsilon_{end})= \q_{x,D}^{z,z',\alpha}(\varUpsilon_1)= C_{D,{\cal B}(x,r_{k})}(x) .
\end{equation}

In fact, notice  that $H_{{\cal B}(x,r_{k})}(x, z)=0$ as well as $H_{{\cal B}(x,r_{k})}(x, z')=0$. Applying \eqref{eq:radonzx} to $L=1$  gives that  $\q_{x,D}^{z,z',\alpha}(\varUpsilon_{begin})=   C_{D,{\cal B}(x,r_{k})}(x) $. In the definition of $\varUpsilon(x,r)$ in \eqref{Sigmaxr}, $\xi_{{\cal B}(x,r)}(x, \, \eee_{g},\, \eee_{d}) = \xi_{{\cal B}(x,r)}(x, \, \eee_d,\, \eee_g)$, so by the time-reversal property \eqref{timereversal} for the path from $x$ to $z'$ and another application of \eqref{eq:radonzx}  to $L=1$, we get $\q_{x,D}^{z,z',\alpha}(\varUpsilon_{end})=   C_{D,{\cal B}(x,r_{k})}(x) $.  Finally,   recall     that     $\nu_{D}(x,x)(T_{{\cal C}(x,r_{k})} < T_{x})= C_{D,{\cal B}(x,r_{k})}(x)$  by definition (see  Lemma \ref{l:ito-muxx}). Applying   \eqref{eq:radonxx} to $L=1$  gives that $\q_{x,D}^{z,z',\alpha}(\varUpsilon_1)= C_{D,{\cal B}(x,r_{k})}(x) $ and completes  the justification of \eqref{varUpsilon-moment1}.

%  We see from   \eqref{eq:radonzx} and \eqref{eq:radonxx}   with $L=1$  that they all have first moment equal to $C_{D,{\cal B}(x,r_{k})}(x) $ (remark that $H_{{\cal B}(x,r_{k})}(x, z)=0$ for all large $k$;   we use \eqref{eq:radonzx}  for the path from $x$ to $z'$ by time-reversal; we notice that  for $\varUpsilon_{j}, j\ge 2$, in \eqref{eq:radonxx} we have to condition $\nu_{D}(x,x)$ on the event $\{T_{{\cal C}(x,r_{k})} < T_{x}\}$ which has measure $C_{D,{\cal B}(x,r_{k})}(x)$). 

%We decompose the trace of the Brownian motion on ${\cal C}(x,r_{k})$ as the trace of the excursion from $z$ to $x$, of the excursion from $x$ to $z'$ and of the $N(x,r_{k})$ loops which hit ${\cal C}(x,r_{k})$. We call $\varUpsilon_{0}$, $\varUpsilon_{1}$ and $\varUpsilon_{j},\, 2\le j\le N(x,r_{k})+1$  the contribution of each of these excursions to $\varUpsilon(x,r_{k})$ (we omit the dependence on $k$ in the notation $\varUpsilon_{j}$).  These random variables are independent.   We see from   \eqref{eq:radonzx} and \eqref{eq:radonxx}   with $L=1$  that they all have first moment equal to $C_{D,{\cal B}(x,r_{k})}(x) $ (remark that $H_{{\cal B}(x,r_{k})}(x, z)=0$ for all large $k$;   we use \eqref{eq:radonzx}  for the path from $x$ to $z'$ by time-reversal; we notice that  for $\varUpsilon_{j}, j\ge 2$, in \eqref{eq:radonxx} we have to condition $\nu_{D}(x,x)$ on the event $\{T_{{\cal C}(x,r_{k})} < T_{x}\}$ which has measure $C_{D,{\cal B}(x,r_{k})}(x)$). 

By  Lemma \ref{l:ito-muxx} (ii),  $C_{D,{\cal B}(x,r_{k})}(x) = \int_{{\cal C}(x, r_k)} G_D(x, y) H_{{\cal B}(x,r_{k})}(x, y ) \d y \sim \log(1/r_{k})$ as $k \to \infty$. [A fact already used in the proof of (i)]. % It  yields that    $ \q_{x,D}^{z,z',\alpha}( \varUpsilon_{j}) \sim \log(1/r_{k}) $ uniformly in $j$. 
%Note that $\varUpsilon(x,r_{k}) = \sum_{j=0}^{N(x,r_{k})+1}  \varUpsilon_{j}$. 

 From (i),    the statement  (ii) immediately follows once we have shown  that for  any  deterministic sequence $n_k$ such  that $\liminf_{k\to \infty} \frac{n_k}{k} >0$,  \begin{equation}\label{BorelCantelli3}{1\over n_k^2}  \Big[ \varUpsilon_{begin} + \varUpsilon_{end}+ \sum_{j=1}^{n_k}  \varUpsilon_{j}-  (n_k+2) \, C_{D,{\cal B}(x,r_{k})}(x)    \Big]  \, \to\, 0, \qquad \mbox{$\q_{x,D}^{z,z',\alpha}$-a.s.}\end{equation}

  To get \eqref{BorelCantelli3},  we shall use  the following     inequality  (see Petrov \cite{petrov}, Theorem 2.10):  There exists some constant $c_1>0$ such that for  any  sequence of   independent  real-valued integrable  random variables $(\eta_i)_{i\ge1}$,    \begin{equation} \label{sumiid}
  \e \Big| \sum_{i=1}^n (\eta_i - \e(\eta_i))\Big|^3 
  \le
  c_1\, n^{1/2}\, \sum_{i=1}^n \e[ |\eta_i|^3], \qquad \forall\, n\ge 1. 
  \end{equation}
 
%  \begin{comment}
%We claim the existence of some constant $c_2>0$ such that for all large $k$, uniformly in $j\ge 0$,  $$\q_{x,D}^{z,z',\alpha}(\varUpsilon_{j}^3) \le c_2 \, \Big( C_{D,{\cal B}(x,r_{k})}(x)\Big)^3.$$
% 
% %\sum_{  {\underline{\eee}_{{\cal B}(x, r_k) }}  }\xi_{{\cal B}(x, r_k)}(x,{\mathfrak e}_{g}, {\mathfrak e}_{d}) 1__{\{\eee_{g}^1\neq z, \eee_{d}^L\neq x\}}
% 
% \noindent In fact,  notice that $\xi_{D_{1}}(x,\eee_{g},\eee_{d})$ is bounded by $4$ (see \eqref{xibound}).  It follows that  \begin{eqnarray*}
% &&
% \q_{x,D}^{z,z',\alpha}(\varUpsilon_{0}^3) 
% \\
% &=&
% \e^{z, x}_D \Big( \sum_{{\mathcal E}_{{\cal B}(x, r_k) }}  \xi_{{\cal B}(x, r_k)}(x,{\mathfrak e}_{g}, {\mathfrak e}_{d})  \Big)^3
% \\
% &\le&
% 6\,  \e^{z, x}_D \Big( 4^2 \sum_{{\underline{\eee}_{{\cal B}(x, r_k) }}  }\xi_{{\cal B}(x, r_k)}(x,{\mathfrak e}_{g}, {\mathfrak e}_{d})+
%  4 \,   \sum_{{\underline{\eee}^2_{{\cal B}(x, r_k) }}  } \prod_{j=1}^2 \xi_{{\cal B}(x, r_k)}(x,{\mathfrak e}^j_{g}, {\mathfrak e}^j_{d})
%  +   \sum_{{\underline{\eee}^3_{{\cal B}(x, r_k) }}  } \prod_{j=1}^3 \xi_{{\cal B}(x, r_k)}(x,{\mathfrak e}^j_{g}, {\mathfrak e}^j_{d})\Big)
% \end{eqnarray*}
% \end{comment}
 
Using the fact that  $\xi_{{\cal B}(x, r_k)}(x,\eee_{g},\eee_{d})$ is bounded by   $2$  (see \eqref{xibound}), we deduce from the definition of $ \varUpsilon(x, r_k)$ in \eqref{Sigmaxr} that for some numerical constant $c_2$ [we may take $c_2=  24$],  $$ \varUpsilon(x, r_k)^3 \le c_2 \,  \sum_{L=1}^3\, \sum_{ {\underline{\eee}_{{\cal B}(x, r_k)}^L}} \prod_{j=1}^{L}\xi_{{\cal B}(x, r_k)}(x,{\mathfrak e}_{g}^j, {\mathfrak e}_{d}^{j}) 1_{\{\eee_{g}^1\neq x, \eee_{d}^L\neq x\}}, $$ 
  
  \noindent where as in the previous section,  $\sum_{{\underline \eee}_{{\cal B}(x, r_k)}^L}$ is a short way to denote sum over all ordered (distinct) excursions $(\eee^1,\ldots,\eee^L)   \in ({\mathcal E}_{{\cal B}(x, r_k)})^L$. 
  
{F}rom  the  equations \eqref{eq:radonzx} and \eqref{eq:radonxx}   with $L \in \{1, 2, 3\}$, $D_1={\cal B}(x, r_k)$  there, we see that there exists some positive constant $c_3 = c_3(x)$ such that    the third moments of $\varUpsilon_{begin}$, of $\varUpsilon_{end}$ and of $  \varUpsilon_1$   are  less than $c_3 \, (\log(1/r_{k}))^3$.
 
%  From  equations \eqref{eq:radonzx} and \eqref{eq:radonxx}   with $L \in \{1, 2, 3\}$, $D_1={\cal B}(x, r_k)$  there,  and  using the fact that $\xi_{D_{1}}(x,\eee_{g},\eee_{d})$ is bounded by $4$ (see \eqref{xibound}),  we see that for all large $k$, uniformly in $j\ge 0$,    $\q_{x,D}^{z,z',\alpha}(\varUpsilon_{j}^3) \le c_2 \, \big( C_{D,{\cal B}(x,r_{k})}(x)\big) ^3 \le c_3 \, (\log(1/r_{k}))^3$.
  
  Recalling \eqref{varUpsilon-moment1}.   It follows from   \eqref{sumiid}   that 
\begin{eqnarray*}
&& 
	\q_{x,D}^{z,z',\alpha} \Big[\Big| {1\over n_k^2}  \big( \varUpsilon_{begin} + \varUpsilon_{end}+ \sum_{j=1}^{n_k}  \varUpsilon_{j}-  (n_k+2) \, C_{D,	{\cal B}(x,r_{k})}(x)    \big) \Big|^3\Big] \\
&\le&
	c_1 c_3\, \frac{(n_k+2)^{3/2}}{n_k^6} \,  (\log \frac{1}{r_k})^{3}, 
\end{eqnarray*}

\noindent which is summable in $k$ thanks to the assumption on $n_k$: $\liminf_{k\to\infty} \frac{n_k}{k}>0$.  The Borel--Cantelli lemma yields \eqref{BorelCantelli3}   and  completes the proof of (ii).  \hfill $\Box$

\subsection{Proof of Proposition \ref{p:Mfini}}\label{sub:proof}

We fix $D$, $z\neq z'$ nice points of $\overline D$.  Let $0\le \alpha< 2$ and $0 \le \varrho<  2-\alpha$. Take $\gamma=\gamma(\varrho)$ and $\varepsilon=\varepsilon(\varrho)$ such that $\gamma>\alpha$, $\varepsilon>0$ and 
\begin{equation}\label{gamma}
2(1+\varepsilon)\sqrt{\gamma\alpha} -\alpha +\varrho<2.
\end{equation}

\noindent Let $\eta \in (0, 1) $ be the constant in Lemma \ref{l:poissonkernel} associated with our choice of $\varepsilon$. Let $K\ge 5$ be such that $2^{-K}<\frac{\eta}{16}$. Constants $c_5$, $c_6$, $\ldots$ in the proof can depend on $D$, $z$, $z'$, $\alpha$, $\varrho$, $\gamma$, $\varepsilon$, $K$ even if not specified.

Let $x\in D$ be a suitable point, meaning, as before, that $x$ is in some square $D_{n}\in \D_{n}$ at any level $n$, and is different from $z$ and $z'$ (recall that the set of suitable points has full Lebesgue measure). Let $r_{0}\in (0, \, 1)$ and $r_{k}= \frac{r_0}{2^k}$, let 
\begin{eqnarray}
 {E}^{(x)}(r_{0},\gamma)    
 &:=&
 \big\{d(x,\partial D\cup\{z,z'\})>r_{0}\big\} \cap  \nonumber
 \\
 && \cap_{k\ge 0} \big\{N(x,r_{k}) \le \gamma \log \frac{1}{r_k},\, \varUpsilon(x,r_k) \le \gamma  (\log \frac{1}{r_k})^2 \big\}.  \label{goodevent}
\end{eqnarray}
% \begin{align}\label{goodevent}
% \begin{split}
 %{E}^{(x)}(r_{0},\gamma) &:= \\
% &\big\{d(x,\partial D\cup\{z,z'\})>r_{0}\big\}
% \bigcap_{k\ge 0} \big\{N(x,r_{k}) \le \gamma \ln(1/r_{k}),\, \varUpsilon(x,r_k) \le \gamma  \log(r_{k})^2 \big\}.
% \end{split}
% \end{align}
%Recall that a point $x$ is suitable if for any $n\ge0$, $x$ belongs to some open square of mesh $r_n$, denoted by  $D_n^{(x)}$.  %Recall that  the set of suitable points has the full Lebesgue measure. 

\noindent We assume, for the moment, that for all suitable $x\in D$,
\begin{equation}
      \limsup_{n\to+\infty}  1_{\{d(x,\partial D_{n}^{(x)})\ge {2^{-n} \over 4} \}} \e_{\q^{z,z',\alpha}_{x,D}}\left[ \int_{D} \, \frac{{\mathcal M}_{\D_{n}}^\alpha(\d u)}{|u-x|^\varrho},\, {E}^{(x)}(r_{0},\gamma)  \right]  \le c_5(r_{0}), 
    \label{c_5(N)}
\end{equation}

\noindent where $c_5(r_{0})$  is some positive constant depending on  $r_{0}$. We claim that Proposition \ref{p:Mfini} will follow from Lemma \ref{l:spine} and (the case $\varrho=0$ of)  \eqref{c_5(N)}. In fact,  for any  $c>0$, we deduce from  \eqref{c_5(N)} (with $\varrho=0$)  that
%\begin{eqnarray*}
%&&\liminf_{n\to+\infty}  \e_{\q^{z,z',\alpha}_{x,D}}\left[ \min\left(K, \int_{D} \, \frac{{\mathcal M}_{\D_{n}}^\alpha(\d u)}{|u-x|^\varrho} \right),\, \mathcal{E}^{(x)}_{N}(\gamma_{1},\varepsilon_{1})  \right] 
%\\
%&\le&
% c_5(N,\varepsilon) + K\liminf_{n\to+\infty}  1_{\{d(x,\partial D_{n}^{(x)})< {r_{n} \over 4} \}}
% \\
% &=&c_5(N,\varepsilon),
 %\end{eqnarray*}
\begin{eqnarray*}
&& 
 	\liminf_{n\to+\infty}  \e_{\q^{z,z',\alpha}_{x,D}}\left[ \min\left(c,  {\mathcal M}_{\D_{n}}^\alpha( D)  \right),\, {E}^{(x)}(r_{0}, \gamma)  \right] \\
&\le&
	c_5(r_{0}) + c\liminf_{n\to+\infty}  1_{\{d(x,\partial D_{n}^{(x)})< {r_{n} \over 4} \}} \\
 &=&
 	c_5(r_{0}),
 \end{eqnarray*}

\noindent for Lebesgue-a.e.\ $x\in D$.  By Fatou's lemma, %\footnote{For any nonnegative measurable function $f$,  we have that  $\liminf_{n\to \infty}  \int  f(u) {\mathcal M}_{\D_{n}}^\alpha(\d u)  \ge \int  f(u) {\mathcal M}_{\infty}^\alpha(\d u)$, $\q^z_{x,D}$-a.s. In fact,   let $(f_k)$ be a non-decreasing sequence of nonnegative step functions such that $f=\lim_{k\to \infty}\uparrow  f_k$.  For any fixed $k\ge1$,     $\int  f_k(u) {\mathcal M}_{\D_{n}}^\alpha(\d u) \to \int  f_k(u) {\mathcal M}_{\infty}^\alpha(\d u)  $, a.s.  We conclude by the monotone convergence theorem and the fact that $\int  f(u) {\mathcal M}_{\D_{n}}^\alpha(\d u)  \ge \int  f_k(u) {\mathcal M}_{\D_{n}}^\alpha(\d u)$.}
  it gives that for any $c>0$, 
$$
 \e_{\q^{z,z',\alpha}_{x,D}}\left[\min\left(c,  {\mathcal M}_\infty^\alpha( D)  \right),\, {E}^{(x)}(r_{0},\gamma)  \right] \le c_5(r_{0}), 
$$
then, by monotone convergence,
\begin{equation} \label{cnuxa}
 \e_{\q^{z,z',\alpha}_{x,D}}\left[ {\mathcal M}_\infty^\alpha( D)  ,\, {E}^{(x)}(r_{0}, \gamma)  \right] \le c_5(r_{0}).
\end{equation}

\noindent We deduce that  for Lebesgue-a.e.\ $x\in D$, $\q^{z,z',\alpha}_{x,D}$-a.s., $  {\mathcal M}_\infty^\alpha(D) <\infty$ on the event ${E}^{(x)}(r_{0},\gamma)$. By Lemma \ref{l:spine}, we have for any $x\in D\backslash\{z,z'\}$, 
\begin{equation}
    \q_{x,D}^{z,z',\alpha} \Big( \bigcup_{\ell=1}^\infty {E}^{(x)}(2^{-\ell}, \, \gamma) \Big)
    =
    1 \, .
    \label{unionen}
\end{equation}

\noindent Hence we get \eqref{eq:Mfini}. %As $\varepsilon$ can be chosen as small as desired,  Proposition \ref{p:Mfini} follows. 
[The case $0< \varrho < 2-\alpha$ in \eqref{c_5(N)} will be used in Section \ref{s:thickpoints}.]

It remains to prove \eqref{c_5(N)}. We fix $r_{0}>0$ and suitable $x\in D$ such that $d(x,\partial D\cup\{z,z'\})>r_{0}$. We write ${\mathcal B}_{0}:=D$ and ${\mathcal B}_{k}:={\mathcal B}(x,r_{k})$ for any $k\ge 1$. We distinguish three possible situations: (i) $|u-x|>r_{K}$ with $u\notin D_{n}^{(x)}$, (ii) $|u-x|<r_{K}$ with $u\notin D_{n}^{(x)}$, and (iii) $u\in D_{n}^{(x)}$. We suppose that $d(x,\partial D_{n}^{(x)})\ge {2^{-n} \over 4}$, which means that $x$ is not too close to the boundary of $D_{n}^{(x)}$.

\paragraph{First case: $|u-x|>r_K$ with $u \notin D_n^{(x)}$} Since we are going to integrate over $u$ with respect to the Lebesgue measure, we suppose, without loss of generality, that $u$ is a suitable point. With the notation of Proposition \ref{p:martingale2}, observe that $\widetilde M_{x,D_{n}^{(u)}}(u,\alpha)=M_{D_{n}^{(u)}}(u,\alpha)$ since the excursions inside $D_{n}^{(u)}$ cannot hit $x$. Recall the definition of $\varUpsilon_{x,D_{1},S,D_{2}}(u)$ in \eqref{defSigma}. We are going to take $D_{1}=D$, $S={\cal B}_{K}$ and $D_{2} = D_{n}^{(u)}$. Then $\varUpsilon_{x,D,{\cal B}_{K},D_{n}^{(u)}} (u)$ is (possibly) contributed by excursions from $z$ to $x$, from $x$ to $z'$ and by the loops at $x$ in $D$ which hit ${\cal C}(x,r_{K})$. By definition,
% of $\xi_{x,D,{\cal B}_{K},D_{n}^{(u)} }$,   
\begin{eqnarray*}
 && 
 	\varUpsilon_{x,D,{\cal B}_{K},D_{n}^{(u)}} (u) \\
&=&
	\frac{2\pi}{H_D(z, x)} (H_D(u, z)- H_{D_n^{(u)}}(u, z)) (H_D(u, x)- H_{D_n^{(u)}}(u, x))  \\
&&
	+ \frac{2\pi}{H_D(x, z')} (H_D(u, z')- H_{D_n^{(u)}}(u, z')) (H_D(u, x)- H_{D_n^{(u)}}(u, x))
\\
&& 
	+  \frac{(2\pi)^2}{C_{D,{\cal B}_{K}}(x)} (H_D(u, x)- H_{D_n^{(u)}}(u, x))^2 \, N(x,r_{K}), 
\end{eqnarray*}

\noindent by recalling that $N(x,r_{K})$ denotes the number of  loops at $x$ which hit ${\cal C}(x, r_K)$. This implies that
\begin{equation}
    \label{cxr}
    \varUpsilon_{x,D,{\cal B}_{K},D_{n}^{(u)}} (u)
    \le
    \xi_{D}(u,x,z)+\xi_{D}(u,x,z') + {G_{D}(u,x)^2\over C_{D,{\cal B}_{K}}(x)}N(x,r_{K})
    =:
    c_{\eqref{cxr}} \, ,
\end{equation}

\noindent recalling that $G_{D}(u,x) = 2\pi H_D(u, x)$. 

By applying Proposition \ref{p:martingale2} to $D_{1}=D$, $S={\cal B}_{K}$ and $D_{2} = D_{n}^{(u)}$, we get that
\begin{eqnarray}
    \e_{\q^{z,z',\alpha}_{x,D}} [ M_{D_{n}^{(u)}}(u,\alpha)\, | \,  N(x,r_{K}) ]
 &=& \e_{\q^{z,z',\alpha}_{x,D}} [ \widetilde M_{x,D_{n}^{(u)}}(u,\alpha) \, | \, N(x,r_{K}) ]
   \nonumber \\
 &\le& 
c_{\eqref{cxr}} \, \ee^{-\alpha\,  C_{D}(u)}\, \ee^{2\sqrt{\alpha \, c_{\eqref{cxr}}}}. \label{condNxrK}
\end{eqnarray}

Let us control the three terms in $c_{\eqref{cxr}}$.  For the last term in $c_{\eqref{cxr}}$, we remark that since $|u-x| >r_K$, we have, by Lawler \cite{Lawler} Proposition 2.36, $G_D(u, x)\le \log \mathrm{diam}(D) + \log \frac1{r_K}$ (where $\mathrm{diam}(D)$ denotes, as before, the diameter of $D$). For any $|y-x| = r_K$, $G_D(x, y) \ge   \log d (x, \partial D) + \log \frac1{r_K}$ (again by Lawler \cite{Lawler} Proposition 2.36). It follows that $C_{D,{\cal B}_K}(x) = \int_{{\cal C}(x, r_K)} G_D(x, y) H_{{\cal B}(x,r_{K})}(x, y ) \d y \ge  \log(r_0/r_{K})$. On the event ${E}^{(x)}(r_{0},\gamma)$, we have $N(x,r_{K})\le \gamma (\log r_{K})^2$, hence for some constant $c_6 = c_6(r_0, \, K, \, D, \, \gamma)$,
$$
{G_{D}(u,x)^2\over C_{D,{\cal B}_{K}}(x)}N(x,r_{K}) \le c_6 .
$$

To control $\xi_{D}(u,x,z)+\xi_{D}(u,x,z')$ in $c_{\eqref{cxr}}$, we discuss separately two cases: If $z \in \partial D$, then we apply \eqref{xibound-2} to see that $\xi_{D}(u,x,z)\le c_4$. If $z \in D$, $H_D(u,\, x)$ is bounded as seen in the previous paragraph (recalling that $|u-x|>r_K$), $H_D(x, \, z) \ge \inf_{y: \, d(y, \partial D \cup \{ z\}) \ge r_0} H_D(y, \, z) =: c_7 (r_0, \, z, \, D)>0$, hence $\xi_{D}(u,x,z)=2 \pi \frac{H_D(u, x)}{H_D(x, z)} H_D(u, z) \le 2\pi c_7\,  H_D(u, z)$. We get that for $z\in \overline D$,  $\xi_D(u, x, z) \le c_4 1_{\{z\in \partial D\}} + 2\pi c_7 H_D(u,z)1_{\{z\in   D\}} =:f_z(u)$. A similar bound holds for $\xi_{D}(u,\, x,\, z')$. Therefore we have shown that on the event ${E}^{(x)}(r_{0},\gamma)$, %for some constant $c_8 =c_8 (r_0, \, z, \, z', \, K, \, D, \, \gamma)$, 
$$
c_{\eqref{cxr}}
\le
c_6 + f_z(u)+ f_{z'}(u)=: f_{z, z'}(u) .
$$

Going back to \eqref{condNxrK}, we see that  
\begin{eqnarray*}
&& 	
	\e_{\q^{z,z',\alpha}_{x,D}}\left[ M_{D_{n}^{(u)}}(u,\alpha) ,\, {E}^{(x)}(r_{0},\gamma)  \right] \\
& \le&
	 \e_{\q^{z,z',\alpha}_{x,D}}\left[ M_{D_{n}^{(u)}}(u,\alpha) ,\, N(x,r_{K})\le \gamma (\log r_{K})^2 \right] 
 \\
 &\le&
 	f_{z, z'} (u) \,  \ee^{-\alpha\,  C_{D}(u)}\, \ee^{2\sqrt{\alpha  f_{z, z'} (u)}},
 \end{eqnarray*}
 
 \noindent  which implies that \begin{eqnarray} 
 &&  \e_{\q^{z,z',\alpha}_{x,D}}\left[ \int_{|u-x|>r_K} \, \frac{{\mathcal M}_{\D_{n}}^\alpha(\d u)}{|u-x|^\varrho},\, {E}^{(x)}(r_{0},\gamma)  \right] 
\nonumber
\\
& \le&
 \int_{|u-x|>r_K} \, \frac{ \d u}{|u-x|^\varrho}  f_{z, z'} (u) \,  \ee^{-\alpha\,  C_{D}(u)}\, \ee^{2\sqrt{\alpha  f_{z, z'} (u)}}
 \nonumber
 \\
 &\le&
 c_9, \label{Mxr}
 \end{eqnarray}

 \noindent with $c_9 :=   r_K^{-\varrho}\, \int_D    f_{z, z'} (u) \,  \ee^{-\alpha\,  C_{D}(u)}\, \ee^{2\sqrt{\alpha  f_{z, z'} (u)}}  \d u < \infty$.

  %for any $|u-x|>r_{K}$ with $u\notin D_{n}^{(x)}$, $u\neq z$ and $u\neq z'$,  
%\begin{equation}\label{Mxr}
%\e_{\q^{z,z',\alpha}_{x,D}} \Big[ M_{D_{n}^{(u)}}(u,\alpha),\, {E}^{(x)}(r_{0},\gamma)\Big]
%\le 
%c_5(r_{0}).
%\end{equation}

\paragraph{Second case: $r_{k+K+1}\le |u-x|<r_{k+K}$ for some $k\ge 0$, with $u\notin D_{n}^{(x)}$} We still have $\widetilde M_{x,D_n^{(u)}}(u,\alpha)=M_{D_n^{(u)}}(u,\alpha)$. We claim
\begin{equation}
    \overline{D_n^{(u)}} \subset {\cal B}_k \, ,
    \qquad
    \partial D_{n}^{(u)} \subset {\cal B}_{k+K+5}^c \, .
    \label{k+K+5}
\end{equation}

To see why $\overline{D_n^{(u)}} \subset {\cal B}_k$ holds: It suffices to check $|u-x| + \sqrt{2} \times 2^{-n} < r_k$, which is easy. Indeed, $|u-x| <r_{k+K} < \frac12 \, r_k$ (since $K>1$), and since $d(x,\partial D_{n}^{(x)})\ge \frac14 \, 2^{-n}$, we also have $\sqrt{2} \times 2^{-n} < \sqrt{2} \times 4 d(x,\partial D_{n}^{(x)}) \le \sqrt{2} \times 4|u-x| < \sqrt{2} \times 4 r_{k+K}$ which is smaller than $\frac12 \, r_k$ (since $K\ge 4$).

We now prove the second inclusion in \eqref{k+K+5} by discussing on two possible situations. If $r_{k+K+1} \ge \frac{4}{2^n}$, then trivially $|u-x| \ge r_{k+K+1} > r_{k+K+5} +  \sqrt{2} \times 2^{-n}$, which yields $\partial D_{n}^{(u)} \subset {\cal B}_{k+K+5}^c$. If, on the other hand, $r_{k+K+1} < \frac{4}{2^n}$, then $r_{k+K+5} < \frac14 \, 2^{-n}$, which yields ${\cal B}_{k+K+5} \subset   (\partial D_{n}^{(u)})^c$ (because $d(x,\partial D_{n}^{(x)})\ge \frac14 \, 2^{-n}$), which, in turn, implies $\partial D_{n}^{(u)} \subset {\cal B}_{k+K+5}^c$. As such, \eqref{k+K+5} is proved.

Recall the definition of $\varUpsilon_{x,D_1,S,D_2}(u)$ in \eqref{defSigma}; we take $D_{1}={\cal B}_k$, $S={\cal B}_{k+K+5}$ and $D_{2}=D_n^{(u)}$. Then $\varUpsilon_{x,{\cal B}_{k},{\cal B}_{k+K+5},D_n^{(u)}}(u)$ is (possibly) contributed by excursions from ${\cal C}(x,r_{k})$ to $x$, from $x$ to ${\cal C}(x,r_{k})$,   by the loops in ${\cal B}_{k}$ which hit  ${\cal C}(x, r_{k+K+5})$, and by excursions from ${\cal C}(x,r_{k})$ to itself without hitting $x$. [The latter excursions make the most significant contribution to $\varUpsilon_{x,{\cal B}_{k},{\cal B}_{k+K+5},D_n^{(u)}}(u)$.] We claim that for some positive constant $c_{10}$ only depending on $K$ and all $\eee \in {\mathcal E}_{x,{\cal B}_k,{\cal B}_{k+K+5}}$,
\begin{equation}
    \label{c_{10}}
    \xi_{x,{\cal B}_{k},{\cal B}_{k+K+5},D_{n}^{(u)}}(u,\eee_{g},\eee_{d})
    \le
    c_{10}. 
\end{equation}

To prove \eqref{c_{10}}, we start by noting that $H_{ D_{n}^{(u)}}(u, x)=0$ (because $x\notin \overline{D_{n}^{(u)}}$) and $H_{ D_{n}^{(u)}}(u, y)=0$ for any $ y \in \partial {\cal B}_{k} = {\cal C}(x, r_k)$. Hence by \eqref{xi},
\begin{align*}
	\xi_{x,{\cal B}_{k},{\cal B}_{k+K+5},D_{n}^{(u)}}(
	&u,\eee_{g},\eee_{d}) = 
\\ 
 	&\begin{cases}
	(2\pi)^2 \big( H_{{\cal B}_{k}}(u, x)  \big)^2/C_{{\cal B}_{k}, {\cal B}_{k+K+5}}(x), \qquad & \mbox{if }  \eee_{g}=\eee_{d}=x,  \\
	\xi_{{\cal B}_{k}}(u, \eee_{g},\eee_{d}), \qquad &\mbox{otherwise}.
	\end{cases}
\end{align*}

\noindent If $  \eee_{g}, \eee_{d} \in {\cal C}(x, r_k)$, we use  \eqref{xibound-2} to see that $\xi_{x,{\cal B}_{k},{\cal B}_{k+K+5},D_{n}^{(u)}}(u,\eee_{g},\eee_{d}) \le 2$; otherwise, we use the following 
   explicit computations:   $H_{{\cal B}_{k}}(u, x)= \frac1{2\pi} \log \frac{r_k}{|u-x|}$,  $C_{{\cal B}_{k}, {\cal B}_{k+K+5}}(x)=\log \frac{r_k}{r_{k+K+5}}= (K+5)\log 2$, $H_{{\cal B}_{k}}(u, y) = \frac{1}{2\pi r_k} \frac{r_k^2- |u-x|^2}{|u-y|^2}$ and $H_{{\cal B}_{k}}(x, y) = \frac{1}{2\pi r_k}$ for any $y \in {\cal C}(x, r_k)$.   Since $r_{k+K+1}\le |u-x|<r_{k+K}$,   we easily get \eqref{c_{10}}.

For $\eee$ inside ${\cal C}(x,r_{k})$ which does not hit $x$, we have by our choice of $K$ and Lemma \ref{l:poissonkernel}, $\xi_{{\cal B}_{k}}(u,\eee_{g},\eee_{d}) \le (1+\varepsilon)^2  \xi_{{\cal B}_{k}}(x,\eee_{g},\eee_{d}) $. 

With the notation of \eqref{Nxrbeta}, the number of excursions from ${\cal C}(x,r_{k})$ to $x$, from $x$ to ${\cal C}(x,r_{k})$, and the loops in ${\cal B}_{k}$ which hit  ${\cal C}(x, r_{k+K+5})$, is less than $2 N(x,r_{k+K+5})+2$. [The presence of $+2$ is due to the path from $z$ to $x$, and to the path from $x$ to $z'$.] In the notation of \eqref{Sigmaxr}, we have
\begin{equation}
    \label{cxrk}
    \varUpsilon_{x,D,{\cal B}(x,r),D_{n}^{(u)}} (u)
    \le
    (1+\varepsilon)^2\, \varUpsilon(x,r_{k}) + 2 \, c_{10}\, (N(x,r_{k+K+5})+1)
    =:
    c_{\eqref{cxrk}} \, .
\end{equation}

\noindent In view of \eqref{k+K+5}, we are entitled to apply Proposition \ref{p:martingale2} to $D_{1}={\cal B}_{k}$, $S={\cal B}_{k+K+5}$ and $D_{2}=D_n^{(u)}$, to see that,   
\begin{eqnarray*}
 &&\e_{\q^{z,z',\alpha}_{x,D}} [ M_{D_{n}^{(u)}}(u,\alpha)\, | \, N(x,r_{k+K+5}), \, \varUpsilon(x,r_{k}) ]
    \\
 &=& \e_{\q^{z,z',\alpha}_{x,D}} [ \widetilde M_{x,D_{n}^{(u)}}(u,\alpha) \, | \, N(x,r_{k+K+5}),\, \varUpsilon(x,r_{k})]
    \\
 &\le& 
    c_{\eqref{cxrk}} \, \ee^{-\alpha  \, C_{{\cal B}_{k}}(u)}\, \ee^{2\sqrt{\alpha \, c_{\eqref{cxrk}}}}.
\end{eqnarray*}

\noindent Recall the definition of $C_{S}(u)$ in Lemma \ref{l:ito-muxx}. We see that $C_{{\cal B}_{k}}(u)\ge \log(1/(2 r_{k}))$. By the choice of the event ${E}^{(x)}(r_{0},\gamma)$, we deduce the existence of a constant $c_{11}$ depending on $(K, \, \varepsilon, \, \alpha, \, \gamma)$ such that uniformly in $u$ satisfying $|u-x|<r_{K}$ and $u\notin D_{n}^{(x)}$,
\begin{equation}
    \label{Mrk}
\e_{\q^{z,z',\alpha}_{x,D}} \Big[ M_{D_{n}^{(u)}}(u,\alpha),\, {E}^{(x)}(r_{0},\gamma)\Big]
    \le 
    \frac{c_{11}}{|u-x|^{ 2(1+\varepsilon)\sqrt{\alpha \gamma}-\alpha}} \Big(\log \frac{1}{|u-x|}\Big)^{ \! 2}.
\end{equation}

\paragraph{Third (and last) case: $u\in D_{n}^{(x)}$}  Here we assume that $n$ is large enough so that $D_{n}^{(x)}\subset {\cal B}(x,r_{0})$. Recall that  $d(x,\partial D_{n}^{(x)})\ge {2^{-n} \over 4}$.  Let ${\mathfrak a}_n$ be the smallest integer $j\ge K$ such that $r_j <  {2^{-n} \over 4}$. In particular,  ${\cal B}_{{\mathfrak a}_n} \subset D_{n}^{(x)}$, so the number of excursions in ${\cal E}_{D_{n}^{(x)}}$ which hit $x$ is smaller than $N:=N(x,r_{{\mathfrak a}_n})+1$.  For all $0\le m\le L$, and for all $\widetilde \eee^{1}, ..., \widetilde \eee^{m} $ distinct and ordered excursions in ${\cal E}_{D_{n}^{(x)}}$ which {\it do not} hit $x$,  consider the set of  $\eee^{1}, ...,   \eee^{L} $ distinct and ordered excursions  in ${\cal E}_{D_{n}^{(x)}}$ such that  $\{\widetilde \eee^{1}, ..., \widetilde \eee^{m} \} \subset  \{\eee^{1}, ...,   \eee^{L}\}$ and any $\eee \in \{\eee^{1}, ...,   \eee^{L}\}\backslash  \{\widetilde \eee^{1}, ..., \widetilde \eee^{m} \}  $ hits $x$.   We remark that the cardinality of this set is less than $ \binom{N}{L-m}$. 

Recall from \eqref{xibound} that $ \xi_{D_{n}^{(x)}}(u,y,y')\le 2$ for any $y,y'\in \partial D_{n}^{(x)}$.   It follows    that 
\begin{equation*} 
  \sum_{\underline \eee_{D_{n}^{(x)}}^L  } \prod_{j=1}^{L}\xi_{D_{n}^{(x)}}(u,\eee_{g}^j, \eee_{d}^j)  
  \le 
 \sum_{ m=0}^L     2^{L-m} \, \binom{N}{L-m} \,  \Theta^{(m)}_{n}(u) , 
\end{equation*}
 
 \noindent where $  \Theta^{(0)}_{n}(u)  :=1$ and for any $m\ge1$  \footnote{$\binom{N}{L-m}=0$ if $N < L-m$.}
 $$ \Theta^{(m)}_{n}(u) := \sum_{\underline \eee^m_{D_{n}^{(x)}}  } \prod_{j=1}^{m}\xi_{D_{n}^{(x)}}(u,\eee_g^j, \eee_{d}^j)1_{\{ x \notin \eee^j  \}}.$$

\noindent We get that for any $u \in D_n^{(x)}$,  
\begin{eqnarray}
&&  
	M_{D_n^{(x)}}(u, \alpha)   \nonumber \\
&\le&
	\ee^{-\alpha C_{D_{n}^{(x)}}(u)}\sum_{L=1}^\infty \frac{\alpha^{L-1}}{(L-1)!} \,  \sum_{ m =0}^L 2^{L-m} \binom{N}{L-m}   \,  \Theta^{(m)}_{n}(u)   		\nonumber
\\
&\le&
\ee^{-\alpha C_{D_{n}^{(x)}}(u)} \sum_{ m = 0}^\infty  \Theta^{(m)}_{n}(u)  \, \sum_{L=\max(m, 1)}^\infty 2^{L-m}   \frac{\alpha^{L-1}}{(L-1)!} \,  \binom{N}{L-m}  .  \label{eq:YDN1}
\end{eqnarray}

To estimate the   sum $\sum_{L=\max(m, 1)}^\infty (\cdots)$  in  \eqref{eq:YDN1}, we shall use the following elementary inequalities:  for any $b$ and $j$ two integers and for all $s\ge 0$,  \begin{equation}\label{inegalitebinomiale}
\sum_{\ell = 0}^\infty {s^\ell \over (b+\ell)!}  \binom{j}{\ell}
\le
{1\over b!} \sum_{\ell = 0}^\infty {s^\ell \over \ell! }   \binom{j}{\ell}
\le
{1\over b!} \sum_{\ell = 0}^\infty {(s j )^\ell \over (\ell! )^2}   
\le
{1\over b!} \ee^{2\sqrt{s\, j}}.
\end{equation}

When $m=0$, the   sum $\sum_{L=\max(m, 1)}^\infty (\cdots)$  in \eqref{eq:YDN1} is equal to 
\begin{eqnarray*}
	\sum_{L= 1}^\infty 2^{L}   \frac{\alpha^{L-1}}{(L-1)!} \,  \binom{N}{L} 
&=& 
	2 \sum_{\ell=0}^\infty \frac{(2 \alpha)^\ell}{\ell!} \,  \binom{N}{\ell+1} 
\\
&\le& 
	2 \, N \, \sum_{\ell=0}^\infty \frac{(2\alpha)^\ell}{\ell!} \,  \binom{N-1}{\ell } 
\\
&\le& 
	2 N \ee^{2 \sqrt{2\alpha (N-1)}},
	\end{eqnarray*}

\noindent by  using  \eqref{inegalitebinomiale} for the last inequality.  When $m\ge1$,   the   sum $\sum_{L=\max(m, 1)}^\infty (\cdots)$  in \eqref{eq:YDN1} is equal to, after a change of variables $L=m+\ell$,  $$ \alpha^{m-1} \, \sum_{\ell=0}^\infty \frac{(2 \alpha)^\ell}{(m-1+\ell)!} \,  \binom{N}{\ell } 
\le
\frac{\alpha^{m-1}}{(m-1)!} \ee^{2 \sqrt{2\alpha N}},$$

\noindent by using again \eqref{inegalitebinomiale} for the last inequality.  It follows from \eqref{eq:YDN1} that  \begin{eqnarray}
  M_{D_n^{(x)}}(u, \alpha)  
& \le&
\ee^{-\alpha C_{D_{n}^{(x)}}(u)}\Big[2\, N + \sum_{m=1}^\infty \, \frac{\alpha^{m-1}}{(m-1)!}\, \Theta^{(m)}_{n}(u)\Big] \, \ee^{2 \sqrt{ 2\alpha \, N}}  \nonumber \\
&=& \Big[2\, N \ee^{-\alpha C_{D_{n}^{(x)}}(u)}+ \widetilde M_{x,D_{n}^{(x)}}(u,\alpha) \Big] \, \ee^{2 \sqrt{2 \alpha \, N}},  \label{eq:YDN}
\end{eqnarray}

\noindent with the notation of Proposition \ref{p:martingale2}. % (recalling that $N:=N(x,r_{{\mathfrak a}_n})+1$). 
We want to use Proposition \ref{p:martingale2} with $D_{1} = {\cal B}_{{\mathfrak a}_n-K}$, $D_{2}=D_{n}^{(x)}$ and $S={\cal B}_{{\mathfrak a}_n}$ (Note that   $D_{n}^{(x)}\subset {\cal B}_{{\mathfrak a}_n-K}$ as $K\ge 4$).   Remark that   uniformly in $u \in D_n^{(x)}  $,  \begin{equation}\label{01082018} 0 \le H_{ {\cal B}_{{\mathfrak a}_n-K}}(u, x) - H_{D_{n}^{(x)}}(u, x) \le c_{12}(K). 
\end{equation}

\noindent   In fact, if $|u-x| \ge \frac{2^{-n}}{8}$, then $H_{ {\cal B}_{{\mathfrak a}_n-K}}(u, x)= \frac1{2\pi} \log \frac{r_{{\mathfrak a}_n-K}}{|u-x|} \le \frac{\log 2}{2\pi}  (K+1)$, by using the fact that $r_{{\mathfrak a}_n} < \frac{2^{-n}}{4}$. Now if $|u-x| < \frac{2^{-n}}{8}$, then $d(u, \partial D_n^{(x)}) >  \frac{2^{-n}}{8}$ (recall that $d(x, \partial D_n^{(x)}) > \frac{2^{-n}}{4}$).  By Lawler \cite{Lawler}, Proposition 2.36, $H_{ {\cal B}_{{\mathfrak a}_n-K}}(u, x) - H_{D_{n}^{(x)}}(u, x) =\frac1{2\pi} \e^x \log |B_{T_{\partial {\cal B}_{{\mathfrak a}_n-K}}}- u |  - \frac1{2\pi} \e^x \log |B_{T_{\partial D_n^{(x)}}}- u | \le \frac1{2\pi} \log \frac{2  r_{{\mathfrak a}_n-K}}{   (2^{-n}/8)} \le \frac{\log 2}{2\pi}  (K+2)=: c_{12}(K).$  This proves \eqref{01082018}. 

 Using  \eqref{01082018},  the similar computations leading to   \eqref{c_{10}}     show that  there  is some  positive constant $c_{13}$ only depending on $K$ such that  for all $u \in D_n^{(x)}$, for  all $\eee \in {\mathcal E}_{x,{\cal B}_{{\mathfrak a}_n-K},{\cal B}_{{\mathfrak a}_n}}$,  $$\xi_{x,{\cal B}_{{\mathfrak a}_n-K},{\cal B}_{{\mathfrak a}_n}, D_{n}^{(x)}}(u,\eee_{g},\eee_{d})
\le
c_{13}.$$

Recall the definition of $\varUpsilon_{x,D_1,S,D_2}(u)$ in \eqref{defSigma}; we take  $D_{1} = {\cal B}_{{\mathfrak a}_n-K}$, $D_{2}=D_{n}^{(x)}$ and $S={\cal B}_{{\mathfrak a}_n}$.    Remark that for all $u \in D_n^{(x)}$, $|u-x| < \sqrt{2} 2^{-n} < \sqrt{2}  r_{{\mathfrak a}_n -1} < \eta\,  r_{{\mathfrak a}_n- K}$ by the choice of $K$.   By  Lemma \ref{l:poissonkernel}, for all $y \in {\cal C}(x,  r_{{\mathfrak a}_n- K})$,  $H_{{\cal B}_{{\mathfrak a}_n- K}}(u, y) \le (1+\varepsilon) H_{{\cal B}_{{\mathfrak a}_n- K}}(x, y)$.  Then 
similarly to \eqref{cxrk}, we have   that  
\begin{equation}\label{cxrn0}
 \varUpsilon_{x,{\cal B}_{{\mathfrak a}_n-K},{\cal B}_{{\mathfrak a}_n}, D_{n}^{(x)}}(u) 
 \le
  (1+\varepsilon)^2\varUpsilon(x,r_{{\mathfrak a}_n-K})+ 2\, c_{13} \, (N(x,r_{{\mathfrak a}_n})+1)
  =: c_{\eqref{cxrn0}}.
\end{equation}

  Applying Proposition \ref{p:martingale2} to  $D_{1} = {\cal B}_{{\mathfrak a}_n-K}$, $D_{2}=D_{n}^{(x)}$ and $S={\cal B}_{{\mathfrak a}_n}$   gives that 
\begin{eqnarray*}
&&
	\e_{\q^{z,z',\alpha}_{x,D}} \Big[ \widetilde M_{D_{n}^{(u)}}(u,\alpha)  \, \big| \,  N(x,r_{{\mathfrak a}_n}) , (1+\varepsilon)^2\varUpsilon(x,r_{{\mathfrak a}_n-K}) \Big]
\\
&\le&
	c_{\eqref{cxrn0}}\,  \ee^{-\alpha \, C_{{\cal B}_{{\mathfrak a}_n-K}}(u)}\, \ee^{2\sqrt{\alpha \, c_{\eqref{cxrn0}}}} .
\end{eqnarray*}

On the event ${E}^{(x)}(r_{0},\gamma)$, $N(x,r_{{\mathfrak a}_n}) \le \gamma \, \log \frac1{r_{{\mathfrak a}_n}}$ and $\varUpsilon(x,r_{{\mathfrak a}_n-K}) \le \gamma \, ( \log \frac1{r_{{\mathfrak a}_n-K}})^2$, hence  $$
c_{\eqref{cxrn0}} 
\le
(1+\varepsilon)^2   \gamma \, ( \log  \frac1{r_{{\mathfrak a}_n-K}})^2+ 2\, c_{13} \, (\gamma \, \log \frac1{r_{{\mathfrak a}_n}}+1) 
\le 
(1+\varepsilon)^2   \gamma \, (\log \frac{c_{14}}{r_{{\mathfrak a}_n}})^2,
$$

\noindent where $ c_{14}= c_{14}(K, \gamma, r_0)>0$ denotes some constant in the last inequality. Therefore  
\begin{eqnarray*}
&& 
	\e_{\q^{z,z',\alpha}_{x,D}} \Big[ \widetilde M_{D_{n}^{(u)}}(u,\alpha) , \, {E}^{(x)}(r_{0},\gamma)\Big] \\
&\le&
	(1+\varepsilon)^2   \gamma \, (\log \frac{c_{14}}{r_{{\mathfrak a}_n}})^2 \,   \ee^{-\alpha \, C_{{\cal B}_{{\mathfrak a}_n-K}}(u)}\, \big(\frac{c_{14}}	{r_{{\mathfrak a}_n}}\big)^{2 (1+\varepsilon) \sqrt{\alpha \gamma}} .
\end{eqnarray*}

By definition, $C_{{\cal B}_{{\mathfrak a}_n-K}}(u)= \int_{{\cal C}(x, r_{{\mathfrak a}_n-K})} \log \frac1{|u-y|} H_{{\cal B}_{{\mathfrak a}_n-K}}(u, y ) \d y  \ge \log \frac1{2  r_{{\mathfrak a}_n-K}}  $ and in the same way, $C_{D_n^{(x)}}(u) \ge \log \frac{1}{\sqrt{2}\,  2^{-n}}$.  Using  \eqref{eq:YDN} yields  that uniformly in $u\in D_{n}^{(x)}$, 
\begin{eqnarray}\label{Mrn0a}
&&
	\e_{\q^{z,z',\alpha}_{x,D}} \Big[ M_{D_{n}^{(u)}}(u,\alpha),\, {E}^{(x)}(r_{0},\gamma)\Big] \\
&\le& 
	c_{15} \, (\log \frac1{r_{{\mathfrak a}_n}})^2 \,  \ee^{ 2 \sqrt{2\alpha \gamma \log  {1\over r_{{\mathfrak a}_n}}}} \, \Big({1\over r_{{\mathfrak a}_n}} 	\Big)^{2(1+\varepsilon)\sqrt{\alpha \gamma}-\alpha }. \nonumber
\end{eqnarray}

Recall that  $ r_{{\mathfrak a}_n} < \frac14 2^{-n} \le  r_{{\mathfrak a}_n-1}$.  We deduce from \eqref{Mrn0a}  that  \begin{equation}\label{Mrn0} \e_{\q^{z,z',\alpha}_{x,D}}\left[ \int_{u \in D_n^{(x)}} \, \frac{{\mathcal M}_{\D_{n}}^\alpha(\d u)}{|u-x|^\varrho},\, {E}^{(x)}(r_{0},\gamma)  \right] 
\to 
0, \qquad n \to \infty.
\end{equation}

\noindent Then equation \eqref{c_5(N)} comes from \eqref{gamma}, \eqref{Mxr}, \eqref{Mrk} and \eqref{Mrn0}.  This completes the proof of Proposition \ref{p:Mfini}. \hfill$\Box$

 \section{Image of of ${\mathcal M}_{\infty}^{\alpha}$ under a conformal mapping}\label{s:conformal}

 The following result is analogous to \cite{BBK}, Theorem 5.2, with the extensions to all $\alpha \in (0, 2)$. Recall that $C_D(x)= - \int_{\partial D} \log (| x- y|) \, H_D(x, y) \d y$ and $\xi_D(x, z, z') = \frac{2\pi H_D(x, z) H_D(x, z')}{H_D(z, z')}$.     

\begin{proposition}\label{p:spine} Let $0\le \alpha< 2$. Recall that $D$ is a simply connected nice domain and $z,z'$ nice points of $\overline D$. For any nonnegative measurable function $f$, we have \begin{equation}\label{spine}
\e^{z,z'}_{D} \int_{D}  f(x,  B) {\mathcal M}_\infty^\alpha(\d x)=\int_{D}  \e_{\q^{z,z',\alpha}_{x, D}} \Big( f(x, B)\Big)  \, \ee^{- \alpha C_D(x)}\, \xi_{D}(x,z,z') \d x.
\end{equation}
Consequently, if  for a family of events $(E_x)_{x \in D}$, $\q^{z,z',\alpha}_{x, D}(E_x)=1$ for Lebesgue-a.e.  $x\in D$, then with probability one, the event $E_x$ holds for ${\mathcal M}_\infty^\alpha$-almost all $x$. 
\end{proposition}

{\noindent\it Proof of  Proposition \ref{p:spine}.} By the standard  monotone class argument,   it is enough to prove the equality    \eqref{spine} for $f(x, B)= 1_A(x) \, 1_F$, where $A\subset D$ a Borel set and  $F \in \F_{{\mathscr D}_n}$ for an arbitrary  $n\ge1$.   Then the left-hand-side of \eqref{spine}  is equal to 
\begin{eqnarray*}
	\e^{z,z'}_{D}  \big( {\mathcal M}_\infty^\alpha(A) 1_F\big)
&=&
	 \int_A \e^{z,z'}_{D}  \big( M_{D_n^{(x)}}(x, \alpha) 1_F\big) \d x \\
 &=&
  \int_A \q^{z,z',\alpha}_{x, D}(F) \, \ee^{- \alpha C_D(x)}\, \xi_{D}(x,z,z') \d x,
 \end{eqnarray*}

\noindent where the first equality follows from Theorem \ref{t:conv-p} and the second from Corollary \ref{c:radon}. \hfill$\Box$

%\medskip
%\begin{corollary}\label{c:portmanteau} Let $0< \alpha< 2$. Recalling $D:= (0, 1)^2$ and $z,z'$ nice points of $\partial D$.   Then $\p^{z,z'}_{D}$-a.s. 
%\end{corollary} 

\medskip

 \medskip 
 % Obviously  for any open square $\widetilde D$,  we may construct the measure ${\mathcal M}_{\infty}^{\alpha}$ in the domain $\widetilde D$ exactly by  the same way as for the square $D=(0, 1)^2$. 

%We still take $D=(0,1)^2$. 

As pointed out in \cite{BBK}, equation \eqref{p:spine} characterizes the measure ${\mathcal M}_\infty^\alpha$. To be more precise, under the assumptions of Proposition \ref{p:spine}, suppose that there is a random finite measure ${\mathfrak m}$ on the Borel sets of ${\mathbb R}^2$, measurable with respect to the Brownian $B:=(B_t)_{0\le t\le T_{z'}}$, and  which also verifies equation \eqref{spine}, replacing ${\mathcal M}_\infty^\alpha$ by ${\mathfrak m}$ there. Then ${\mathfrak m}=M_\infty^\alpha$ \, $\p^{z,z'}_D$-a.s. To see it, define ${\mathfrak m}':= {\mathfrak m} -M_\infty^\alpha$. Equation \eqref{spine} applied for ${\mathfrak m}$ and $M_\infty^\alpha$ yields that, for any bounded measurable function $f$,
$$
\e^{z,z'}_{D} \int_{D}  f(x,  B) {\mathfrak m}'(\d x)= 0.
$$
Let $A$ a Borel set of ${\mathbb R}^2$. For $c>0$, taking $f(x,B):={\mathfrak m}'(A)1_{x\in A}1_{|{\mathfrak m}'(A)| < c}$, we get that $\e^{z,z'}_{D}[{\mathfrak m}'(A)^21_{|{\mathfrak m}'(A)| < c}]=0$, hence by monotone convergence, $\e^{z,z'}_{D}[{\mathfrak m}'(A)^2]=0$. We deduce the claim.

\medskip

Let $D$, $D'$ be simply connected nice domains, and $\Psi:D\to D'$ a conformal transformation.
% which maps $D$ onto $D'$. 
Let $z$, $z'$ be distinct nice points of $\overline D$ such that $\Psi(z)$ and $\Psi(z')$ are nice points of $\overline{D'}$ ($\Psi$ can be extended to a  conformal transformation   of  $D \cup {\cal B}(z, \varepsilon) \cup {\cal B}(z', \varepsilon)$ for some $\varepsilon>0$, see Lawler \cite{Lawler}, p.~48).

\begin{proposition}\label{p:invariance}
The image measure of ${\mathcal M}_{\infty}^{\alpha}$ by $\Psi$ under $\p^{z,z'}_{D}$ has the same law as the measure $|\Psi'(\Psi^{-1}(x))|^{-2-\alpha} {\mathcal M}_{\infty}^{\alpha}(\d x)$ under $\p^{\Psi(z),\Psi(z')}_{D'}$.
\end{proposition}
 
\noindent {\it Proof of Proposition \ref{p:invariance}}.  It suffices to show that for any nonnegative measurable function $f$,
\begin{eqnarray}\label{eq:inv0}
\nonumber && \e^{z,z'}_{D} \int_{D}  f(\Psi(x),  \Psi(B)) {\mathcal M}_\infty^\alpha(\d x)  
\\
 &=& \e^{\Psi(z),\Psi(z')}_{D'} \int_{D'}  f(x,  B) |\Psi'(\Psi^{-1}(x))|^{-2-\alpha}  {\mathcal M}_\infty^\alpha(\d x).
\end{eqnarray}
From Proposition \ref{p:spine}, the left-hand side is
\begin{equation} \label{eq:inv1}
\int_{D}  \e_{\q^{z,z',\alpha}_{x, D}} \Big( f(\Psi(x), \Psi(B))\Big)  \, \ee^{- \alpha C_D(x)}\, \xi_{D}(x,z,z') \d x.
\end{equation}

\noindent The conformal invariance of $\q^{z,z',\alpha}_{x, D}$ (which results from the conformal invariance of the Brownian measures $\p^{z,x}_D$, $\p^{x,z'}_D$ and $\nu_D(x,x)$) implies that 
$$
\e_{\q^{z,z',\alpha}_{x, D}} \Big( f(\Psi(x), \Psi(B))\Big) = \e_{\q^{\Psi(z),\Psi(z'),\alpha}_{\Psi(x), D'}} \Big( f(\Psi(x), B)\Big).
$$

\noindent Moreover, $\xi_{D}(x,z,z')=\xi_{D'}(\Psi(x),\Psi(z),\Psi(z'))$ and $C_{D}(x)=C_{D'}(\Psi(x))+ \log |\Psi' (x)|$. By the change of variables $y=\Psi(x)$, equation \eqref{eq:inv1} becomes
$$
\int_{D'}  \e_{\q^{\Psi(z),\Psi(z'),\alpha}_{y, D'}} \Big( f(y, B)\Big)  \, \ee^{- \alpha C_{D'}(y)}\, \xi_{D'}(y,\Psi(z),\Psi(z'))  |\Psi'(\Psi^{-1}(y))|^{-2-\alpha}\d y
$$

\noindent which proves \eqref{eq:inv0} by another use of Proposition \ref{p:spine}. $\Box$

%By the conformal invariance of the Brownian motion, and the fact that $C_{\Psi(D)}(\Psi(y))=C_{D}(y)-\log |\Psi' (y)|$, we see that the image measure of ${\mathcal M}_{\D_{n}}^{\alpha}$ by $\Psi$ under $\p^{z,z'}_{D}$ has the same law as the measure $|\Psi'(y)|^{-2-\alpha}M_{\Psi(D_{n}^{(y)})}(x,\alpha)\d x$ under $\p^{\Psi(z),\Psi(z')}_{D'}$, where $y:=\Psi^{-1}(x)$. Therefore, we only have to show that the measure $|\Psi'(y)|^{-2-\alpha}M_{\Psi(D_{n}^{(y)})}(x,\alpha) \d x$ converges as $n\to+\infty$ to $|\Psi'(y)|^{-2-\alpha}{\mathcal M}_{\infty}^{\alpha}(\d x)$ under $\p^{\Psi(z),\Psi(z')}_{D'}$.  But this follows from an application of  \eqref{eq:partition} to  ${{\mathbb K}}_{n}$   the image of the partition $\D_{n}$ by $\Psi$. \hfill$\Box$ 

%Let $A$ a Borel set of $D$. Define the partition ${{\mathbb K}}_{n}$ as the image of the partition $\D_{n}$ by $\Psi$.  Applying  \eqref{eq:partition}, we get  that %$$
%   \e^{\Psi(z),\Psi(z')}_{D}\left[ \int_{A} |\Psi'(\Psi^{-1}(x))|^{-2-\alpha} {\mathcal M}_{\infty}^{\alpha}(\d x)\, \bigg| \, \F_{{\mathbb K}_{n}} \right] = \int_{A} |\Psi'(\Psi^{-1}(x))|^{-2} M_{\Psi(D_{n}^{(\Psi^{-1}(x))})}(x,\alpha) \d x.
 %$$

 \begin{proposition}\label{p:loizeroun}    Let   $z, z'$ be distinct nice points of $\overline D$.    We have $\p^{z,z'}_{D}({\mathcal M}_{\infty}^{\alpha}(D)>0)=1$.
 \end{proposition}
 
 \noindent {\it Proof of Proposition \ref{p:loizeroun}}.   To stress the dependence of $D$ in ${\mathcal M}_{\infty}^{\alpha}$,  we write ${\mathcal M}_{\infty, D}^{\alpha}$ in this proof.  We want to prove that  $\p^{z,z'}_{D}({\mathcal M}_{\infty, D}^{\alpha} =0)=0$. By Proposition \ref{p:invariance},  $c:=\p^{z,z'}_{D}({\mathcal M}_{\infty, D}^{\alpha} =0)$ is  independent of  $z, z'$ and $D$. Consider $D=(0,1)^2$, and $\widetilde D \in \D_2$ a square of mesh $\frac14$. The event $\{{\mathcal M}_{\infty, D}^{\alpha} =0\}$ implies that $\{{\mathcal M}_{\infty, D}^{\alpha}(\widetilde D) =0\}$. Conditionally on  $\F^+_{\widetilde D}$ and on  $y, y'$ the starting and ending points of an excursion inside $\widetilde D$ (if   it exists),  $$ \p^{z,z'}_{D} \Big({\mathcal M}_{\infty, D}^{\alpha} (\widetilde D) =0 \, \big| \, \F^+_{\widetilde D} \Big) \le 
 \p^{y, y'}_{\widetilde D}({\mathcal M}_{\infty, \widetilde D}^{\alpha} =0) 
 =c.$$

\noindent Let $\mathscr K:= \#\{\widetilde D \in \D_2 : \mbox{there is an excursion inside } \widetilde D\}$.   It follows that $$ c= \p^{z,z'}_{D} \Big({\mathcal M}_{\infty, D}^{\alpha} (\widetilde D) =0   \Big)  
 \le
 \e^{z,z'}_{D} \Big( c^{\mathscr K}\Big). $$
 
 \noindent  Note that $c<1$,   $\mathscr K  \in    \{1, 2, 3, 4\}$,  $\p^{z,z'}_{D}$-a.s.   and $\mathscr K \ge 2$ holds with positive probability.   Then the only possibility is that $c=0$.  \hfill $\Box$

\section{Thick points}\label{s:thickpoints}

Recalling $D$ a simply connected nice domain and $z\neq z'$ nice points of $\overline D$.  By using  Proposition \ref{p:spine}, we immediately deduce from    Lemma \ref{l:spine} (i)  the following result.   Let $N(x,r)$ be as in \eqref{Nxrbeta}  the number of loops from $x$ that hit  ${\cal C}(x, r)$. 

\begin{corollary}\label{c:nxr} Let $0< \alpha< 2$.      With  $\p^{z,z'}_{D}$-probability one, the measure $ {\mathcal M}_\infty^\alpha$ is supported on the set of points $x$ such that $$ 
\lim_{r \to 0^+} \frac{N(x,r)}{\log 1/r}=\alpha.
$$
\end{corollary}

Corollary \ref{c:nxr} is an extension of  \cite{BBK}, Corollary 5.1, to all $\alpha\in (0,2)$. 

In the sequel, we establish the relationship between $ {\mathcal M}_\infty^\alpha$ and  the    thick points defined in   \eqref{perfectlythick}.

\begin{theorem}\label{p:dprz} Let $\alpha \in (0, 2)$.   With $\p^{z,z'}_{D}$-probability one, the measure $ {\mathcal M}_\infty^\alpha$ is supported by the set of $\alpha$-thick points. 
\end{theorem}

By assuming Theorem  \ref{p:dprz}, we are able to complete the proof of Theorem \ref{t:main}:

{\noindent\it Proof of Theorem \ref{t:main}.} By  Corollary \ref{c:nxr} and Theorem  \ref{p:dprz}, it remains to check \eqref{carryingdimension}.

The upper bound of  \eqref{carryingdimension} follows from that of Theorem B. In fact,    for any sufficiently small $r>0$, the law of $(B_{T_{{\cal C}(z, r)\cap D}+t}, 0 \le t \le T_{{\cal C}(z', r)\cap D}- T_{{\cal C}(z, r)\cap D})$ under $\p^{z, z'}_D$ is absolutely continuous with respect  to that of a planar Brownian motion.   Then we deduce   from  Theorem B     that   $$
\mathrm{dim_H}\left\{ \alpha\mbox{-thick points}\right\}=2-\alpha,  \qquad \p^{z,z'}_{D}\mbox{-a.s.,}
 $$ 

\noindent which  in view of Theorem  \ref{p:dprz}  yields  the upper bound in \eqref{carryingdimension}.

For the lower bound in \eqref{carryingdimension},   let $ {E}^{(x)}(r_{0},\gamma) $ be the event defined in \eqref{goodevent}.    Corollary \ref{c:Minfini}  says   that   
$\p^{z, z'}_D$-a.s., ${\mathcal M}_{\D_{n}}^\alpha$ converges weakly to ${\mathcal M}_\infty^\alpha$.  By Theorem \ref{t:conv-p},  the probability $Q_{D}^{z, z', \alpha}$    defined in \eqref{QAzz'} with $A=D$ there, is absolutely continuous with respect to $\p^{z, z'}_D$.  Hence  ${\mathcal M}_{\D_{n}}^\alpha$ converges weakly to ${\mathcal M}_\infty^\alpha$,    $Q_{D}^{z, z', \alpha}$-a.s. In other words,  for Lebesgue-a.e. $x\in D$, $\q^{z,z',\alpha}_{x,D}$-a.s., ${\mathcal M}_{\D_{n}}^\alpha$ converges weakly to ${\mathcal M}_\infty^\alpha$.

 Let $0< \varrho < 2-\alpha$.    Using  the same arguments leading to \eqref{cnuxa} and by replacing ${\mathcal M}_{\D_{n}}^\alpha( D)$ by $ \int_{D} \, \frac{{\mathcal M}_{\D_{n}}^\alpha(\d u)}{|u-x|^\varrho}$,   we deduce from \eqref{c_5(N)}    that  $$ \e_{\q^{z,z',\alpha}_{x,D}}\left[ \int_{D} \, \frac{{\mathcal M}_\infty^\alpha(\d u)}{|u-x|^\varrho},\, {E}^{(x)}(r_{0}, \gamma)  \right]  \le  c_5(r_{0}), $$
 
 \noindent for Lebesgue-a.e. $x\in D$. 

%\begin{remark}  Let $0< \varrho < 2- \alpha$.   Once we have proven  Theorem \ref{t:conv-p},  Corollary \ref{c:Minfini} will imply that    $\p^{z, z'}_D$-a.s., ${\mathcal M}_{\D_{n}}^\alpha$ converges weakly to ${\mathcal M}_\infty^\alpha$, which is also the case    $Q_{D}^{z, z', \alpha}$-a.s., where  $Q_{D}^{z, z', \alpha}$ is the probability defined in \eqref{QAzz'} with $A=D$ there.  It follows that  for Lebesgue-a.e. $x\in D$, $\q^{z,z',\alpha}_{x,D}$-a.s., ${\mathcal M}_{\D_{n}}^\alpha$ converges weakly to ${\mathcal M}_\infty^\alpha$ as well. We use the same argument leading to \eqref{cnuxa}, by replacing ${\mathcal M}_{\D_{n}}^\alpha( D)$ by $ \int_{D} \, \frac{{\mathcal M}_{\D_{n}}^\alpha(\d u)}{|u-x|^\varrho}$,  deduce from \eqref{c_5(N)}    that for Lebesgue-a.e. $x\in D$  with $d(x,\partial D)>\varepsilon$, we have \begin{equation} \label{cnuxarho} \e_{\q^{z,z',\alpha}_{x,D}}\left[ \int_{D} \, \frac{{\mathcal M}_\infty^\alpha(\d u)}{|u-x|^\varrho},\, \mathcal{E}^{(x)}_{N}(\gamma_{1},\varepsilon_{1})  \right]  \le  c_5(N,\varepsilon).  \end{equation}
%\end{remark}

 Let $D^{(r_{0})}$ be the set of points $x\in D$ such that $d(x,\partial D)> r_{0}$. 
It follows from  \eqref{spine} that  for any $r_{0}>0$,  
\begin{eqnarray*}
&& \e^{z,z'}_{D} \left[ \int_{D\times D } \, 1_{{E}^{(x)}(r_{0},\gamma)}  \, \frac{{\mathcal M}_\infty^\alpha(\d u) {\mathcal M}_\infty^\alpha(\d x)}{|u-x|^\varrho} \right] \\
&=& \int_{D}  \, \d x \,  \ee^{-\alpha C_D(x)} \, \xi_{D}(x,z,z')
 \e_{\q^{z,z',\alpha}_{x,D}}\left[ \int_{D} \, \frac{{\mathcal M}_\infty^\alpha(\d u)}{|u-x|^\varrho} ,\, {E}^{(x)}(r_{0},\gamma)  \right]  
 \\
 &\le&
 c_5(r_{0})\,   \int_{D}  \, \d x \,  \ee^{-\alpha C_D(x)} \, \xi_{D}(x,z,z')
 \\
 &<&   \infty.
 \end{eqnarray*}

\noindent  Hence $\p^{z,z'}_{D}$-almost surely, for all $r_{0}>0$, $$\int_{D\times D}  \,1_{{E}^{(x)}(r_{0},\gamma)} \, \frac{{\mathcal M}_\infty^\alpha(\d u) {\mathcal M}_\infty^\alpha(\d x)}{|u-x|^\varrho}  < \infty.$$
 
By \eqref{unionen} and Proposition \ref{p:spine},  $\p^{z,z'}_{D}$-almost surely, ${\mathcal M}_\infty^\alpha$-almost all $x$,  $\cup_{\ell=1}^\infty {E}^{(x)}(2^{-\ell}, \gamma)$ holds. Then by Proposition \ref{p:loizeroun},  for some $\ell \ge 1$, the measure  ${\mathfrak m}_\ell (\d x):= 1_{{E}^{(x)}(2^{-\ell}, \gamma)}  {\mathcal M}_\infty^\alpha(\d x) $ is not trivial   and $\int_{D\times D} \, \frac{ {\mathfrak m}_\ell (\d u)  {\mathfrak m}_\ell (\d x)}{|u-x|^\varrho} < \infty$.  

For any Borel set $A$ such that ${\mathcal M}_\infty^\alpha(A^c)=0$,  we have   ${\mathfrak m}_\ell (A^c)=0$ and $\int_{A \times A} \, \frac{ {\mathfrak m}_\ell (\d u)  {\mathfrak m}_\ell (\d x)}{|u-x|^\varrho} = \int_{D\times D} \, \frac{ {\mathfrak m}_\ell (\d u)  {\mathfrak m}_\ell (\d x)}{|u-x|^\varrho} < \infty.$ It follows from Frostman's lemma  that $\mathrm{dim_H}(A)\ge \varrho$. This yields that a.s., $\mathrm{Dim}({\mathcal M}_\infty^\alpha) \ge \varrho$. As $\varrho$ can be as close as possible  to $2-\alpha$, we get  the lower bound in \eqref{carryingdimension}. This completes  the proof of Theorem \ref{t:main}. \hfill$\Box$

\medskip
 Denote by $\zeta({\mathfrak e})$ the lifetime of a loop  ${\mathfrak e}$ and let \begin{equation}\label{Lambda}
\Lambda_s:=  \sum_{v \le s} \zeta({\mathfrak e}_v), \qquad 0 \le s \le \alpha.
\end{equation}

\noindent  Note that $\Lambda$ is a subordinator on $[0, \alpha]$  whose sample paths are strictly increasing.   The law of $\zeta$ under $\nu_D(x,x)$ can easily be computed:  it follows from  \eqref{def-muDxx} that    \begin{eqnarray} \nu_D(x,x) \big(\zeta \in \d t \big)  
&=&
\pi  \mu_D(x,x; t) \big(\zeta \in \d t \big) 
\nonumber \\
&=&
     \lim_{\varepsilon\to0^+} \frac1{\varepsilon^2} \p^x\big( t< T_{\partial D}, |B_t -x | < \varepsilon\big)  \d t 
\nonumber \\
&=&
\pi  \, p_D(t, x, x) \, \d t . \label{lawofzeta}
\end{eqnarray}

Let us define  $(\ell_t^x)_{t\ge 0}$ by the inverse of $\Lambda$:   \begin{equation}\label{elltx}
\ell_t^x:=
\begin{cases}
0,  \qquad &\mbox{if } t\le T_x, 
\\
\inf\{s>0: \Lambda_s >t-T_x \}, \qquad & \mbox{if } T_x \le t < T_x+\Lambda_\alpha, 
\\
\alpha, \qquad &\mbox{if } t\ge T_x+ \Lambda_\alpha  . 
\end{cases}
\end{equation}

\noindent   By construction,   $\q^{z,z',\alpha}_{x, D}$-almost surely, $t \to \ell_t^x$ is  continuous.  Moreover,  by imitating the proof in the one-dimensional Brownian motion case (for instance the proof of Proposition VI.2.5, \cite{RY}), we  get  that $\q^{z,z',\alpha}_{x, D}$-almost surely,  the support of $   \d \ell_t^x$  is exactly equal to  the level set  $\{T_x \le t\le  T_x+ \Lambda_\beta: B_t=x\}$.  %which is also identified as the image of the subordinator $\{T_x+ \Lambda_s, T_x+ \Lambda_{s-}, 0 \le s \le \alpha\}$.  
We call   $(\ell_t^x)_{t\ge 0}$ the family of local times of $B$ at $x$. We mention that after $ T_x+ \Lambda_\beta$, $B$ is a Brownian motion started at $x$ and conditioned to hit  $z'$.   Theorems   \ref{p:dprz} and   \ref{t:main2} will be a consequence of   the following result:

\begin{proposition} \label{p:thick} For any $x\in   D$,  $\q^{z,z',\alpha}_{x, D}\mbox{-a.s.}$ for all $0\le t \le T_{z'},$
\begin{equation}\label{thick1}
\lim_{r\to 0^+} \frac1{r^2 (\log r)^2} \, \int_0^t 1_{\{ B_s \in {\cal B}(x, r)\}} \d s = \ell_t^x,
\end{equation}
where $ \ell_t^x$,  defined in \eqref{elltx}, denotes the local time at $x$ up to time $t$ (under $\q^{z,z',\alpha}_{x, D}$). 
\end{proposition}

By assuming Proposition \ref{p:thick} for the moment, we give the proofs of Theorems   \ref{p:dprz} and  \ref{t:main2}:

\medskip
{\noindent\it Proof of Theorem   \ref{p:dprz}.}  By definition,   $\ell_{T_{z'}}^x=\alpha$. It follows from  \eqref{thick1} that $\q^{z,z',\alpha}_{x, D}\mbox{-a.s.}$, $\lim_{r\to 0^+} \frac1{r^2 (\log r)^2} \, \int_0^{T_{z'}} 1_{\{ B_s \in {\cal B}(x, r)\}} \d s = \alpha$, which in view of Proposition \ref{p:spine} yields Theorem   \ref{p:dprz}. \hfill$\Box$

\medskip
{\noindent\it Proof of Theorem   \ref{t:main2}.} For any   $x \in D$, we define $${\mathfrak L}_t^x:= \lim_{r\to 0^+} \frac1{r^2 (\log r)^2} \, \int_0^{t\wedge T_{z'}} 1_{\{ B_s \in {\cal B}(x, r)\}} \d s,  \; \mbox{if the limit exists for all $t\ge0$},$$

\noindent and ${\mathfrak L}_t^x:= 0, \, \forall\, t\ge0,$ otherwise.  By Proposition \ref{p:thick}, for any $x\in   D$,  $\q^{z,z',\alpha}_{x, D}\mbox{-a.s.}$, ${\mathfrak L}_t^x=\ell_{t\wedge T_{z'}}^x$  for all $t\ge0$. Then  $\q^{z,z',\alpha}_{x, D}\mbox{-a.s.}$,  $t \to {\mathfrak L}_t^x$   is a continuous additive functional and such that $t \to \d {\mathfrak L}_t^x$ is supported by $\{t\in [0, T_{z'}]: B_t=x\}$. By Proposition \ref{p:spine}, we get Theorem   \ref{t:main2}. \hfill$\Box$

\medskip

The rest of this section is devoted to the proof of Proposition \ref{p:thick}. Recall $\Lambda$   defined in \eqref{Lambda}.   Let   for any $r>0$ and $0\le \beta< \alpha$,   $$
I_\beta(r)  :=   \int_0^{T_x+ \Lambda_\beta} 1_{\{ B_s \in {\cal B}(x, r)\}} \d s 
\qquad \mbox{and } \quad 
I_\alpha(r)  :=   \int_{0}^{T_{z'}} 1_{\{ B_s \in {\cal B}(x, r)\}} \d s.
 $$

Notice that $I_0(r)= \int_0^{T_x } 1_{\{ B_s \in {\cal B}(x, r)\}} \d s $ and $I_\alpha(r)- I_{\alpha-}(r)=  \int_{T_x+ \Lambda_\alpha}^{T_{z'}} 1_{\{ B_s \in {\cal B}(x, r)\}} \d s$ correspond to the times spent  in ${\cal B}(x, r)$ by  a conditioned Brownian motion  under $\p^{z, x}_D$ and under $\p^{x, z'}_D$ respectively. In the literature, there is a law of the iterated logarithm  for their asymptotics as $r\to0^+$, see Ray \cite{Ray63} and Le Gall \cite{Legall85}. 

%In other words,    $I_1(r)$ corresponds to the time spent by the   Brownian motion conditioned to hit $x$,  $I_2^{(\beta)}(r)$  to  the concatenation of Brownian loops in $D$ at $x$ up to time $\beta$ (in the time scale of Poisson point process),  and  $I_3(r)$  to a conditioned Brownian motion started at $x$ conditioned to hit $z'$. Since the law of the Brownian motion under $\p^{x,z'}_{D}$ up to $T_{{\cal C}(z',\eta)}$, for $\eta>0$ small, is absolutely continuous with respect to a standard Brownian motion starting from $x$,  the law of iterated logarithm  in Ray \cite{Ray63}  says that   $$ \limsup_{r \to 0^+} \frac1{r^2   (\log \frac1{r}) Â  \log  \log \log \frac1{r}} I_3(r) = 1 , \qquad \q^{z,z',\alpha}_{x, D}\mbox{-a.s.}, $$

%\noindent the constant $1$ is due to Le Gall \cite{Legall85}.   Hence  $I_{3}(r)$  does not contribute in the limit \eqref{thick1}. By time-reversal (see \eqref{timereversal}), we can ignore $I_{1}(r)$ as well.     We proved that   \begin{equation}\label{I_1+I_3} \lim_{r\to 0^+}  \frac1{r^2 (\log r)^2} (I_1(r)+I_3(r))= 0, \qquad \q^{z,z',\alpha}_{x, D}\mbox{-a.s.} \end{equation}

\medskip
{\noindent \it Proof of Proposition \ref{p:thick}.} 
 Proposition \ref{p:thick} will follow  from the following statement:  $\q^{z,z',\alpha}_{x, D}\mbox{-a.s.}, $ for all  $\beta \in [0, \alpha]$, 
\begin{equation}\label{thickbeta}
\lim_{r\to 0^+} \frac{I_\beta(r)}{r^2 (\log r)^2}    =\beta.
\end{equation}

In fact,  take $\beta=0$ in \eqref{thickbeta} implies  \eqref{thick1}   for all $t\in [0,T_x]$. For $t\in [\Lambda_\alpha +T_x, T_{z'}]$, we remark that $I_{\alpha-}(r) \le \int_0^t 1_{\{ B_s \in {\cal B}(x, r)\}} \d s  \le I_\alpha(r)$, then by       \eqref{thickbeta} we see that \eqref{thick1} holds for these $t$.    Now consider  $T_x < t < T_x + \Lambda_\alpha$, $\beta:= \ell_t^x \in (0, \alpha)$. We have $\Lambda_{\beta-} + T_x \le t \le \Lambda_\beta +T_x$.  Then $ I_{\beta-}(r) \le \int_0^t 1_{\{ B_s \in {\cal B}(x, r)\}} \d s \le I_\beta(r)$.  Applying       \eqref{thickbeta} gives  \eqref{thick1} for any $T_x < t < \Lambda_\beta +T_x$.

 In order to prove \eqref{thickbeta}, by the monotonicity, we only need to show that for any fixed $\beta \in (0, \alpha]$,  $\q^{z,z',\alpha}_{x, D}\mbox{-a.s.}$,  \begin{equation}\label{thickbeta2}
\lim_{r\to 0^+} \frac{I_\beta(r)}{r^2 (\log r)^2}    =\beta.
\end{equation}
 
 To this end,  let   $r_{0}>0$ such that ${\cal B}(x,r_{0}) \subset D\backslash \{z,z'\}$.  For any $0< r< r_0$,  denote by  $L(r)$   the local time at ${\cal C}(x, r)$ of $(B_t)_{t\ge0}$  (as the    occupation time density at $r$   of the process  $(|B_t-x|)_{t\ge0}$)  till $T_x + \Lambda_\beta$ if $\beta < \alpha$ and till $T_{z'}$ if $\beta=\alpha$. Then $$ I_\beta(r) = \int_0^r L(u) \d u. $$

We will prove that  $\q_{x, D}^{z, z', \alpha}$-a.s., $$ \lim_{r \to 0^+}  \frac{L(r)}{r (\log r)^2} = 2 \beta ,$$ which gives \eqref{thickbeta2} by integrating $L(\cdot)$. 

Without loss of generality, we suppose that $x=0$.  We can decompose the trajectories in ${\cal B}(0,r_0)$ as 
 
 $\bullet$ Brownian excursions from ${\cal C}(0,r_0)$ to ${\cal C}(0,r_0)$ which do not hit $0$,
 
 $\bullet $  Brownian excursions from ${\cal C}(0,r_0)$ to $0$ or from $0$ to ${\cal C}(0,r_0)$, 
 
 $\bullet$ and Brownian loops at $0$ in ${\cal B}(0,r_0)$.

 The Brownian excursions which do not hit $0$ can be ignored to understand the  asymptotics of the local time around $0$ since none will hit  ${\cal C}(0,r)$ for $r$ small enough. 
 
 Let us consider the Brownian excursions from $0$ to ${\cal C}(0,r_0)$.  These excursions are of  finite number  almost  surely.  The norm of an excursion is a two-dimensional Bessel process stopped at the first  hitting time of $r_0$. Its local time process at   $r\in(0,r_0)$ is equal in law to   the process $(r U_{\log(r_0/r)})_{r \in (0, r_0)}$ where $U_{\cdot}$ is the square of    Bessel processes of dimension $2$  starting from $0$ (see Exercise 2.6, Chapter XI in Revuz and Yor \cite{RY}).  By the classical law of the iterated logarithm for the process $U$, almost surely as $r\to 0^+$, $U_{\log(1/r)}= O(\log(1/r) (\log \log \log 1/r))$.  Hence  the contribution to $L(r)$ by the excursions from $0$ to ${\cal C}(0,r_0)$ are $o(r  (\log r)^2)$. The same is true, by time-reversal, of excursions from ${\cal C}(0,r_0)$ to $0$.

 If we denote by $L_{(\mathfrak e_\cdot)}(r)$ the local time at ${\cal C}(0, r)$ of Brownian loops $({\mathfrak e}_s)_{s\le \alpha}$ in ${\cal B}(0, r_0)$, then we have shown that almost surely, $L(r) - L_{(\mathfrak e_\cdot)}(r)=o(r  (\log r)^2)$ as $r\to0^+$.   It remains to prove that  almost surely,  $ \lim_{r \to 0^+}  \frac{L_{(\mathfrak e_\cdot)}(r)}{r (\log r)^2} = 2 \beta $, or equivalently $$ \lim_{t \to \infty}  \frac{\ee^t \, L_{(\mathfrak e_\cdot)}(\ee^{-t})}{t^2} = 2 \beta . $$
 
  %  At this stage it is clear that we may take $r_0=1$ because almost surely there are only a finite number of loops in  ${\cal B}(0, 1)$ which touch  ${\cal C}(0, r_0)$, and the local times at ${\cal C}(0, r)$ of these loops  are  of order $ O(\log(1/r) (\log \log \log 1/r))$ as shown before. 
 
 By scaling, we may take $r_0=1$. 
  For a loop $\mathfrak e$, denote by $h({\mathfrak e}):= \max_{u\ge 0} |{\mathfrak e}(u)| $ its maximum norm (recalling that $x=0$).  %%The restriction of $\nu_D(0,0)$ on the set of loops in  ${\cal B}(0, r_0)$ is equal to $\nu_{{\cal B}(0, r_0)}(0,0)$. 
   By Lemma \ref{l:ito-muxx} (ii), we get  that  $$ \nu_{{\cal B}(0, 1)}(0,0) \Big( h({\mathfrak e}) > r\Big) = \log (1/r), \qquad 0< r\le 1 ,$$

   \noindent which implies that the point process ${\cal P}(\beta):= \sum_{s\le \beta : {\mathfrak e_s} \in {\cal B}(0,1)} \delta_{\{\log (1/ h({\mathfrak e_s}))\}}$ is a Poisson point process   on $(0, \infty)$ with intensity $\beta \d t$.  Write $0<{\mathfrak u}_1< {\mathfrak u}_2< ...< {\mathfrak u}_k <...$ for the points of ${\cal P}(\beta)$ and let  $N_t:= \sum_{k=1}^\infty 1_{\{ {\mathfrak u}_k \le t\}}, t\ge0.$ Then $(N_t)_{t\ge0} $ is a Poisson process  with parameter $\beta$.

 Conditionally on $\{h({\mathfrak e}) = r\}$, the norm of a Brownian loop $|{\mathfrak e}(\cdot)|$  can be decomposed as a two-dimensional Bessel process  from $0$ to  its first hitting time of $r$,  followed by an independent copy of the same process going backwards in time, see Pitman and Yor \cite{PY96}.  Therefore, conditionally on $\{h({\mathfrak e}) = \ee^{-s}\}$, the process of  local times at $\ee^{-t}$ of  $|{\mathfrak e}(\cdot)|$, for $t\ge s$,    has the law of the sum of two independent copies of $(\ee^{-t} U_{t-s})_{t\ge s}$ which  by the additivity of the square of Bessel processes, is equal in law to the process $(\ee^{-t} U^{(4)}_{t-s})_{t\ge s}$, where $U^{(4)}$ denotes the square of Bessel processes of dimension $4$, starting from $0$.  
 
 Let $U^{(4)}_k, k\ge1,$ be an i.i.d. copies of $U^{(4)}$, independent of $\{{\mathfrak u}_k\}_{k\ge1}$.  Then the process $\big( \ee^t \, L_{(\mathfrak e_\cdot)}(\ee^{-t}) \big)_{t\ge0} $ is equal in law to  the process $(X_t)_{t\ge0}$ where $$X_t:= \sum_{k=1}^\infty 1_{\{{\mathfrak u}_k \le t\}} \, U^{(k)}_{t- {\mathfrak u}_k}, \qquad t\ge0.$$

 \noindent It is enough to show that almost surely, $$ \lim_{t \to \infty} \frac{X_t}{t^2} = 2\beta.$$

\noindent   To this end,  write  $$X_t= 4  \sum_{k=1}^\infty  (t-  {\mathfrak u}_k) 1_{\{{\mathfrak u}_k \le t\}}  + \widehat X_t
  = 4 \int_{(0, t]}  (t-s) \d N_s  +  \widehat X_t,$$
  
  \noindent where $ \widehat X_t:= \sum_{k=1}^\infty 1_{\{{\mathfrak u}_k \le t\}} \, (U^{(k)}_{t- {\mathfrak u}_k} - 4 (t-  {\mathfrak u}_k)).$  Notice that $\int_{(0, t]}  (t-s) \d N_s= \int_{(0, t]}  N_s \d s$.  By the law of large numbers $\lim_{s\to \infty} \frac{N_s}{s} = \beta$, which yields that $\frac{4}{t^2}   \int_{(0, t]}  (t-s) \d N_s \to 2 \beta$ a.s. 
  
  To complete the proof, we only need to check that $\widehat X_t = o(t^2)$ a.s.   It is well-known that $U_s- 4 s$ is a martingale and $\mathrm{Var}(U_s)= 8 s^2$ (see Chapter XI in Revuz-Yor \cite{RY}).   Conditionally on $\{  {\mathfrak u}_k\}_{k\ge1}$,  $  \widehat X$ is a (finite) sum of independent martingale hence is a martingale. By Doob's $L^2$-inequality, we get that for any $t>0$, $$
  \e  ( \sup_{0\le s \le t}   \widehat X_s^2)
  \le
  4 \e  (\widehat X_t ^2)
  =
  32 \,  \e  \sum_{k=1}^\infty 1_{\{{\mathfrak u}_k \le t\}} \, (t- {\mathfrak u}_k)^2
  =
  \frac{32 \beta}{3} t^3.$$
  
   \noindent Using the Borel--Cantelli lemma gives that almost surely, $  \sup_{0\le s \le t} |\widehat X_s|  = O(t^{3/2 +o(1)})$ as $t\to \infty$. Hence  $\widehat X_t = o(t^2)$ a.s.,  which  completes the proof of   Proposition \ref{p:thick}. \hfill$\Box$

 %Observe that $U_{t-s}$ is an exponential random variable with mean $2(t-s)$. Finally, the image by $a\to -\log(a)$ of the point process of the maximal norms of the Brownian loops is a Poisson point process with intensity $\alpha dt $. We call $L_{t}$ the local time accumulated at distance $\ee^{-t}$ by all the Brownian loops, and $X_{t}:= L_{t}/(\ee^{-t}t)$. We deduce that, for any real $\lambda<1/2$, $\ee^{\lambda X_{t}}$ has expectation $ \exp \alpha t \big( {2\lambda\over 1-2\lambda}  \big)$. Standard large deviations principle gives that $X_{k}/k$ converges almost surely to $2 \alpha$ as integer $k\to+\infty$. Doob's inequality shows that $X_{t}/t$ converges to $2\alpha$ as real $t\to\infty$.

%At first, we show  that  in the limit \eqref{thick1}, nor $I_1(r)$ and $I_3(r)$ does contribute. In  fact,   the law of iterated logarithm  in Ray \cite{Ray63}  says that  
%$$
%\limsup_{r \to 0^+} \frac1{r^2   (\log \frac1{r}) ÃÂ  \log  \log \log \frac1{r}} I_3(r) = 1 , \qquad \q^z_{x, D}\mbox{-a.s.},
%$$

%\noindent the constant $1$ is due to Le Gall \cite{Legall85}.   

\section{Further discussions}\label{s:discussion}

\subsection{The case $\alpha=0$}

 In the case $\alpha=0$, $\q^{z,z',\alpha}_{x,D}$ is the law of a Brownian excursion from $z$ to $x$ in $D$, followed by an independent  Brownian excursion from $x$ to $z'$ in $D$.   We claim that  for any $z ,   z'$ two distinct  nice  points of $\overline D$,     \begin{equation}\label{Malpha=0} {\cal M}_{\infty}^0 (\bullet) =
\pi \,  \mu(\bullet \cap  D) ,  \qquad \p^{z, z'}_D\mbox{-a.s.},
\end{equation}

\noindent where $\mu$ denotes  the occupation measure of the Brownian motion (under $\p^{z, z'}_D$)  defined by $\mu(A) := \int_0^{T_{z'}} 1_{\{B_{t} \in A\}} \d t  $ for any Borel set $A$. To identify the two measures, we only need to show that   \begin{equation}\label{spinealpha0}
\pi \, \e^{z,z'}_{D} \int_{D}  f(x,  B)  \mu(\d x)=   \int_{D}  \e_{\q^{z,z',0}_{x, D}} \Big( f(x, B)\Big)  \,   \xi_{D}(x,z,z') \d x, \end{equation}

\noindent  for $f(x, B) = \tilde h(x) \ee^{- \int_0^\infty  h(t,B_{t}) dt }$ with   nonnegative Borel  functions $h$ and $\tilde h$. 

It is enough to prove  \eqref{spinealpha0} for $z  \in D$; the case  when $z\in \partial D$  follows  by considering $z_n:= z + \varepsilon_n {\bf n}_z$ with $\varepsilon_n \to 0^+$ as $n \to \infty$.

Consider the right-hand side of \eqref{spinealpha0}.  By the strong Markov property at time $T_{x}$,  
$$
\e_{\q^{z,z',0}_{x, D}} \Big( f(x, B)\Big)  
=
\tilde h(x) \, \e_{D}^{z,x} \Big[ \ee^{-\int_{0}^{T_{x}} h(t,B_{t}) \d t }\Big]\,  \e_{D}^{x,z'}\Big[ \ee^{-\int_{0}^{T_{z'}} h(t,B_{t}) \d t }\Big].  
$$

\noindent The law $\p_{D}^{z,x}$ is  the normalized excursion measure denoted by $\mu^{\#}_{D}(z,x)$ in Lawler (\cite{Lawler}, Chapter 5.2). It  is     $\mu_{D}(z,x)$ divided by  the  mass $2 H_D(z, x)$,  and 
$$
\mu_{D}(z,x) := \int_0^\infty  \mu_{D}(z,x; s)\d s ,
$$   %%:= \lim_{\varepsilon\to 0^+} \frac1{\pi \varepsilon^2} \p^z(\bullet, |B_s-x|< \varepsilon, T_{\partial D } > s)%%
where   $\mu_{D}(z,x; s)$ denotes  the measure   on the Brownian paths in $D$ of length $s$, from $z$ to $x$ (the factor $2$ comes from the renormalization  of $H_D$ in \eqref{HDzz'}).
 %Moreover,    
%$$ |\mu_{D}(z,x)|  =
%\begin{cases}
% 2 H_D(z, x) ,  &\qquad \mbox{if } z \in D, \\
%H_D(z, x),  & \qquad \mbox{if $z\in \partial D$},
%\end{cases}
%$$
%\noindent where the factor $2$ comes from the renormalization  of $H_D$ in \eqref{HDzz'}. 
We refer to  Lawler (\cite{Lawler}, Chapter 5.2) for the precise definition of $\mu_{D}(z,x; s)$  and the following equality:  For any nonnegative Borel  function $F$ and for any $s>0$, $$
 \int_{D}  \, \mu_{D}(z,x; s)(F(B_s) \ee^{-\int_{0}^{s} h(t,B_{t}) \d t }) \, \d x
 =
 \e^z \Big[ F(B_{s}) \, \ee^{-\int_{0}^{s} h(t,B_{t}) \d t }\, , s< T_{\partial D}\Big].$$

\noindent It follows that 
\begin{eqnarray*}
&&
	\int_{D}  \, F(x)\, \e_{D}^{z,x} \Big[ \ee^{-\int_{0}^{T_{x}} h(t,B_{t}) \d t }\Big]\, \d x \\
&=&
	\int_{D}  \,   \frac{F(x)}{2 H_D(z, x)}\,  \mu_{D}(z,x)(\ee^{-\int_{0}^{T_{x}} h(t,B_{t}) \d t }) \, \d x 
\\ 
 &=& 
	\e^z_{D}\Big[\int_{0}^{T_{\partial D}} \,  \frac{F(B_{s})}{2 H_D(z, B_s)} \, \ee^{-\int_{0}^{s} h(t,B_{t}) \d t }\,  \d s\Big], 
\end{eqnarray*}

%$$
%\int_{D}  F(x) \mu_{D}(z,x)(A) \d x = \int_0^\infty \int_{D} F(x) \mu_{D}(z,x,s)(A)\d x \d s = \e^z_{D}\Big[\int_{0}^{T_{\partial D}} F(B_{s})1_{A} \d s\Big]
%$$

\noindent which implies that,  with $g(x):= \e_{D}^{x,z'}\Big[ \ee^{-\int_{0}^{T_{z'}} h(t,B_{t}) \d t }\Big]$, 
\begin{eqnarray*}
&&
	\int_{D}  \e_{\q^{z,z',0}_{x, D}} \big( f(x, B)\big)  \, \xi_{D}(x,z,z') \d x \\
&=& 
	\pi \, 
	\e^z_{D}\Big[\int_{0}^{T_{\partial D}} \tilde h(B_{s})g(B_{s})\ee^{-\int_{0}^{s} h(t,B_{t}) \d t } {H_{D}(B_{s},z') \over H_{D}(z,z')} \d s\Big].
\end{eqnarray*}

\noindent It is, by the $h$-transform \eqref{htransform2},
\begin{eqnarray*}
&&
	\pi \, \e^{z,z'}_{D}\Big[\int_{0}^{T_{z'}} \tilde h(B_{s})g(B_{s})\ee^{-\int_{0}^{s} h(t,B_{t}) \d t } \d s\Big] \\
&=&
	 \pi \, \e^{z,z'}_{D}\Big[\int_{0}^{T_{z'}} \tilde h(B_{s})\ee^{-\int_{0}^{T_{z'}} h(t,B_{t}) \d t }\d s \Big]
 \\
 &=& 
	 \pi \, \e^{z,z'}_{D}\Big[\int_{0}^{T_{z'}} f(B_{s},B)\d s \Big],
\end{eqnarray*}

\noindent  where the first equality follows from   the Markov property at time $s$. This 
 proves  \eqref{spinealpha0}. \hfill$\Box$

%Applying the case $\alpha=0$   gives that for any $z, z' \in \overline D$ nice points, for any nonnegative Borel function $f$, 
% \begin{equation}\label{spinealpha0}
%  \e^{z,z'}_{D}   \int_0^{T_{z'}} f(B_s) \d s =\frac1\pi \, \int_{D}  \, f(x) \, \xi_{D}(x,z,z') \d x.
%\end{equation}

\subsection{The intersection local times  }

Fix $p\ge2$   an integer. %  let   $(z_1, z'_1), ..., (z_p, z'_p)$  be $p$ pairs  of distinct nice points of  $\overline D$. We consider  $\otimes_{i=1}^p \p^{z_i, z'_i}_D$ the law of  $p$ independent excursions  $B^{(i)}$, $1\le i\le p$, inside $D$, from $z_i$ to $z'_i$. For any $x \in D$ and $ \alpha\in [0, 2)$, we consider   $\otimes_{i=1}^p \q^{z_i, z'_i, \alpha}_{x, D}$ the product measure. 
We take $\alpha_1,\alpha_2,\ldots,\alpha_p\ge 0$ such that $\sum_{i=1}^p \alpha_i <2$   and $p$ pairs of distinct points $(z_1,z'_1),\ldots, (z_p,z'_p)$ in $\overline{D}$.    The law of  $p$ independent excursions  $B^{(i)}$, $1\le i\le p$, inside $D$, from $z_i$ to $z'_i$ is denoted by $\otimes_{i=1}^p \p^{z_i, z'_i}_D$. For any $x \in D$, we consider   $\otimes_{i=1}^p \q^{z_i, z'_i, \alpha_i}_{x, D}$ the product measure.

Recall the notation ${M}_{D_{n}^{(x)}}(x,\alpha_i)$ in Section \ref{s:construction}.  For any Borel set $A$, the Radon--Nykodim derivative of $\int_A \otimes_{i=1}^p \q^{z_i, z'_i, \alpha_i}_{x, D}  \prod_{i=1}^p  {M}_{D}^{(i)}(x, \,\alpha_i)\d x$ with respect to  $\otimes_{i=1}^p \p^{z_i, z'_i}_D$ on the $\sigma$-field $\F_{\D_{n}}$ is given by 
$$
\mathcal{M}_{\D_{n}}^{\alpha_1,\ldots,\alpha_p}(A) := \int_{A} \prod_{i=1}^p {M}_{D_{n}^{(x)}}^{(i)}(x, \,\alpha_i) \,\d x,
$$

\noindent where ${M}_{D_{n}^{(x)}}^{(i)}(x, \,\alpha_i)$ is associated to the trajectory $B^{(i)}$. It is a nonnegative martingale (in $n$, adapted to the filtration $(\F_{\D_{n}})_{n\ge 0}$) hence has an almost sure limit ${\mathcal M}^{\alpha_1,\ldots,\alpha_p}_\infty$.  We want to show the uniform integrability of this martingale. We need to follow the proof of Proposition \ref{p:Mfini} in Section \ref{sub:proof}.  

We let $0 \le \varrho<  2-\sum_{i=1}^p \alpha_i$. Take $\gamma_i=\gamma_i(\varrho), i=1\ldots p$ and $\varepsilon=\varepsilon(\varrho)$ such that $\gamma_i>\alpha_i$, $\varepsilon>0$ and 
$$
\varrho + \sum_{i=1}^p 2(1+\varepsilon)\sqrt{\gamma_i\alpha_i} -\alpha_i <2.
$$

We keep the notation $\eta$, $K$, $r_{0}\in (0, \, 1)$ and $r_{k}= \frac{r_0}{2^k}$. We define ${E}^{(x),i}(r_{0},\gamma_i)$ as in equation \eqref{goodevent}, associated to the trajectory $B^{(i)}$ and the parameters $\alpha_i$, $\gamma_i$. We let  ${E}^{(x)}(r_{0},\gamma) := \bigcap_{i=1}^p {E}^{(x),i}(r_{0},\gamma_i)$. We will prove the analog of equation  \eqref{c_5(N)}:
$$
    \limsup_{n\to+\infty}   1_{\{d(x,\partial D_{n}^{(x)})\ge {2^{-n} \over 4} \}} \e_{\otimes_{i=1}^p \q^{z_i,z'_i,\alpha_i}_{x,D}}\left[ \int_{D} \, \frac{  {\mathcal M}_{\D_{n}}^{\alpha_1,\ldots,\alpha_p}(\d u)}{|u-x|^\varrho},\, {E}^{(x)}(r_{0},\gamma)  \right]  \le c_5(r_{0}).
$$
 
\noindent  The case $\varrho=0$ implies uniform integrability. We have
\begin{eqnarray*}
&&  \e_{\otimes_{i=1}^p \q^{z_i,z'_i,\alpha_i}_{x,D}}\left[ \int_{D} \, \frac{ {\mathcal M}_{\D_{n}}^{\alpha_1,\ldots,\alpha_p}(\d u)}{|u-x|^\varrho},\, {E}^{(x)}(r_{0},\gamma)  \right] \\
&=&
\int_{D} \prod_{i=1}^p \e_{\q^{z_i,z'_i,\alpha_i}_{x,D}}\left[  \, \frac{   M_{D_{n}^{(u)}}^{(i)}(u,\alpha_i)}{|u-x|^\varrho},\, {E}^{(x),i}(r_{0},\gamma_i)  \right] \d u.
\end{eqnarray*}

\noindent We now use for each trajectory $B^{(i)}$: 
\begin{itemize}
\item if $|u-x|>r_K$, the equation before equation \eqref{Mxr}:
$$
 \e_{\q^{z_i,z'_i,\alpha_i}_{x,D}}\left[ M_{D_{n}^{(u)}}^{(i)}(u,\alpha_i) ,\, {E}^{(x),\alpha_i}(r_{0},\gamma_i)  \right] 
 \le
 f_{z_i, z'_i} (u) \,  \ee^{-\alpha_i\,  C_{D}(u)}\, \ee^{2\sqrt{\alpha_i  f_{z_i, z'_i} (u)}},
$$
\item if $r_{k+K+1}\le |u-x|<r_{k+K}$ for some $k\ge 0$, with $u\notin D_{n}^{(x)}$, equation \eqref{Mrk}:
$$
\e_{\q^{z_i,z'_i,\alpha_i}_{x,D}} \Big[ M_{D_{n}^{(u)}}^{(i)}(u,\alpha_i),\, {E}^{(x)}(r_{0},\gamma_i)\Big]
    \le 
    \frac{c_{11}}{|u-x|^{ 2(1+\varepsilon)\sqrt{\alpha_i \gamma_i}-\alpha_i}} \Big(\log \frac{1}{|u-x|}\Big)^{ \! 2}.
$$
\item if $u\in D_n^{(x)}$, equation  \eqref{Mrn0a}, (where ${\mathfrak a}_n$ is the smallest integer $j\ge K$ such that $r_j <  {2^{-n} \over 4}$):
\begin{eqnarray*}
&&\e_{\q^{z_i,z'_i,\alpha_i}_{x,D}} \Big[ M_{D_{n}^{(u)}}^{(i)}(u,\alpha_i),\, {E}^{(x)}(r_{0},\gamma_i)\Big]
\\
&\le& 
c_{15} \, (\log \frac1{r_{{\mathfrak a}_n}})^2 \,  \ee^{ 2 \sqrt{2 \alpha_i \gamma_i \log  {1\over r_{{\mathfrak a}_n}}}} \, \Big({1\over r_{{\mathfrak a}_n}} \Big)^{2(1+\varepsilon)\sqrt{\alpha_i \gamma_i}-\alpha_i }.
\end{eqnarray*}
\end{itemize}

We can then conclude by integrating over $u$.

\medskip

By construction, the measure ${\mathcal M}_\infty^{\alpha_1,\ldots,\alpha_p}$ satisfies that  for any nonnegative measurable functions $f_i$, $1\le i \le p$,   \begin{eqnarray}\label{spineintersection}
 && \otimes_{i=1}^p \e^{z_i, z'_i}_D \int_{D}  \prod_{i=1}^p f_i(x,  B^{(i)}) {\mathcal M}_\infty^{\alpha_1,\ldots,\alpha_p}(\d x)
\\
 &=&\int_{D}  \e_{\otimes_{i=1}^p \q^{z_i, z'_i, \alpha_i}_{x, D}} \Big(  \prod_{i=1}^p f_i(x,  B^{(i)})\Big)  \, \ee^{- C_D(x) \sum_{i=1}^p \alpha_i }\, \prod_{i=1}^p \xi_{D}(x,z_i,z'_i) \d x. \nonumber
\end{eqnarray}

 We refer to Le Gall (\cite{Legall90}, Chapter 9, Theorem 1) for a similar decomposition  involving the self-intersection local times. The measure ${\mathcal M}_\infty^{\alpha_1,\ldots,\alpha_p}$ was constructed by other means by Jego \cite{jego19}, see Proposition 1.2 there.

When $\alpha_1=...=\alpha_p=0$,  $\pi^{-p}  {\mathcal M}_\infty^{0, ..., 0}$ is nothing else than the image of the $p$-th  intersection local times measure. In fact,  let  $\gamma(\d s_1 ...\, \d s_p)$  be the  intersection local times of  $p$ independent excursions $B^{(1)}, ..., B^{(p)}$, which means that $\gamma$ is a measure supported by $\{(s_1, ...., s_p) \in \r_+^p: B^{(1)}_{s_1}=...=B^{(p)}_{s_p}\}$ and formally defined by   $$ \gamma(\d s_1 ...\, \d s_p):=  \delta_{\{0\}}(B^{(1)}_{s_1}-B^{(2)}_{s_2})  \cdots  \delta_{\{0\}}(B^{(p-1)}_{s_{p-1}}-B^{(p)}_{s_p}) \, \d s_1 \cdots \d s_p.$$

\noindent Then we may check as in the previous subsection that  $\otimes_{i=1}^p \p^{z_i, z'_i}_D$-a.s.,  $\pi^{-p}  {\mathcal M}_\infty^{0, ..., 0}$ is exactly the image of $\gamma$ by the application $ (s_1, ...., s_p) \to B^{(1)}_{s_1}$.  The details are omitted. 

\medskip 
\medskip
 
 \appendix
 
 \section{Appendix}
 
   \subsection{Justification of \eqref{conv-functionalzz'}} \label{A:s1}
 
 The existence of the limiting probability measure  $\p^{z,z'}_D$  follows from   the Kolmogorov extension theorem once we have proven that  \begin{equation}\label{conv-functionalF}
  \lim_{\varepsilon\to 0^+  }\e^{z + \varepsilon {\bf n}_z,z'}_D\left( F( B_{T_{D\cap {\cal C}(z,r)+t}}, \, t\ge 0) \right) \qquad \mbox{exists},
\end{equation}

  \noindent for any  $0< r < | z- z'|$ and for any bounded measurable function $F$ on ${\mathcal K}$.

Let $\nu_{D}(z)$ be the excursion measure in $D$ at $z$ defined by   \begin{equation} \label{mudy}
 \nu_{D}(z):= \lim_{\varepsilon \to 0^+ } \frac1{\varepsilon}  \p^{z + \varepsilon {\bf n}_z}_{D}, 
 \end{equation}

 \noindent in the sense that for  any     $r>0$ such that $D \cap {\cal C}(z, r) \neq \emptyset$, any ${\bf A} \in \sigma \{B_{T_{D \cap {\cal C}(z,r)}+t}, t\ge 0\}$, \begin{equation}\label{conv-functional-y}
 \nu_{D}(z)\left( {\bf A} \cap \{ T_{D \cap {\cal C}(z,r)} < T_{\partial D}\}\right)
 =   \lim_{\varepsilon \to 0^+ } \frac1{\varepsilon} \p_D^{z + \varepsilon {\bf n}_z}  \left({\bf A} \cap \{ T_{D \cap {\cal C}(z,r)}< T_{\partial D}\}\right).
 \end{equation}

 \noindent We refer to Burdzy (\cite{Burdzy}, Theorem   4.1 and p. 34--35)  for a justification  of \eqref{conv-functional-y} and the fact that  $\nu_D(z)(T_{D \cap {\cal C}(z,r)} < T_{\partial D}) < \infty$.  By linearity, \eqref{conv-functional-y} is true when replacing $1_{\bf A}$ by simple functions, then by any bounded measurable function (any measurable bounded function being limit of simple functions for the uniform norm).

 Now let $0< r < | z- z'|$.     Using the strong Markov property at time $T_{D\cap {\cal C}(z,r)}$ it is enough to show \eqref{conv-functionalF}  for $f( B_{T_{D\cap {\cal C}(z,r)}})$ instead of $F( B_{T_{D\cap {\cal C}(z,r)+t}}, t\ge 0)$, where $f$ denotes some bounded Borel function.    
 
  To this end,  we have by the definition of $\p^{y,z'}_{D}$ in \eqref{htransform2}  that for any $0< r'< | z- z'|- r$, 
\begin{align*}
&
	\lim_{\varepsilon\to 0^+  }\e^{z + \varepsilon {\bf n}_z,z'}_{D} \left( f( B_{T_{D\cap {\cal C}(z,r)}}) \right)
\\
&=
	\lim_{\varepsilon\to 0^+  }\e_D^{z + \varepsilon {\bf n}_z}\Big[ \frac{H_D(B_{T_{D\cap {\cal C}(z',r')}}, z')}{H_D(z + \varepsilon {\bf n}_z, z')} 		f( B_{T_{D\cap {\cal C}(z,r)}}) , T_{{\cal C}(z',r')}<T_{\partial D} \Big]
\\
&=
	\lim_{\varepsilon\to 0^+  } \frac1{\varepsilon H_D(z, z')} \e_D^{z + \varepsilon {\bf n}_z}\Big[   H_D(B_{T_{D\cap {\cal C}(z',r')}}, z')  f( B_{T_{D\cap 	{\cal C}(z,r)}}), T_{{\cal C}(z',r')}<T_{\partial D} \Big]
\\
&=
	\frac1{H_D(z, z')} \nu_D(z) \Big[   H_D(B_{T_{D\cap {\cal C}(z',r')}}, z')  f( B_{T_{D\cap {\cal C}(z,r)}}) , T_{{\cal C}(z',r')}<T_{\partial D}\Big], 
\end{align*}

 \noindent by \eqref{mudy}. This completes the proof of \eqref{conv-functionalF} and  justifies the definition of $\p^{z, z'}_D$  in \eqref{conv-functionalzz'}.

  \subsection{Justification of \eqref{conv-functional4}} \label{A:s2}

Let  $r>0$ be  small such that ${\cal B}(x, r) \subset D$ and let  $\Phi$ be a bounded    measurable function on ${\mathcal K}$.  By the strong Markov property at $T_{{\cal C}(x,r)}$, \begin{eqnarray*}
\e^{z,x}_D\left[ \Phi(B_{T_{{\cal C}(x,r)}+s}, s\ge0) ,   T_{{\cal C}(x,r)}< T_{x} \right]
=
\e^{z,x}_D\left[f(B_{T_{{\cal C}(x,r)}}), T_{{\cal C}(x,r)}< T_{x} \right],
\end{eqnarray*}

\noindent with $f(b):= \e^{b}_{x,D}\left[ \Phi(B_s, s\ge 0)\right]$.  By the $h$-transform, 
\begin{eqnarray*}
&&
	\e^{z,x}_D\left[f(B_{T_{{\cal C}(x,r)}}), T_{{\cal C}(x,r)}< T_{x} \right] \\
&=&
	\lim_{\varepsilon\to 0^+}  \e^{z,x}_D \left[f(B_{T_{{\cal C}(x,r)}}), T_{{\cal C}(x,r)}< T_{{\cal C}(x,\varepsilon)} \right]
\\
&=&
	\lim_{\varepsilon\to 0^+}  \e^{z}\left[ \frac{G_D(x, B_{T_{{\cal C}(x,r)}})}{G_D(x, z)} f(B_{T_{{\cal C}(x,r)}}), T_{{\cal C}(x,r)}< T_{{\cal C}(x,			\varepsilon)} \right]
\\
&=&
	 \e^{z}\left[ \frac{G_D(x, B_{T_{{\cal C}(x,r)}})}{G_D(x, z)} f(B_{T_{{\cal C}(x,r)}}) \right].
\end{eqnarray*}

\noindent It follows that 
\begin{eqnarray*}
&&
	\lim_{z \to x}  \left( \log  \frac1{|z-x|}\right) \,  \e^{z,x}_D\left[ \Phi(B_{T_{{\cal C}(x,r)}+\cdot}) ,   T_{{\cal C}(x,r)}< T_{x} \right]
\\
&=&
	\e^{z}\left[  G_D(x, B_{T_{{\cal C}(x,r)}}) f(B_{T_{{\cal C}(x,r)}}) \right].
\end{eqnarray*}

Therefore the proof of \eqref{conv-functional4} reduces to show that 
\begin{equation}\label{A-mudx1}
\nu_{D}(x, x)\left[ \Phi({\mathfrak e}_{T_{{\cal C}(x,r)}+\cdot}) ,   T_{{\cal C}(x,r)}< T_{x} \right]
= \e^{z}\left[  G_D(x, B_{T_{{\cal C}(x,r)}}) f(B_{T_{{\cal C}(x,r)}}) \right], 
\end{equation}

\noindent where as before  $f(b):= \e^{b}_{x,D}\left[ \Phi(B)\right]$ for any $b \in {\cal C}(x,r)$.  By the monotone class theorem, it is enough to consider $\Phi({\mathfrak e}_{T_{{\cal C}(x,r)}+\cdot})$ of form $\Phi_{r, \delta}({\mathfrak e}):= \Phi({\mathfrak e}_{T_{{\cal C}(x,r)}+s}, s \le T_{{\cal C}(x,\delta)}\circ \theta_{T_{{\cal C}(x,r)}} )$ where $0< \delta <r$ is arbitrary and $\theta$ denotes the usual shift operator.  Then 
\begin{eqnarray*}
 	f(b)
&=&
	\e^{b,x}_{D}\left[ \Phi(B_s, s\le T_{{\cal C}(x,\delta)})\right] \\
 &=&
 	\e^{b}\left[ \frac{G_D(x,  B_{T_{{\cal C}(x,\delta)}})}{G_D(x, b)} \Phi(B_s, s\le T_{{\cal C}(x,\delta)}), T_{{\cal C}(x,\delta)}< T_{\partial D}\right],
 \end{eqnarray*}
 
 \noindent where the second equality  is due to the $h$-transform. 
 By the definition of $\mu_D(x, x; t)$ and the strong Markov property at $T_{{\cal C}(x,\delta)}\circ \theta_{T_{{\cal C}(x,r)}}$, 
 \begin{align*}
 	&\mu_D(x, x; t) \left[ \Phi_{r, \delta}({\mathfrak e}), \, T_{{\cal C}(x,r)}< T_x\right]
 \\
 	&=
 	\lim_{\varepsilon\to0^+} \frac1{\pi \varepsilon^2} \, \e^x  \left[ \Phi_{r, \delta}(B) , \, T_{{\cal C}(x,\delta)}\circ \theta_{T_{{\cal C}(x,r)}} <  t \wedge 		T_{\partial D} ,  \p^{a} \big( |B_s - x | < \varepsilon, s < T_{\partial D}\big) \right]
 \\
 	&=  \e^x  \left[ \Phi_{r, \delta}(B) \,  p_D(s, x, a) , \, T_{{\cal C}(x,\delta)}\circ \theta_{T_{{\cal C}(x,r)}} <  t \wedge T_{\partial D} \right] , 
 \end{align*}

\noindent where  $s:=t-T_{{\cal C}(x,\delta)}\circ \theta_{T_{{\cal C}(x,r)}} , a:=B_{T_{{\cal C}(x,\delta)}\circ \theta_{T_{{\cal C}(x,r)}} }$ and the last equality follows from the definition of the transition probabilities $p_D$.  Taking the integral over $t$ gives that 
\begin{eqnarray*}
&& 
	\nu_D(x, x) \left[ \Phi_{r, \delta}({\mathfrak e}), \, T_{{\cal C}(x,r)}< T_x\right] \\
&=&
	\e^x\left[ \Phi_{r, \delta}(B) \, G_D( x , B_{T_{{\cal C}(x,\delta)}\circ \theta_{T_{{\cal C}(x,r)}} }), \, T_{{\cal C}(x,\delta)}\circ \theta_{T_{{\cal C}(x,r)}} <   	T_{\partial D} \right].
\end{eqnarray*} 

 By the strong Markov property at $T_{{\cal C}(x,r)}$ and  the $h$-transform,  the above expectation term $\e^x[\, \cdots]$  is equal to $$\e^x\left[ G_D(x, B_{T_{{\cal C}(x,r)}}) f(B_{T_{{\cal C}(x,r)}}), \,  T_{{\cal C}(x,r)} <   T_{\partial D}   \right], $$

\noindent yielding \eqref{A-mudx1} because $T_{{\cal C}(x,r)} <   T_{\partial D}$ holds with probability one. \hfill$\Box$

 \bigskip
{\noindent\bf Acknowledgements.}
We are grateful to Chris Burdzy for raising the question about existence of local times for thick points. The project was partly supported by ANR MALIN (ANR-16-CE93-0003); E.A.\   acknowledges supports from ANR GRAAL (ANR-14-CE25-0014) and  ANR Liouville (ANR-15-CE40-0013), Y.H.   acknowledges supports from ANR  SWiWS (ANR-17-CE40-0032).

\end{document}